\documentclass[12pt]{article}
\usepackage{latexsym}
\usepackage{amsbsy}
\usepackage{graphics}
\usepackage{tikz}
\usepackage{pgf}
\usepackage{pgflibraryarrows}

\usepackage{amssymb}

\PassOptionsToPackage{hyphens}{url}
\usepackage{hyperref}

\definecolor{greenish}{rgb}{0,0.5,0}
\definecolor{gold}{rgb}{0.85,.66,0}
\newtheorem{thm}{Theorem}
\newtheorem{cor}[thm]{Corollary}
\newtheorem{propn}[thm]{Proposition}
\newtheorem{lemma}[thm]{Lemma}

\newcommand{\pf}{\noindent\textit{Proof.\ \ }}
\newcommand{\eopf}{\hfill\hspace*{4pt}\hfill
\fbox{\rule[-1pt]{0pt}{4pt}\hspace*{2pt}}}
\newcommand{\caseSpace}{\rule[0pt]{0pt}{15pt}}

\newcommand{\dd}{\ensuremath{.\,.}}

\newcommand{\vr}[2]{
\left( \begin{array}{c}#1\\#2\end{array} \right)}
\newcommand{\textvr}[2]{
\hbox{\footnotesize$\vr{#1}{#2}$}}

\newcommand{\tri}{\hbox{\rm tri}}
\newcommand{\med}{\hbox{\rm med}}
\newcommand{\cfaces}{\hbox{\rm\scriptsize cf}}
\newcommand{\afaces}{\hbox{\rm\scriptsize af}}
\newcommand{\alt}{\hbox{\rm alt}}
\newcommand{\alta}{\alt_{\hbox{\scriptsize a}}}
\newcommand{\altc}{\alt_{\hbox{\scriptsize c}}}
\newcommand{\alti}{\alt_{\hbox{\scriptsize i}}}

\newcommand{\minor}[1]{\,\|_{{}_{#1}}}
\newcommand{\mino}[1]{[#1]}    %  8/4/12

\newlength{\hatsize}
\newlength{\hataltitude}

\title{Minors for alternating dimaps}
\author{
G. E. Farr\thanks{Part of the work of this paper was done
while the author was a Visiting Fellow (Combinatorics and Statistical Mechanics
programme) at the Isaac Newton
Institute for Mathematical Sciences, Cambridge, U.K.,
Jan.--Feb.\ 2008, and on sabbatical at the Department of Mathematics and Statistics,
University of Melbourne (Jan.--June\ 2011).  Parts have been presented at the Newton Institute
(Feb.\ 2008), Queen Mary,
University of London (Feb.\ 2011),
and the 36th Australasian Conference
on Combinatorial Mathematics and Combinatorial Computing (36 ACCMCC, Dec.\ 2012).  This research was supported in part by ARC Discovery Grant DP110100957.}
\thanks{Email: \href{mailto:Graham.Farr@monash.edu}{\texttt{Graham.Farr@monash.edu}}}   \\
Clayton School of Information Technology   \\
Monash University, Clayton, Victoria 3800   \\
Australia  \\
}
\date{30 November 2013}     %  started this file on 6 April 2012

\begin{document}

\maketitle

\begin{abstract}
We develop a theory of minors for alternating dimaps --- orientably
embedded digraphs where, at each vertex, the incident edges (taken in the order
given by the embedding) are directed alternately into, and out of, the vertex.
We show that they are related by the triality relation of Tutte.  They do
not commute in general, though do in many circumstances, and we characterise the situations
where they do.  The relationship with triality is reminiscent of similar relationships for binary
functions, due to the author, so we characterise those alternating dimaps which correspond
to binary functions.
We give a characterisation of alternating dimaps of at most a given genus, using
a finite set of excluded minors.  We also use the minor operations to define simple Tutte invariants
for alternating dimaps and characterise them.  We establish a connection with the Tutte
polynomial, and pose the problem of characterising universal
Tutte-like invariants for alternating dimaps based on these minor operations.
\end{abstract}


\section{Introduction}

The \textit{minor} relation is one of the most important order relations on graphs.
A graph $H$ is a \textit{minor} of a graph $G$ if it can be obtained from $G$ by some
sequence of deletions and contractions of edges.  Many important classes of graphs
can be characterised by the exclusion of some finite set of minors.
These include forests, series-parallel graphs \cite{dirac52,duffin65},
outerplanar graphs \cite{chartrand-harary67}, planar graphs \cite{kuratowski30,wagner37} ---
and, in fact, \textit{any} minor-closed class of graphs, by Robertson and Seymour's proof
of Wagner's conjecture \cite{robertson-seymour04}.
Minors also play a central role in enumerative graph theory:
the Tutte-Whitney polynomials, which contain information on a great variety of counting
problems on graphs or matroids, satisfy recurrence relations using deletion and contraction
(see, e.g., \cite{brylawski-oxley92,ellismonaghan-merino2011,farr07a,welsh93}).

The theory of minors derives much of its richness and beauty from the fact that the
deletion and contraction operations are dual (in the sense of planar graph duality or,
more generally, matroid duality \cite{oxley92}) and commute.

In this paper, we introduce and study a minor relation on alternating dimaps.
An \textit{alternating dimap} is a directed graph
without isolated vertices, 2-cell-embedded in a disjoint union
of orientable 2-manifolds,
where each vertex has even degree and, for each vertex $v$,
the edges incident with $v$ are
directed alternately into, and out of, $v$ (when considered in
the order in which they appear around $v$ in the embedding).
An alternating dimap may have loops and/or multiple edges,
but cannot have a bridge.  We allow the \textit{empty alternating dimap} with no
vertices, edges or faces.

For alternating dimaps, we have three minor operations, instead of two.
We show in \S\ref{sec:minors}
that they are related by a triality relation of Tutte \cite{tutte48}, in a manner
analogous to the duality between deletion and contraction.  The form of the relationship
is the same as that found by the author for some other combinatorial objects
(binary functions) on which minors
and triality can be defined \cite{farr2013}.

One property of ordinary minor operations
(and also of the minor operations in \cite{farr2013}) is that they commute.  We show
in \S\ref{sec:non-commut} that
minor operations on alternating dimaps do not commute in general, although they do in
most circumstances, and we determine exactly when they do.

Although duality, with associated minors, appears in many forms for many different kinds
of objects, there are far fewer settings with natural minor operations related by triality.
Two of these are binary functions \cite{farr2013} and alternating dimaps.
It is natural, then, to ask what relationship there may be between these settings.
Some alternating dimaps certainly cannot be represented by binary functions, or by
transition matroids, since
alternating dimap minor operations may not commute, unlike those in the
other settings.  In \S\ref{sec:bin-fns}
we determine those alternating dimaps that can be represented faithfully by binary functions.
(Other settings with natural minor operations and triality are multimatroids
(including isotropic systems) \cite{bouchet87,bouchet98} and the related
transition matroids \cite[pp.\ 8, 10]{traldi2013}.
Connections between them and alternating dimaps are a topic for future work, though
the settings seem to differ significantly since, in these other settings, the minor operations
commute.)

As seen in the first paragraph, two of the main themes of the classical theory of minors are
excluded minor characterisations and Tutte invariants.  The remainder of this paper takes
these themes up for alternating dimaps and their minor operations.
In \S\ref{sec:excl-minors-fixed-genus} we give an excluded
minor characterisation of alternating dimaps of at most a given genus, using a finite set
of excluded minors in every case.  In \S\ref{sec:tutte-invariants}
we define simple Tutte invariants for alternating dimaps,
and show that there are only a few of them and that they do not contain much information,
in contrast to the situation for graphs, matroids and binary functions.  We then define extended
Tutte invariants and raise the question of how many and varied they might be.  We show that
they are much richer than the simple Tutte invariants, as they include, in a sense, the Tutte
polynomial of a planar graph.

\S\S\ref{sec:non-commut}--\ref{sec:tutte-invariants} may be read independently of each other.

\subsection{History and related work}

Plane alternating dimaps were studied by Tutte \cite{tutte48,tutte75}.
He showed that they come in triples, with the three members of
a triple being derivable from a single larger structure, a \textit{bicubic map} (see below).
The relationship among the members of such a triple
is called \textit{triality} \cite{tutte48} or \textit{trinity} \cite{tutte75},
and each is \textit{trial} or \textit{trine} to the others.
This relationship extends ordinary duality.
Tutte's original motivation was to determine when equilateral triangles
can be tesselated by smaller equilaterial triangles of different sizes
or orientations.  He also proved his Tree Trinity Theorem, on spanning
arborescences of such maps.
In the second paper \cite{tutte75}
he noted the possibility of extending this theory to other surfaces.
Berman \cite{berman80} showed explicitly how to construct the trial
of an alternating dimap without reference to the bicubic
map from which three trial maps are derived, and gave alternative proofs
of some of Tutte's results.  Tutte reviewed some aspects of his theory in \cite{tutte99}.

It is interesting to note that this stream of research, first seen in Tutte's 1948 paper \cite{tutte48},
can be traced back to the same source that eventually
gave rise to Tutte's work on minor operations
and his eponymous polynomial.  Historically, the source of both streams
was the famous paper on ``squaring the square'' \cite{brooks-etal40}.
The 1948 paper extended the theory to ``triangulating the triangle'' (where all
triangles are equilateral) and introduced triality, among other things.  However, this stream
has not previously seen the development of minor
operations or Tutte-like invariants for alternating dimaps.

Tutte showed that a ``triangulated triangle'' gives rise to a plane bipartite
cubic graph (a \textit{plane bicubic map}), which in turn has, as its dual,
an Eulerian plane triangulation.  The triple of trial alternating dimaps is derived from
this bicubic map.
Although bicubic maps are plane in Tutte's work, they may more generally be taken to
be embedded in some orientable surface, so each has a genus.

A separate stream of research concerns \textit{latin bitrades}, which are pairs of partial
latin squares of the same shape and with the same symbol set in each row and column.
These may also be given a genus.
Cavenagh and Lison\v{e}k \cite{cavenagh-lisonek08} established a correspondence between
spherical latin bitrades and 3-connected planar Eulerian triangulations
(dually, 3-connected plane bicubic maps), while the relationship between spherical latin bitrades
and triangulated triangles is described in \cite{drapal91,drapal09}
(see also \cite{cavenagh-wanless09}).  Batagelj \cite{batagelj89}
introduced two operations on plane Eulerian triangulations by means of which larger such maps
can be generated from smaller ones.  These operations have been used and extended
in several papers (via the aforementioned correspondences)
to generate latin bitrades \cite{drapal09,grubman-wanlessXX}.
The inverse of one of these operations,
translated to alternating dimaps, corresponds to (technically, restricted versions of) our minor
operations.\footnote{I thank Tony Grubman and Ian Wanless for pointing out this link.}
This inverse also appears as a reduction on plane bicubic maps due to Jaeger \cite{jaeger92}.

\subsection{Definitions and notation}

If $G$ is an alternating dimap then $kG$ is the disjoint union
of $k$ copies of $G$.  Each of these copies is regarded as being
embedded in a different surface, with all these $k$ surfaces
being disjoint from each other.

An edge $e$ from $u$ to $v$ is sometimes written $e(u,v)$.

Let $G$ is an alternating dimap, viewed topologically as embedded
in a surface.  Let $C\subseteq E(G)$ be a circuit of $G$, and let $S$ be the connected
surface in which the component of $G$ containing $C$ is embedded.
Then the \textit{sides} of $C$ are
the components of $S-C$.
A side is \textit{planar} if it is homeomorphic to the open unit disc.   % **** Do I still need this?

The \textit{genus} $\gamma(G)$ of an alternating dimap $G$
is given by
\[
|V(G)| - |E(G)| + |F(G)| = 2(k(G)-\gamma(G)),
\]
where $F(G)$ is the set of faces of $G$ and $k(G)$ is the number
of components of $G$.

The edges around a face all go in the
same direction, and we say the face is \textit{clockwise} or
\textit{anticlockwise} according to the direction of the edges
around it.  We identify a clockwise (respectively, anticlockwise) face
with its cyclic sequence
of edges, and call it a \textit{c-face} (resp., \textit{a-face})
for short.  Observe that the (edge sets of the)
c-faces partition $E(G)$, as do the a-faces.
So every edge $e$ belongs to one c-face, denoted by $C(e)=C_G(e)$,
and one a-face, denoted by $A(e)=A_G(e)$.
If two faces share a common edge, then one of the faces
is clockwise and the other is anticlockwise.
The \textit{left successor} (respectively, \textit{right successor})
of an edge $e$ is the next edge
along from $e$, going around $A_G(e)$ (resp., $C_G(e)$) in the direction given by $e$.
(This direction is anticlockwise for the left successor, and clockwise for the right successor.)
Often, c-faces and a-faces are simple cycles, but this is not always
the case.  If $v$ is a cutvertex of $G$, then one face incident with $v$
consists of two or more
edge-disjoint cycles.

The numbers of clockwise and anticlockwise faces of $G$ are denoted by $\hbox{cf}(G)$ and
$\hbox{af}(G)$, respectively.

An \textit{in-star} is the set of all edges directed into some vertex.
So in-stars are in one-to-one correspondence with vertices.
The in-star of edges directed into vertex $v$ is denoted by $I(v)=I_G(v)$.
Observe that the in-stars partition $E(G)$, so every edge $e$ also belongs
to one in-star, denoted by $I(e)=I_G(e)$ (overloading notation slightly).

If an alternating dimap is disconnected, then we treat its components
as being embedded in separate, disjoint surfaces.

Throughout, we set $\omega=\exp(2\pi i/3)$.

Alternating dimaps extend ordinary embedded graphs, in that replacing each edge of an
embedded graph by a pair of directed edges, forming a clockwise face of size 2, gives
an alternating dimap \cite{tutte75}.

An alternating dimap $G$ defines three permutations
$\sigma_{G,1},\sigma_{G,\omega},\sigma_{G,\omega^2}\colon E(G)\rightarrow E(G)$
(abbreviated $\sigma_{1},\sigma_{\omega},\sigma_{\omega^2}$
where $G$ is clear from the context), as follows.
For each $e\in E(G)$, its image under
$\sigma_{G,1}$, $\sigma_{G,\omega}$ and $\sigma_{G,\omega^2}$ is
the next edge in clockwise order around $I_G(e)$, $A_G(e)$ and $C_G(e)$,
respectively.  So the left successor of $e$ is $\sigma_{G,\omega}^{-1}(e)$,
while the right successor of $e$ is $\sigma_{G,\omega^2}(e)$.
Note that, in going around an in-star in clockwise order,
we skip outgoing edges at the vertex as these do not belong to the in-star.

Any two of these permutations determine the other.  This follows from
the relation $\sigma_{1}\sigma_{\omega}\sigma_{\omega^2}=\hbox{id}_{E(G)}$,
the identity permutation on $E(G)$.
% *** Associates from right, i.e., as if using  \circ.  Should I put  \circ  in?

Let $S$ be the set of all
triples $(\sigma_0,\sigma_1,\sigma_2)$ of permutations, all acting
on the same set $E$, such that
$\sigma_0\sigma_1\sigma_2=\hbox{id}_E$.  This is called a \textit{3-constellation}
or a \textit{hypermap} \cite{lando-zvonkin04}.  When one of the permutations is an
involution, we may take these involutions to correspond to undirected edges
(using the aforementioned
representation of embedded graphs by alternating dimaps), and we have a
standard combinatorial representation of an orientably embeddeded graph (see, e.g.,
\cite[\S2.2]{bonnington-little95}).  In the general case, we have an equivalence with
alternating dimaps which seems to be well known (see, e.g., \cite{cori-penaud80})
although I have not seen it stated explicitly.

\begin{propn}
The map $\{\hbox{alternating dimaps}\}\rightarrow S$ given by
$G\mapsto (\sigma_{G,1},\sigma_{G,\omega},\sigma_{G,\omega^2})$ is a bijection.
\eopf
\end{propn}

%  \pf
%  To show that any $(\sigma_0,\sigma_1,\sigma_2)\in S$ arises from an
%  alternating dimap in this way, define $G$ as follows.  Its vertices
%  represent the cycles of $\sigma_0$.  For any $e\in E$, write
%  $A_e$ and $C_e$ for the
%  cycle of $\sigma_1$ and $\sigma_2$, respectively, which contains $e$.
%  There is an edge in $G$ from vertex $u$ to
%  vertex $v$ (representing cycles which we shall call $I_u$ and $I_v$,
%  respectively, in $\sigma_0$) if and only if there exists $e\in I_v$
%  such that $\sigma_1(e)$ and $\sigma_2^{-1}(e)=\sigma_0\sigma_1(e)$ both
%  belong to $I_u$.  The in-stars, a-faces and c-faces of $G$ are the
%  cycles of $\sigma_0$, $\sigma_1$ and $\sigma_2$, respectively.  The
%  rest of the proof is left as an exercise.
%  \eopf   \\

%  If $F$ is a c-face or an a-face containing vertices $v$ and $w$, then
%  $F[v,w]$ is the path from $v$ to $w$ around $F$.  If $I$ is an in-star
%  containing edges $e$ and $f$, then $I[e,f]$ is the subgraph of $I$
%  consisting of the edges from $e$ to $f$ in clockwise order around $I$,
%  together with their endpoints.

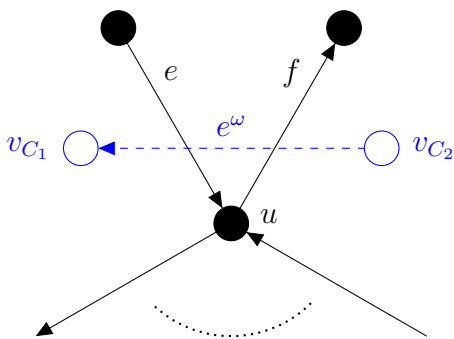
\begin{figure}[ht]
\begin{center}
\begin{tikzpicture}
\begin{scope}[>=triangle 45]
%  before taking trial ...
\coordinate (u) at (3,2);
\coordinate (a) at (1.5,4.6);
\coordinate (b) at (4.5,4.6);
\coordinate (c) at (0.4,0.5);
\coordinate (d) at (5.6,0.5);
\draw[fill] (u) circle (0.23cm);
\draw[fill] (a) circle (0.23cm);
\draw[fill] (b) circle (0.23cm);
\draw[->>] (a) -- (u);
\draw[->>] (u) -- (b);
\draw[->] (u) -- (c);
\draw[->>] (d) -- (u);
\draw[dotted,thick] (4.06,0.94) arc (315:225:1.5);
\draw (3.5,2.1) node {$u$};
\draw (2.2,4) node {$e$};
\draw (3.8,4) node {$f$};
%  trial
\coordinate (vC1) at (1,3);
\coordinate (vC2) at (5,3);
\draw[blue] (0.3,3) node {$v_{C_1}$};
\draw[blue] (5.7,3) node {$v_{C_2}$};
\draw[blue] (3,3.3) node {$e^{\omega}$};
\draw[dashed,blue,->>] (vC2) -- (vC1);
\draw[draw=blue,fill=white,opaque] (vC1) circle (0.23cm);
\draw[draw=blue,fill=white,opaque] (vC2) circle (0.23cm);
\end{scope}
\end{tikzpicture}
\end{center}
\caption{Construction of trial map.}
\label{fig:trial}
\end{figure}

If $G$ is an alternating dimap, then the \textit{trial}
%  (pronounced with a short `i')
$G^{\omega}$ of $G$
%, and the\textit{trial map} $t:E(G)\rightarrow E(G^{\omega})$,
is defined as follows.  Its vertices represent the c-faces of $G$.
We denote the vertex of $G^{\omega}$ representing c-face $C$ by $v_C$.
(Think of $v_C$ being placed inside $C$ in the embedding.)
Edges of $G^{\omega}$ are constructed as follows.
Suppose two c-faces $C_1$ and $C_2$ of $G$
share a vertex $v$, and that there is an
a-face $A$ containing edges $e$ and $f$ going into and out of $v$
respectively, with $e$ and $f$ also belonging to $C_1$ and $C_2$
respectively.  (See Figure~\ref{fig:trial}.  We do not require $e$ and $f$
to be distinct, or $C_1$ and $C_2$ to be distinct.)
Then we put an edge $e^{\omega}$ from $v_{C_2}$ to $v_{C_1}$ in $G^{\omega}$.
%  In the above line:  was  1-to-2, but fixed.  Error noticed by L Traldi.
(Think of $e^{\omega}$ as being drawn by a curve from $v_{C_2}$,
%   the vertex of $G^{\omega}$
inside $C_2$, to its destination $v_{C_1}$ inside $C_1$,
in such a way that
it crosses $f$ in its ``first half'' (i.e., closer to its start than its end)
and crosses $e$ in its ``last half''.)
These edges of $G^{\omega}$ are ordered around $C_1$ according to the order
of the edges $f$ around $C_1$.  Similarly, they are ordered around $C_2$
according to the order of the edges $e$ around $C_2$.
It is routine to show that
the map $\bullet^{\omega}:E(G)\rightarrow E(G^{\omega})$, $e\mapsto e^{\omega}$
is a bijection, and that the
c-faces, a-faces and in-stars of $G^{\omega}$ are the a-faces, in-stars
and c-faces, respectively, of $G$.  We can also express this relationship
in the language of the
permutation triples.

\begin{propn}
\label{propn:trial-sigma}
\begin{eqnarray}
\sigma_{G,1}(e)^{\omega}  & = & \sigma_{G^{\omega},\omega}(e^{\omega})
\label{eq:trial-sigma-1}   \\
\sigma_{G,\omega}(e)^{\omega}  & = & \sigma_{G^{\omega},\omega^2}(e^{\omega})
\label{eq:trial-sigma-omega}   \\
\sigma_{G,\omega^2}(e)^{\omega}  & = & \sigma_{G^{\omega},1}(e^{\omega})
\label{eq:trial-sigma-omega2}
\end{eqnarray}
\eopf
\end{propn}

If $G$ is represented by $(\sigma_{G,1},\sigma_{G,\omega},\sigma_{G,\omega^2})$,
then its trial $G^{\omega}$ is represented by
$(\sigma_{G,\omega},\sigma_{G,\omega^2},\sigma_{G,1})$.

It is clear that Proposition \ref{propn:trial-sigma} still holds if
$\sigma$ is replaced by $\sigma^{-1}$ throughout.

The trial operation on a component of $G$
is independent of the other components.   \\

%  When discussing an alternating dimap $G$ and its trial $G^{\omega}$ and
%  double trial  *****

We write $G^{\omega^2}$ for $(G^{\omega})^{\omega}$.  From the way triality
changes c-faces to in-stars to a-faces to c-faces, we find that
$G^{\omega^3}=(G^{\omega^2})^{\omega}=G^1=G$.

%  ***** remarks on some things known about triality/trinity *****

%  **** terms with their usual undirected meaning: degree, cutset, loop? ****

\subsection{Loops and semiloops}

A \textit{(standard) loop} is just an edge of $G$ which is a loop
in the undirected version of $G$.  If it is a separating circuit of the embedding,
then it divides its component of the
embedding surface into two
sides, its \textit{clockwise side} and its \textit{anticlockwise side}.  If it is non-separating, then
it has just one side, which we take to be both its clockwise and anticlockwise side.

A \textit{$1$-loop} is an edge whose head has degree 2.
This is an edge whose left and right successors are identical, and in such cases
we can refer unambiguously to its \textit{successor}.
An \textit{$\omega$-loop} is an edge forming a single-edge a-face.
An \textit{$\omega^2$-loop} is an edge forming a single-edge c-face.
For any alternating dimap $G$,
an edge $e$ is a 1-loop in $G$ if and only if $e^{\omega}$ is an $\omega$-loop in
$G^{\omega}$, which in turn holds if and only if $e^{\omega^2}$ is an $\omega^2$-loop in
$G^{\omega^2}$.

An \textit{triloop} is an edge which is a $\mu$-loop for some
$\mu\in\{1,\omega,\omega^2\}$.

An \textit{ultraloop} is a triloop which (together with its vertex)
constitutes a component of the graph.  It has two faces, a c-face
and an a-face, and its vertex has degree 2, so it is simultaneously
a 1-loop, an $\omega$-loop and an $\omega^2$-loop.  In fact, if an edge
is a $\mu$-loop for $\mu$ equal to any two of $\{1,\omega,\omega^2\}$,
then it is an ultraloop.

Note that $\omega$- and $\omega^2$-loops are also standard loops,
but the converse is not necessarily true, as we will see when
considering 1-semiloops shortly.
Also, a 1-loop is typically not a standard loop,
since its vertices do not have to coincide.
A 1-loop is only a standard loop if it is an ultraloop.

A $\mu$-loop is a \textit{proper $\mu$-loop} if it is not also an ultraloop.
In such a case, it is not a $\nu$-loop for any $\nu\in\{1,\omega,\omega^2\}\setminus\{\mu\}$.

A \textit{1-semiloop} is just a standard loop.
An \textit{$\omega$-semiloop} is an edge $e$ such that $e$ and its right successor
$\sigma_{\omega^2}(e)$ are equal, or (if distinct) either they form a cutset of
$G$ or deleting both increases its genus.  (After deletion, we no longer have an alternating dimap
in general; we are really referring to the underlying undirected embedded graph here, rather than
$G$ itself.)
This latter condition, on cutset/genus, may be written:
$k(G\setminus\{e,\sigma_{\omega^2}(e)\}) - \gamma(G\setminus\{e,\sigma_{\omega^2}(e)\}) >
k(G) - \gamma(G)$.  Similarly,
an \textit{$\omega^2$-semiloop} is an edge $e$ such that $e$ and
its left successor
$\sigma_{\omega}^{-1}(e)$ together have this same property.

For each $\mu\in\{1,\omega,\omega^2\}$, a \textit{proper $\mu$-semiloop}
is a $\mu$-semiloop that is not a triloop.  A proper 1-semiloop $e$ either gives a non-contractible
closed curve in the embedding, or each of its two sides contains an edge other than $e$ from
the same component as $e$.

We make some basic observations about semiloops.

An edge is a $\mu$-semiloop in $G$ if and only if it is a $\mu\omega$-semiloop in
$G^{\omega}$.

An edge is both a $\mu_1$-semiloop and a $\mu_2$-semiloop if and only if it is
a $(\mu_1\mu_2)^{-1}$-loop.

\section{Minors}
\label{sec:minors}

%  \vspace{0.2cm}
%  \noindent\textbf{Definitions.}
Let $G$ be an alternating dimap and $e=e(u,v)\in E(G)$.

If $e$ is not a loop, then the graph $G\mino{1}e$
is formed by deleting the edge $e$
and identifying its endpoints, while preserving the order of the
edges and faces around vertices.  If $e$ is an $\omega$-loop or
an $\omega^2$-loop, then $G\mino{1}e$ is formed just by deleting $e$.
If $e$ is a 1-semiloop, then $G\mino{1}e$ is formed as follows.
Let the edges incident with $v$, in cyclic order around $v$ starting with $e$ directed into $v$, be
$e, a_1,b_1,\ldots,a_k, b_k,e, c_1,d_1,\ldots,c_l,d_l$.  Here, each $a_i$ and each $d_i$ is
directed out of $v$, while each $b_i$ and $c_i$ is directed into $v$.  We replace $v$ by two
new vertices, $v_1$ and $v_2$, and reconnect the edges $a_i,b_i,c_i,d_i$ as follows.
The tail of each $a_i$ and the head of each $b_i$ becomes $v_1$ instead of $v$, while
the head of each $c_i$ and the tail of each $d_i$ becomes $v_2$ instead of $v$.
The edge $e$ is deleted.  The cyclic orderings of edges around $v_1$ and $v_2$ are those
induced by the ordering around $v$.  Everything else is unchanged.
Observe that if $e$ is a separating 1-semiloop then, in effect,
the shrinking of the loop severs the graph in its clockwise side from
that in its anticlockwise side.  If $e$ is not separating then, in effect, shrinking the loop cuts one
of the handles of the surface component in which it is embedded, so reducing the genus by 1.
%  In any case, $k(G\mino{1}e)-\gamma(G\mino{1}e)=k(G)-\gamma(G)$.

This operation is called \textit{1-reduction} or
\textit{contraction}.  See Figure~\ref{fig:minor-opns}(a),(b).  It just adapts the usual contraction
operation for surface minors to the alternating dimap context.

\begin{figure}[htp]
%   Before taking minor
\begin{tikzpicture}
\begin{scope}[>=triangle 45]
%  vertex  u
\coordinate (u) at (3,6) {};
\draw[fill] (u) circle (0.23cm);
\draw (2.7,5.6) node {$u$};
%  vertex  v
\coordinate (v) at (3,3) {};
\draw[fill] (v) circle (0.23cm);
\draw (2.6,3.3) node {$v$};
%  edge  e(u,v)
\draw[->>] (u) -- (v);
\draw (3.3,4.5) node {$e$};
%  vertices  w1, w2, x
\coordinate (w1) at (0.12,2.16) {};
\draw (0,2.7) node {$w_1$};
\coordinate (w2) at (5.88,2.16) {};
\draw (6,2.7) node {$w_2$};
\coordinate (x) at (0.6,1.2) {};
\draw[fill] (w1) circle (0.23cm);
\draw[fill] (w2) circle (0.23cm);
\draw[fill] (x) circle (0.23cm);
%  other edges at  v
\draw[->>] (v) -- (w1);
\draw[->>] (v) -- (w2);
\coordinate (y) at (2.4,0) {};
\coordinate (z) at (4.6,0) {};
\draw[->>] (x) -- (v);
\draw[->] (v) -- (y);
\draw[->>] (z) -- (v);
\draw[dotted,thick] (3.68,1.12) arc (290:265:2);
%  other edges at  u
\coordinate (a) at (0.12,6.84) {};
\coordinate (b) at (3.8,8.5) {};
\coordinate (c) at (5.88,6.84) {};
\coordinate (d1) at (2.4,6.6) {};
\coordinate (d2) at (2.8,6.8) {};
\draw[->>] (a) -- (u);
\draw[->] (u) -- (b);
\draw[->>] (c) -- (u);
\draw (u) -- (d1);
\draw[->>] (d2) -- (u);
\draw[dotted,thick] (3,8) arc (90:160:1.5);
%  trial
\coordinate (om1) at (1,4.5);
\coordinate (om2) at (4.7,7.5);
\draw[dashed,blue] (om2) -- (3.9,7.2);
\draw[dashed,blue,->>] (4.0,7.8) -- (om2);
\draw[dashed,blue] (om2) -- (5.3,7.0);
\coordinate (om3) at (4.7,1.5);
\draw[dashed,blue,->>] (5.1,2.2) -- (om3);
\draw[dashed,blue] (om3) -- (4.3,1.1);
\draw[dashed,blue,->>] (4.0,2.0) -- (om3);
\coordinate (om4) at (1.8,1.2);
\draw[dashed,blue] (om4) -- (2.4,1.6);
\draw[dashed,blue,->>] (2.3,0.6) -- (om4);
\draw[dashed,blue] (om4) -- (1.3,1.5);
\draw[dashed,blue,->>] (om1) .. controls (3.9,5.1) .. (om2);
\draw[dashed,blue,->>] (om3) .. controls (3.9,3.9) .. (om1);
\draw[blue] (4.3,3.6) node {$e^{\omega}$};
\draw[dashed,blue,->>] (1.4,5.4) -- (om1);
\draw[dashed,blue] (0.7,5.5) -- (om1);
\draw[dashed,blue,->>] (om1) -- (om4);
\draw[dashed,blue,->>] (0.7,3.5) -- (om1);
\draw[dotted,blue,thick] (0.5,5) arc (135:225:0.7071);
\draw[draw=blue,fill=white,opaque] (om1) circle (0.23cm);
\draw[draw=blue,fill=white,opaque] (om2) circle (0.23cm);
\draw[draw=blue,fill=white,opaque] (om3) circle (0.23cm);
\draw[draw=blue,fill=white,opaque] (om4) circle (0.23cm);
\end{scope}
\end{tikzpicture}
~~~~~~~~~~~~~~
%   After taking 1-minor
\begin{tikzpicture}
\begin{scope}[>=triangle 45]
%  one merged vertex,  u = v
\coordinate (uv) at (3,4.5) {};
\draw[fill] (uv) circle (0.23cm);
\draw (2.2,4.5) node {$u=v$};
%  vertices  w1, w2, x
\coordinate (w1) at (0.12,2.16) {};
\draw (0,2.7) node {$w_1$};
\coordinate (w2) at (5.88,2.16) {};
\draw (6,2.7) node {$w_2$};
\coordinate (x) at (0.6,1.2) {};
\draw[fill] (w1) circle (0.23cm);
\draw[fill] (w2) circle (0.23cm);
\draw[fill] (x) circle (0.23cm);
%  other edges that were at  v
\draw[->>] (uv) .. controls (1.7,2.6) .. (w1);
\draw[->>] (uv) .. controls (4.3,2.6) .. (w2);
\coordinate (y) at (2.4,0) {};
\coordinate (z) at (4.6,0) {};
\draw[->>] (x) .. controls (2.2,2.4) .. (uv);
\draw[->] (uv) .. controls (2.7,1.5) .. (y);
\draw[->>] (z) .. controls (3.6,1.75) .. (uv);
%  \draw[->>] (3.8,0.5) -- (3,4.5);
%  \draw (3,4.5) -- (2.8,2.2);
%  \draw[->>] (2.4,2.4) -- (3,4.5);
\draw[dotted,thick] (3.68,1.12) arc (290:265:2);
%  other edges that were at  u
\coordinate (a) at (0.12,6.84) {};
\coordinate (b) at (3.8,8.5) {};
\coordinate (c) at (5.88,6.84) {};
\coordinate (d1) at (2.4,6.6) {};
\coordinate (d2) at (2.8,6.8) {};
\draw[->>] (a) .. controls (1.7,6.4) .. (uv);
\draw[->] (uv) -- (b);
\draw[->>] (c) .. controls (4.3,6.4) .. (uv);
\draw (uv) -- (d1);
\draw[->>] (d2) -- (uv);
\draw[dotted,thick] (3,8) arc (90:160:1.5);
%  trial
\coordinate (om1) at (1,4.5);
\coordinate (om2) at (4.7,7.5);
\draw[dashed,blue] (om2) -- (3.9,7.2);
\draw[dashed,blue,->>] (4.0,7.8) -- (om2);
\draw[dashed,blue] (om2) -- (5.3,7.0);
\coordinate (om3) at (4.7,1.5);
\draw[dashed,blue,->>] (5.1,2.2) -- (om3);
\draw[dashed,blue] (om3) -- (4.3,1.1);
\draw[dashed,blue,->>] (4.0,2.0) -- (om3);
\coordinate (om4) at (1.8,1.2);
\draw[dashed,blue] (om4) -- (2.4,1.6);
\draw[dashed,blue,->>] (2.3,0.6) -- (om4);
\draw[dashed,blue] (om4) -- (1.3,1.5);
\draw[dashed,blue,->>] (om3) -- (om2);
\draw[dashed,blue,->>] (1.4,5.4) -- (om1);
\draw[dashed,blue] (0.7,5.5) -- (om1);
\draw[dashed,blue,->>] (om1) -- (om4);
\draw[dashed,blue,->>] (0.7,3.5) -- (om1);
\draw[dotted,blue,thick] (0.5,5) arc (135:225:0.7071);
\draw[draw=blue,fill=white,opaque] (om1) circle (0.23cm);
\draw[draw=blue,fill=white,opaque] (om2) circle (0.23cm);
\draw[draw=blue,fill=white,opaque] (om3) circle (0.23cm);
\draw[draw=blue,fill=white,opaque] (om4) circle (0.23cm);
\end{scope}
\end{tikzpicture}
\hspace*{0.1cm}   \\   %  hack to get vertical gap
\hspace*{2cm} (a) $G$ and \textcolor{blue}{$G^{\omega}$}
\hspace*{4cm}(b)
$G\mino{1}e$ and \textcolor{blue}{$(G\mino{1}e)^{\omega}=G^{\omega}\mino{\omega^2}e^{\omega}$}  \\
\hspace*{0.1cm}   \\   \\   %  hack to get vertical gap
%   After taking omega-minor
\begin{tikzpicture}
\begin{scope}[>=triangle 45]
%  vertex  u
\coordinate (u) at (3,6) {};
\draw[fill] (u) circle (0.23cm);
\draw (2.7,5.6) node {$u$};
%  vertex  v
\coordinate (v) at (3,3) {};
\draw[fill] (v) circle (0.23cm);
\draw (2.6,3.3) node {$v$};
%  vertices  w1, w2, x
\coordinate (w1) at (0.12,2.16) {};
\draw (0,2.7) node {$w_1$};
\coordinate (w2) at (5.88,2.16) {};
\draw (6,2.7) node {$w_2$};
\coordinate (x) at (0.6,1.2) {};
\draw[fill] (w1) circle (0.23cm);
\draw[fill] (w2) circle (0.23cm);
\draw[fill] (x) circle (0.23cm);
%  other edges at  v
\draw[->>] (v) -- (w1);
\coordinate (y) at (2.4,0) {};
\coordinate (z) at (4.6,0) {};
\draw[->>] (x) -- (v);
\draw[->] (v) -- (y);
\draw[->>] (z) -- (v);
\draw[dotted,thick] (3.68,1.12) arc (290:265:2);
%  other edges at  u
\coordinate (a) at (0.12,6.84) {};
\coordinate (b) at (3.8,8.5) {};
\coordinate (c) at (5.88,6.84) {};
\coordinate (d1) at (2.4,6.6) {};
\coordinate (d2) at (2.8,6.8) {};
\draw[->>] (u) -- (w2);
\draw[->>] (a) -- (u);
\draw[->] (u) -- (b);
\draw[->>] (c) -- (u);
\draw (u) -- (d1);
\draw[->>] (d2) -- (u);
\draw[dotted,thick] (3,8) arc (90:160:1.5);
%  trial
\coordinate (om13) at (2.5,4.5);
\draw[dashed,blue,->>] (1.5,3.5) -- (om13);
\draw[dashed,blue] (1.5,5.5) -- (om13);
\draw[dashed,blue,->>] (2,5.5) -- (om13);
\draw[dotted,blue,thick] (1.634,5) arc (150:210:1);
\coordinate (om2) at (4.7,7.5);
\draw[dashed,blue] (om2) -- (3.9,7.2);
\draw[dashed,blue,->>] (4.0,7.8) -- (om2);
\draw[dashed,blue] (om2) -- (5.3,7.0);
%  \coordinate (om3) at (4.7,1.5);
\draw[dashed,blue,->>] (3.8,4.4) -- (om13);
\draw[dashed,blue] (om13) .. controls (4,4) and (5,2) .. (4.3,1.1);
\draw[dashed,blue,->>] (4.0,2.0) .. controls (3.5,3.5) .. (om13);
\coordinate (om4) at (1.8,1.2);
\draw[dashed,blue] (om4) -- (2.4,1.6);
\draw[dashed,blue,->>] (2.3,0.6) -- (om4);
\draw[dashed,blue] (om4) -- (1.3,1.5);
\draw[dashed,blue,->>] (om13) .. controls (3.9,5.1) .. (om2);
\draw[dashed,blue,->>] (om13) -- (om4);
\draw[draw=blue,fill=white,opaque] (om13) circle (0.23cm);
\draw[draw=blue,fill=white,opaque] (om2) circle (0.23cm);
\draw[draw=blue,fill=white,opaque] (om4) circle (0.23cm);
\end{scope}
\end{tikzpicture}
~~~~~~~~~~~~~~
%   After taking omega^2-minor
\begin{tikzpicture}
\begin{scope}[>=triangle 45]
%  vertex  u
\coordinate (u) at (3,6) {};
\draw[fill] (u) circle (0.23cm);
\draw (3.1,5.5) node {$u$};
%  vertex  v
\coordinate (v) at (3,3) {};
\draw[fill] (v) circle (0.23cm);
\draw (3.3,3.4) node {$v$};
%  vertices  w1, w2, x
\coordinate (w1) at (0.12,2.16) {};
\draw (0,2.7) node {$w_1$};
\coordinate (w2) at (5.88,2.16) {};
\draw (6,2.7) node {$w_2$};
\coordinate (x) at (0.6,1.2) {};
\draw[fill] (w1) circle (0.23cm);
\draw[fill] (w2) circle (0.23cm);
\draw[fill] (x) circle (0.23cm);
%  other edges at  v
\draw[->>] (v) -- (w2);
\coordinate (y) at (2.4,0) {};
\coordinate (z) at (4.6,0) {};
\draw[->>] (x) -- (v);
\draw[->] (v) -- (y);
\draw[->>] (z) -- (v);
\draw[dotted,thick] (3.68,1.12) arc (290:265:2);
%  other edges at  u
\coordinate (a) at (0.12,6.84) {};
\coordinate (b) at (3.8,8.5) {};
\coordinate (c) at (5.88,6.84) {};
\coordinate (d1) at (2.4,6.6) {};
\coordinate (d2) at (2.8,6.8) {};
\draw[->>] (u) -- (w1);
\draw[->>] (a) -- (u);
\draw[->] (u) -- (b);
\draw[->>] (c) -- (u);
\draw (u) -- (d1);
\draw[->>] (d2) -- (u);
\draw[dotted,thick] (3,8) arc (90:160:1.5);
%  trial
\coordinate (om1) at (0.8,5);
\draw[dashed,blue,->>] (1.2,5.9) -- (om1);
\draw[dashed,blue] (0.5,6) -- (om1);
\draw[dashed,blue,->>] (0.5,4) -- (om1);
\draw[dotted,blue,thick] (0.3,5.5) arc (135:225:0.7071);
\coordinate (om2) at (4.7,7.5);
\draw[dashed,blue] (om2) -- (3.9,7.2);
\draw[dashed,blue,->>] (4.0,7.8) -- (om2);
\draw[dashed,blue] (om2) -- (5.3,7.0);
\coordinate (om3) at (4.7,1.5);
\draw[dashed,blue,->>] (5.1,2.2) -- (om3);
\draw[dashed,blue] (om3) -- (4.3,1.1);
\draw[dashed,blue,->>] (4.0,2.0) -- (om3);
\coordinate (om4) at (1.8,1.2);
\draw[dashed,blue] (om4) -- (2.4,1.6);
\draw[dashed,blue,->>] (2.3,0.6) -- (om4);
\draw[dashed,blue] (om4) -- (1.3,1.5);
\draw[dashed,blue,->>] (om1) .. controls (4,4.7) .. (om2);
\draw[dashed,blue,->>] (om3) .. controls (4.7,5) and (1.8,5).. (om4);
\draw[draw=blue,fill=white,opaque] (om1) circle (0.23cm);
\draw[draw=blue,fill=white,opaque] (om2) circle (0.23cm);
\draw[draw=blue,fill=white,opaque] (om3) circle (0.23cm);
\draw[draw=blue,fill=white,opaque] (om4) circle (0.23cm);
\end{scope}
\end{tikzpicture}
\hspace*{0.1cm}   \\   %  hack to get vertical gap
\hspace*{0.2cm} (c)
$G\mino{\omega}e$ and \textcolor{blue}{$(G\mino{\omega}e)^{\omega}=G^{\omega}\mino{1}e^{\omega}$}
\hspace*{2cm}(d)
$G\mino{\omega^2}e$ and
\textcolor{blue}{$(G\mino{\omega^2}e)^{\omega}=G^{\omega}\mino{\omega}e^{\omega}$}
\hspace*{0.1cm}   %  hack to get vertical gap
\caption{Minor operations: $G$ and reductions (black,
solid edges, filled vertices), with their trials,
which equal \textcolor{blue}{$G^{\omega}$ and reductions (blue,
dashed edges, open vertices).}}
\label{fig:minor-opns}
\end{figure}

Let $f=f(v,w_1)$ be the right successor of $e$,
and let $g=g(v,w_2)$ be the left successor of $e$.
(It is possible that $f=g$, in which case $w_1=w_2$.  This occurs when
$v$ has indegree = outdegree = 1, i.e., when $e$ is a 1-loop.)

The graph $G\mino{\omega}e$ is formed by deleting $e$ and $g$ and,
if $e\not=g$, changing $g$ so that it now joins $u$ to $w_2$.
The revised edge $g$ replaces $e$ and $g$ in $A(e)$, and it replaces $g$ in $I(w_2)$.
If $\deg v=2$, then $v$ and $I(v)$ no longer exist in $G\mino{\omega}e$.
If $\deg v\not=2$, then the in-star at $v$ in $G\mino{\omega}e$ is
$I(v)\setminus\{e\}$.  The c-face containing $e'$ in $G\mino{\omega}e$
is $(C(e)\setminus\{e\})\cup C(g)$.
This operation is called \textit{$\omega$-reduction}. 
%   **** Do I need this last?  If not, see also end of para below
%        defining \omega^2-semiloop. *****
See Figure~\ref{fig:minor-opns}(a),(c).

The graph $G\mino{\omega^2}e$ is formed by deleting $e$ and,
if $e\not=f$, changing $f$ so that it now joins $u$ to $w_1$.
The revised
edge $f$ replaces $e$ and $f$ in $C(e)$, and it replaces $f$ in $I(w_1)$.
If $\deg v=2$, then $v$ and $I(v)$ no longer exist in $G\mino{\omega^2}e$.
If $\deg v\not=2$, then the in-star at $v$ in $G\mino{\omega^2}e$ is
$I(v)\setminus\{e\}$.  The a-face containing $e'$ in $G\mino{\omega^2}e$
is $(A(e)\setminus\{e\}\cup A(f))$.
This operation is called \textit{$\omega^2$-reduction}.
See Figure \ref{fig:minor-opns}(a),(d).

For all $\mu\in\{1,\omega,\omega^2\}$, $\mu$-reduction of a proper $\mu^{-1}$-semiloop
either increases the number of connected components or decreases the genus.

Each of these three
operations is a \textit{reduction} or a \textit{minor operation}.

The operations of $\omega$-reduction and $\omega^2$-reduction are
special cases of \textit{lifting} (see, e.g., \cite{mader78}), used
in the immersion relation on graphs \cite{nash-williams65}, here
restricted to cases where the two incident edges are consecutive in a face.

An alternating dimap obtained from $G$ by a sequence of minor operations
is a \textit{minor} of $G$.

These constructions are illustrated in Figure~\ref{fig:minor-opns}.

If $e$ is a triloop, then
$G\mino{1}e = G\mino{\omega}e = G\mino{\omega^2}e$.  We sometimes write
$G\mino{*}e$ for the common result of the three reductions in this case,
in order to avoid being unnecessarily specific.

If $X=(x_1,\ldots,x_k)$ is a \textit{sequence} of edges,
then we write $G\mino{\mu}X$
as a shorthand for $G\mino{\mu}x_1\mino{\mu}x_2\cdots\mino{\mu}x_k$.

It is straightforward to
translate the the above constructions for the minor operations into the language
of permutation triples.

\vbox{
\begin{thm}
\label{thm:minors-perms}
If $G$ is an alternating dimap with permutation triple
$(\sigma_1,\sigma_{\omega},\sigma_{\omega^2})$, then the
permutation triples for the three minors of $G$ are as
given in the following table.   \\

\noindent
\begin{tabular}{c|c|c}
\hline
&&   \\
$G[1]e$  &  $G[\omega]e$  &  $G[\omega^2]e$   \\
&&   \\
\hline
\parbox[t]{2in}{
\begin{eqnarray*}
\sigma_{G[1]e,1}(\sigma_{1}^{-1}(e))
& = & \sigma_{\omega^2}^{-1}(e)   \\
\sigma_{G[1]e,\omega}(\sigma_{\omega}^{-1}(e))
& = & \sigma_{\omega}(e)   \\
\sigma_{G[1]e,\omega^2}(\sigma_{\omega^2}^{-1}(e))
& = & \sigma_{\omega^2}(e)   \\
\sigma_{G[1]e,1}(\sigma_{\omega}(e)) & = & \sigma_{1}(e)   \\
&&   \\
\multicolumn{1}{l}{Otherwise:}   \\
\sigma_{G[1]e,\mu}(f)  & = &  \sigma_{\mu}(f)
\end{eqnarray*}
}
&
\parbox[t]{2in}{
\begin{eqnarray*}
\sigma_{G[\omega]e,1}(\sigma_{1}^{-1}(e))
& = & \sigma_{1}(e)   \\
\sigma_{G[\omega]e,\omega}(\sigma_{\omega}^{-1}(e))
& = & \sigma_{\omega}(e)   \\
\sigma_{G[\omega]e,\omega^2}(\sigma_{\omega^2}^{-1}(e))
& = & \sigma_{\omega}^{-1}(e)   \\
\sigma_{G[\omega]e,\omega^2}(\sigma_{1}(e)) & = & \sigma_{\omega^2}(e)   \\
&&   \\
&&   \\
\sigma_{G[\omega]e,\mu}(f)  & = &  \sigma_{\mu}(f)
\end{eqnarray*}
}
&
\parbox[t]{2in}{
\begin{eqnarray*}
\sigma_{G[\omega^2]e,1}(\sigma_{1}^{-1}(e))
& = & \sigma_{1}(e)   \\
\sigma_{G[\omega^2]e,\omega}(\sigma_{\omega}^{-1}(e))
& = & \sigma_{1}^{-1}(e)   \\
\sigma_{G[\omega^2]e,\omega^2}(\sigma_{\omega^2}^{-1}(e))
& = & \sigma_{\omega^2}(e)   \\
\sigma_{G[\omega^2]e,\omega}(\sigma_{\omega^2}(e)) & = & \sigma_{\omega}(e)   \\
&&   \\
&&   \\
\sigma_{G[\omega^2]e,\mu}(f)  & = &  \sigma_{\mu}(f)
\end{eqnarray*}
}   \\
\hline
\end{tabular}
\eopf
\end{thm}
}

We will refer to a specific equation in the table by
``Theorem \ref{thm:minors-perms}$(r,c)$'', 
where $r$ and $c$ index the row and column in which the equation appears.
For example,
Theorem \ref{thm:minors-perms}(2,1) is 
$\sigma_{G[1]e,\omega}(\sigma_{\omega}^{-1}(e))
= \sigma_{\omega}(e)$, and
Theorem \ref{thm:minors-perms}(5,2) is 
$\sigma_{G[\omega]e,\mu}(f) = \sigma_{\mu}(f)$ which holds when
the pair $\mu,f$ is not covered by any of the previous entries in that
column.

We are now in a position to establish the relationship between minors and triality.

\begin{thm}
\label{thm:trial-minor}
If $e\in E(G)$ and $\mu,\nu\in\{1,\omega,\omega^2\}$ then
\[
G^{\mu}\mino{\nu}e^{\mu} = (G\mino{\mu\nu}e)^{\mu} .
\]
\end{thm}

\pf
We first prove that
\begin{equation}
\label{eq:trial-minor-1}
G^{\omega}\mino{1}e^{\omega}  = (G\mino{\omega}e)^{\omega} .
\end{equation}

% *** DO I NEED THIS NOW? *****
If $e$ is an ultraloop, then so is $e^{\omega}$, and both sides of (\ref{eq:trial-minor-1})
are empty, so the equation is true.  So suppose that $e$ is not an ultraloop.

\begin{eqnarray*}
\sigma_{G^{\omega}\mino{1}e^{\omega},1}(\sigma_{G^{\omega},1}^{-1}(e^{\omega}))
& = &
\sigma_{G^{\omega},\omega^2}^{-1}(e^{\omega})
~~~~~~~~~~~
\hbox{(by Theorem \ref{thm:minors-perms}(1,1))}   \\
& = &
\sigma_{G,\omega}^{-1}(e)^{\omega}
~~~~~~~~~~~
\hbox{(by (\ref{eq:trial-sigma-omega}))}   \\
& = &
\sigma_{G\mino{\omega}e,\omega^2}(\sigma_{G,\omega^2}^{-1}(e))^{\omega}
~~~~~~~~~~~
\hbox{(by Theorem \ref{thm:minors-perms}(3,2))}   \\
& = &
\sigma_{(G\mino{\omega}e)^{\omega},1}(\sigma_{G,\omega^2}^{-1}(e)^{\omega})
~~~~~~~~~~~
\hbox{(by (\ref{eq:trial-sigma-omega2}))}   \\
& = &
\sigma_{(G\mino{\omega}e)^{\omega},1}(\sigma_{G^{\omega},1}^{-1}(e^{\omega}))
~~~~~~~~~~~
\hbox{(by (\ref{eq:trial-sigma-omega2}))} .
\end{eqnarray*}
Similarly, we have
\begin{eqnarray*}
\sigma_{G^{\omega}\mino{1}e^{\omega},\omega}(\sigma_{G^{\omega},\omega}^{-1}(e^{\omega}))
& = &
\sigma_{G^{\omega},\omega}(e^{\omega})
= \sigma_{G,1}(e)^{\omega}
= \sigma_{G\mino{\omega}e,1}(\sigma_{G,1}^{-1}(e))^{\omega}
= \sigma_{(G\mino{\omega}e)^{\omega},\omega}(\sigma_{G,1}^{-1}(e)^{\omega})   \\
& = &
\sigma_{(G\mino{\omega}e)^{\omega},\omega}(\sigma_{G^{\omega},\omega}^{-1}(e^{\omega})) ,   \\
\sigma_{G^{\omega}\mino{1}e^{\omega},\omega^2}(\sigma_{G^{\omega},\omega^2}^{-1}(e^{\omega}))
& = &
\sigma_{G^{\omega},\omega^2}(e^{\omega})
= \sigma_{G,\omega}(e)^{\omega}
= \sigma_{G\mino{\omega}e,\omega}(\sigma_{G,\omega}^{-1}(e))^{\omega}
= \sigma_{(G\mino{\omega}e)^{\omega},\omega^2}(\sigma_{G,\omega}^{-1}(e)^{\omega})   \\
& = &
\sigma_{(G\mino{\omega}e)^{\omega},\omega^2}(\sigma_{G^{\omega},\omega^2}^{-1}(e^{\omega})) ,   \\
\sigma_{G^{\omega}\mino{1}e^{\omega},1}(\sigma_{G^{\omega},\omega}(e^{\omega}))
& = &
\sigma_{G^{\omega},1}(e^{\omega})
= \sigma_{G,\omega^2}(e)^{\omega}
= \sigma_{G\mino{\omega}e,\omega^2}(\sigma_{G,1}(e))^{\omega}
= \sigma_{(G\mino{\omega}e)^{\omega},1}(\sigma_{G,1}(e)^{\omega})   \\
& = &
\sigma_{(G\mino{\omega}e)^{\omega},1}(\sigma_{G^{\omega},\omega}(e^{\omega})) .
\end{eqnarray*}

It remains to consider the cases where
\begin{equation}
\label{eq:exclude-the-four-special-cases1}
(\mu,f^{\omega})\not\in\{(1,\sigma_{G^{\omega},1}^{-1}(e^{\omega})),
(\omega,\sigma_{G^{\omega},\omega}^{-1}(e^{\omega})),
(\omega^2,\sigma_{G^{\omega},\omega^2}^{-1}(e^{\omega})),
(1,\sigma_{G^{\omega},\omega}(e^{\omega}))\} .
\end{equation}
Equivalently,
\begin{equation}
\label{eq:exclude-the-four-special-cases2}
(\mu\omega^2,f)\not\in\{(1,\sigma_{G,1}^{-1}(e)),
(\omega,\sigma_{G,\omega}^{-1}(e)),
(\omega^2,\sigma_{G,\omega^2}^{-1}(e)),
(\omega^2,\sigma_{G,1}(e))\} .
\end{equation}
In these cases, we have
\begin{eqnarray*}
\sigma_{G^{\omega}\mino{1}e^{\omega},\mu}(f^{\omega})
& = &
\sigma_{G^{\omega},\mu}(f^{\omega})
~~~~~~~~~~
\hbox{(by Theorem \ref{thm:minors-perms}(5,1),
using (\ref{eq:exclude-the-four-special-cases1}))}   \\
& = &
\sigma_{G,\mu\omega^2}(f)^{\omega}
~~~~~~~~~~
\hbox{(by
(\ref{eq:trial-sigma-1})--(\ref{eq:trial-sigma-omega2}))}   \\
& = &
\sigma_{G\mino{\omega}e,\mu\omega^2}(f)^{\omega}
~~~~~~~~~~
\hbox{(by Theorem \ref{thm:minors-perms}(5,2),
using (\ref{eq:exclude-the-four-special-cases2}))}   \\
& = &
\sigma_{(G\mino{\omega}e)^{\omega},\mu}(f^{\omega})
~~~~~~~~~~
\hbox{(by
(\ref{eq:trial-sigma-1})--(\ref{eq:trial-sigma-omega2}))} .
\end{eqnarray*}

We have now shown that, for all $\mu\in\{1,\omega,\omega^2\}$
and $f\in E(G)$,
\[
\sigma_{G^{\omega}\mino{1}e^{\omega},\mu}(f^{\omega})
=
\sigma_{(G\mino{\omega}e)^{\omega},\mu}(f^{\omega}) .
\]
It follows that
\[
(\sigma_{G^{\omega}\mino{1}e^{\omega},1},
\sigma_{G^{\omega}\mino{1}e^{\omega},\omega},
\sigma_{G^{\omega}\mino{1}e^{\omega},\omega^2})
=
(\sigma_{(G\mino{\omega}e)^{\omega},1},
\sigma_{(G\mino{\omega}e)^{\omega},\omega},
\sigma_{(G\mino{\omega}e)^{\omega},\omega^2}),
\]
or in other words,
\[
G^{\omega}\mino{1}e^{\omega}  =  (G\mino{\omega}e)^{\omega} .
\]

Similar arguments show that
\begin{eqnarray*}
G^{\omega}\mino{\omega}e^{\omega}  & = &  (G\mino{\omega^2}e)^{\omega} ,   \\
G^{\omega}\mino{\omega^2}e^{\omega}  & = &  (G\mino{1}e)^{\omega} .
\end{eqnarray*}
From these it follows that
\begin{eqnarray*}
G^{\omega^2}\mino{1}e^{\omega^2}  & = &  (G\mino{\omega^2}e)^{\omega^2} ,   \\
G^{\omega^2}\mino{\omega}e^{\omega^2}  & = &  (G\mino{1}e)^{\omega^2} ,   \\
G^{\omega^2}\mino{\omega^2}e^{\omega^2}  & = &  (G\mino{\omega}e)^{\omega^2} .
\end{eqnarray*}
\eopf   \\

% **** remark somewhere on treating  \omega  as complex number?  ****

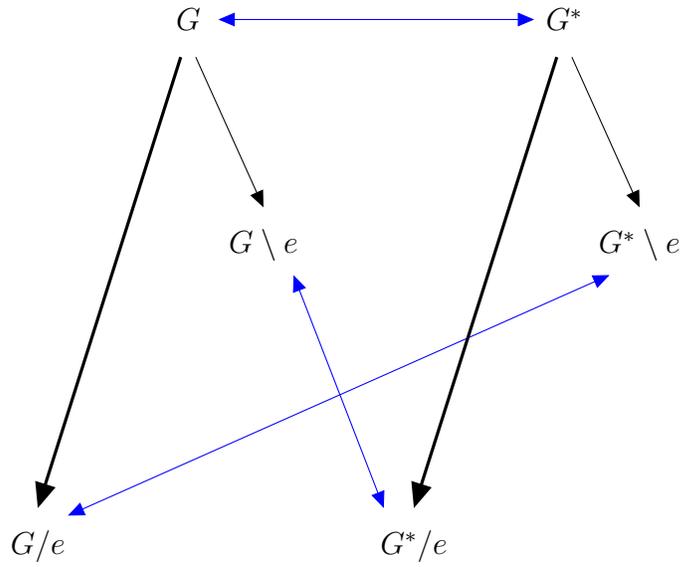
\begin{figure}
\begin{center}
\begin{tikzpicture}
\begin{scope}[>=triangle 45]
\draw (3,8) node {$G$};
\draw (1,1) node {$G/e$};
\draw (4,5) node {$G\setminus e$};
\draw (8,8) node {$G^*$};
\draw (6,1) node {$G^*/e$};
\draw (9,5) node {$G^*\setminus e$};
\draw[blue,<->] (3.4,8) -- (7.6,8);   % duality
\draw[blue,<->] (1.4,1.4) -- (8.6,4.6);   % duality
\draw[blue,<->] (4.4,4.6) -- (5.6,1.5);   % duality
\draw[->,very thick] (2.9,7.5) -- (1,1.5);   % contraction
\draw[->] (3.1,7.5) -- (4,5.5);   % deletion
\draw[->,very thick] (7.9,7.5) -- (6,1.5);   % contraction
\draw[->] (8.1,7.5) -- (9,5.5);   % deletion
\end{scope}
\end{tikzpicture}
\end{center}
\caption{The relationship between the ordinary duality and minor operations.}
\label{fig:duality-minors}

\end{figure}

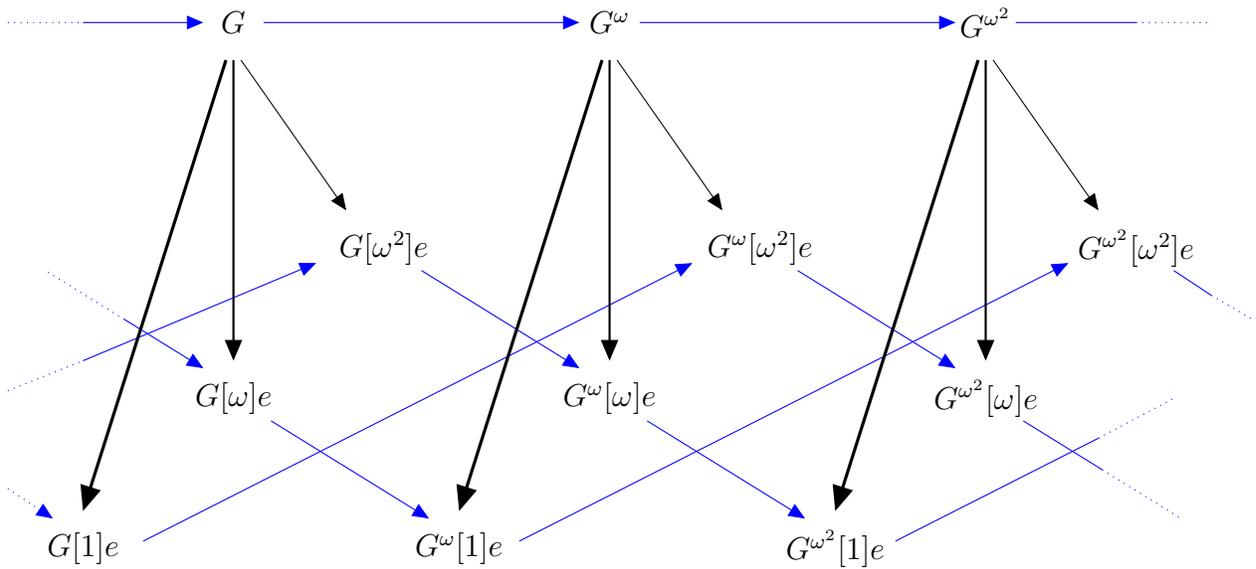
\begin{figure}
\begin{tikzpicture}
\begin{scope}[>=triangle 45]
\draw (3,8) node {$G$};
\draw (1,1) node {$G\mino{1}e$};
\draw (3,3) node {$G\mino{\omega}e$};
\draw (5,5) node {$G\mino{\omega^2}e$};
\draw (8,8) node {$G^{\omega}$};
\draw (6,1) node {$G^{\omega}\mino{1}e$};
\draw (8,3) node {$G^{\omega}\mino{\omega}e$};
\draw (10,5) node {$G^{\omega}\mino{\omega^2}e$};
\draw (13,8) node {$G^{\omega^2}$};
\draw (11,1) node {$G^{\omega^2}\mino{1}e$};
\draw (13,3) node {$G^{\omega^2}\mino{\omega}e$};
\draw (15,5) node {$G^{\omega^2}\mino{\omega^2}e$};
\draw[blue,dotted] (0,8) -- (1,8);   % trinity
\draw[blue,->] (1,8) -- (2.6,8);   % trinity
\draw[blue,->] (3.4,8) -- (7.6,8);   % trinity
\draw[blue,->] (8.4,8) -- (12.6,8);   % trinity
\draw[blue] (13.4,8) -- (15,8);   % trinity
\draw[blue,dotted] (15,8) -- (16,8);   % trinity

\draw[blue,dotted] (0,3.1) -- (1,3.5);
\draw[blue,->] (1,3.5) -- (4.1,4.8);
\draw[blue,dotted,->] (0,1.8) -- (0.6,1.4);
\draw[blue,->] (1.55,4.05) -- (2.6,3.4);
\draw[blue,dotted] (1.55,4.05) -- (0.5,4.7);
%  \draw[blue,->] (0.5,4.7) -- (2.6,3.4);
\draw[blue,->] (1.8,1.1) -- (9.1,4.8);
\draw[blue,->] (3.5,2.7) -- (5.6,1.4);
\draw[blue,->] (5.5,4.7) -- (7.6,3.4);
\draw[blue,->] (6.8,1.1) -- (14.1,4.8);
\draw[blue,->] (8.5,2.7) -- (10.6,1.4);
\draw[blue,->] (10.5,4.7) -- (12.6,3.4);
\draw[blue] (11.8,1.1) -- (14.5,2.47);
\draw[blue,dotted] (14.5,2.47) -- (15.5,3);
\draw[blue] (13.5,2.7) -- (14.55,2.05);
\draw[blue,dotted] (14.55,2.05) -- (15.6,1.4);
\draw[blue] (15.5,4.7) -- (16,4.37);
\draw[blue,dotted] (16,4.37) -- (16.55,4.05);

\draw[->,very thick] (2.9,7.5) -- (1,1.5);   % 1-reduction (contraction)
\draw[->,thick] (3,7.5) -- (3,3.5);   % \omega -reduction
\draw[->] (3.1,7.5) -- (4.5,5.5);   % \omega^2 -reduction
\draw[->,very thick] (7.9,7.5) -- (6,1.5);   % 1-reduction (contraction)
\draw[->,thick] (8,7.5) -- (8,3.5);   % \omega -reduction
\draw[->] (8.1,7.5) -- (9.5,5.5);   % \omega^2 -reduction
\draw[->,very thick] (12.9,7.5) -- (11,1.5);   % 1-reduction (contraction)
\draw[->,thick] (13,7.5) -- (13,3.5);   % \omega -reduction
\draw[->] (13.1,7.5) -- (14.5,5.5);   % \omega^2 -reduction
\end{scope}
\end{tikzpicture}
\caption{The relationship between triality and the three minor operations.  The diagram wraps around at its left and right sides.}
\label{fig:trinity-minors}
\end{figure}

Theorem \ref{thm:trial-minor} extends the classical relationship between duality and minors:
\begin{eqnarray*}
G\setminus e  & = &  (G/e)^* ,   \\
G/e  & = &  (G\setminus e)^* .
\end{eqnarray*}
The classical relationship is illustrated in Figure \ref{fig:duality-minors},
while Theorem \ref{thm:trial-minor} is illustrated in Figure \ref{fig:trinity-minors}

The mappings $\sigma_{H,\mu}$, where $H$ ranges over all minors of $G$
and $\mu\in\{1,\omega,\omega^2\}$, together generate
(under composition) an inverse semigroup,
which we denote by IS($G$).   \\

\noindent\textbf{Problem}

Describe IS($G$), and classify it among known types of inverse semigroup.

\section{Non-commutativity}
\label{sec:non-commut}

Deletion and contraction are well known to commute, in the sense
that, for any graph $G$ and any distinct $e,f\in E(G)$, we have
\begin{eqnarray*}
G\setminus e\setminus f  & = &  G\setminus f\setminus e ,   \\
G/e/f  & = &  G/f/e ,   \\
G\setminus e/f  & = &  G/f\setminus e .
\end{eqnarray*}
The variants of these operations for embedded graphs, where deletion/contraction
of an edge is accompanied by appropriate modifications to the embedding,
also commute \cite[p.\ 103]{mohar-thomassen01}.

Perhaps surprisingly, the reductions we have introduced for
alternating dimaps do not always commute, though they do in most
situations.
In this section we investigate the circumstances under which
the reductions do or do not commute.

We first show that the reductions always commute if one of the
edges involved is a triloop.

\begin{lemma}
\label{lemma:reductions-commute-loops}
If $f$ is a triloop, then for any $\nu\in\{1,\omega,\omega^2\}$,
\[
G\mino{1}e\mino{\nu}f = G\mino{\nu}f\mino{1}e .
\]
\end{lemma}

\pf
If $f$ is an $\omega$-loop or an $\omega^2$-loop,
then $\nu$-reduction of $f$
just amounts to deletion of $f$.
So
\[
G\mino{1}e\mino{\nu}f = G/e\setminus f =
G\setminus f/e = G\mino{\nu}f\mino{1}e ,
\]
where the middle equality follows from the fact that deletion
and contraction commute.  (These deletion and contraction
operations are for embedded graphs, and give surface minors, and so are
not the usual deletion and contraction operations for abstract graphs.
But it is still easy to see that they commute when $f$ is a loop that
is contractible in the surface.)

If $f$ is a 1-loop, then $\nu$-reduction of $f$
just amounts to contraction of $f$.
So
\[
G\mino{1}e\mino{\nu}f = G/e/f =
G/f/e = G\mino{\nu}f\mino{1}e .
\]
\eopf

\begin{thm}
\label{thm:reductions-commute-loops}
If $f$ is a triloop and $\mu,\nu\in\{1,\omega,\omega^2\}$
then
\[
G\mino{\mu}e\mino{\nu}f = G\mino{\nu}f\mino{\mu}e .
\]
\end{thm}

\pf
\begin{eqnarray*}
G\mino{\mu}e\mino{\nu}f
& = &
(G^{\mu}\mino{1}e^{\mu}\mino{\nu\mu^{-1}}f^{\mu})^{\mu^{-1}}
~~~~~~~~~~~ \hbox{(by Theorem \ref{thm:trial-minor})}   \\
& = &
(G^{\mu}\mino{\nu\mu^{-1}}f^{\mu}\mino{1}e^{\mu})^{\mu^{-1}}
~~~~~~~~~~~ \hbox{(by Lemma \ref{lemma:reductions-commute-loops})}   \\
& = &
G\mino{\nu}f\mino{\mu}e
~~~~~~~~~~~~~~~~~~~~~~~~~~ \hbox{(by Theorem \ref{thm:trial-minor})} .
\end{eqnarray*}
\eopf

We next show that two reductions of the same type always commute.

\begin{thm}
\label{thm:same-reductions-commute}
For all $\mu\in\{1,\omega,\omega^2\}$,
\[
G\mino{\mu}e\mino{\mu}f  =
G\mino{\mu}f\mino{\mu}e .
\]
\end{thm}

\pf
We show that
\begin{equation}
\label{eq:contractions-commute}
G\mino{1}e\mino{1}f  =
G\mino{1}f\mino{1}e ,
\end{equation}
which takes up most of the proof, and then use triality to complete it.

To show (\ref{eq:contractions-commute}),
we will show that, for all $\mu\in\{1,\omega,\omega^2\}$ and all
$g\in E(G)\setminus\{e,f\}$,
\begin{equation}
\label{eq:contractions-commute-mu-g}
\sigma_{G\mino{1}e\mino{1}f,\mu}(g)  =
\sigma_{G\mino{1}f\mino{1}e,\mu}(g) .
\end{equation}

We first do this for $\mu\in\{\omega,\omega^2\}$, which we now assume.

Most situations are covered by the following reasoning:
\begin{eqnarray*}
\sigma_{G\mino{1}e\mino{1}f,\mu}(g)
& = &
\sigma_{G\mino{1}e,\mu}(g)
~~~~~~~ \hbox{if $g\not=\sigma_{G\mino{1}e,\mu}^{-1}(f)$}   \\
& = &
\sigma_{G,\mu}(g)
~~~~~~~ \hbox{if $g\not=\sigma_{G,\mu}^{-1}(e)$}   \\
& = &
\sigma_{G\mino{1}f,\mu}(g)
~~~~~~~ \hbox{if $g\not=\sigma_{G,\mu}^{-1}(f)$}   \\
& = &
\sigma_{G\mino{1}f\mino{1}e,\mu}(g)
~~~~~~~ \hbox{if $g\not=\sigma_{G\mino{1}f,\mu}^{-1}(e)$} ,
\end{eqnarray*}
by four applications of Theorem \ref{thm:minors-perms}(5,1),
since the conditions on $g$ ensure that cases (2,1) and (3,1)
(according as $\mu=\omega$ or $\mu=\omega^2$)
of that Theorem do not apply.

We now deal with situations where the above conditions on $g$ are not met.
We have, apparently, four exceptional values of $g$.  We consider each
in turn.

Firstly, suppose $g=\sigma_{G,\mu}^{-1}(e)$.

In this case, we must assume that $f\not=\sigma_{G,\mu}^{-1}(e)$,
else $g=f$ and $g\not\in\mathop{\hbox{dom}}\sigma_{G\mino{1}e\mino{1}f,\mu}$.

Consider $\sigma_{G\mino{1}e\mino{1}f,\mu}(g)$.

If $f=\sigma_{G,\mu}(e)$ then
$\sigma_{G\mino{1}e,\mu}^{-1}(f) = \sigma_{G,\mu}^{-1}(e)$,
by Theorem \ref{thm:minors-perms}(2,1),(3,1).  This justifies the first
step in the following.
\begin{eqnarray*}
\sigma_{G\mino{1}e\mino{1}f,\mu}(\sigma_{G,\mu}^{-1}(e))
& = &
\sigma_{G\mino{1}e\mino{1}f,\mu}(\sigma_{G\mino{1}e,\mu}^{-1}(f))   \\
& = &
\sigma_{G\mino{1}e,\mu}(f)   ~~~~~~~~
\hbox{(by Theorem \ref{thm:minors-perms}(2,1) or (3,1))}   \\
& = &
\sigma_{G,\mu}(f)   ~~~~~~~~
\hbox{(by Theorem \ref{thm:minors-perms}(5,1),
since $f\not=\sigma_{G,\mu}^{-1}(e)$)}.
\end{eqnarray*}
On the other hand, if $f\not=\sigma_{G,\mu}(e)$, then
$\sigma_{G\mino{1}e,\mu}^{-1}(f) = \sigma_{G,\mu}^{-1}(f)$,
by Theorem \ref{thm:minors-perms}(5,1).  This in turn does not equal
$\sigma_{G,\mu}^{-1}(e)$, since $e\not=f$ and $\sigma_{G,\mu}^{-1}$ is
a bijection.  So
\[
\sigma_{G\mino{1}e\mino{1}f,\mu}(\sigma_{G,\mu}^{-1}(e))  =
\sigma_{G\mino{1}e,\mu}(\sigma_{G,\mu}^{-1}(e))  =
\sigma_{G,\mu}(e) ,
\]
by Theorem \ref{thm:minors-perms}(5,1) then (2,1) or (3,1).

Now consider $\sigma_{G\mino{1}f\mino{1}e,\mu}(g)$.
Observe that
$\sigma_{G\mino{1}f,\mu}^{-1}(e) = \sigma_{G,\mu}^{-1}(e)$,
by Theorem \ref{thm:minors-perms}(5,1), since
$f\not=\sigma_{G,\mu}^{-1}(e)$.  Therefore
\begin{eqnarray*}
\sigma_{G\mino{1}f\mino{1}e,\mu}(\sigma_{G,\mu}^{-1}(e))
& = &
\sigma_{G\mino{1}f\mino{1}e,\mu}(\sigma_{G\mino{1}f,\mu}^{-1}(e))   \\
& = &
\sigma_{G\mino{1}f,\mu}(e)  ~~~~~~~~~~~
\hbox{(by Theorem \ref{thm:minors-perms}(2,1) or (3,1))}   \\
& = &
\left\{
\begin{array}{ll}
\sigma_{G,\mu}(f) ,  &  \hbox{if $e=\sigma_{G,\mu}^{-1}(f)$},  \\  
\sigma_{G,\mu}(e) ,  &  \hbox{otherwise},\caseSpace
\end{array}
\right.
\end{eqnarray*}
by Theorem \ref{thm:minors-perms}(2,1) or (3,1), and (5,1).

So $\sigma_{G\mino{1}e\mino{1}f,\mu}(g) = \sigma_{G\mino{1}f\mino{1}e,\mu}(g)$
when $g=\sigma_{G,\mu}^{-1}(e)$.

Secondly, suppose $g=\sigma_{G,\mu}^{-1}(f)$.  This can be treated the
same as the first case, except that $e$ and $f$ are swapped throughout.

Thirdly and fourthly, the remaining two exceptional values of $g$, namely
$\sigma_{G\mino{1}e,\mu}^{-1}(f)$ and $\sigma_{G\mino{1}f,\mu}^{-1}(e)$,
are really nothing new, for application of Theorem \ref{thm:minors-perms}
gives
\begin{eqnarray*}
\sigma_{G\mino{1}e,\mu}^{-1}(f)
& = &
\left\{
\begin{array}{ll}
\sigma_{G,\mu}^{-1}(e) ,  &  \hbox{if $f=\sigma_{G,\mu}(e)$},  \\  
\sigma_{G,\mu}^{-1}(f) ,  &  \hbox{otherwise};\caseSpace
\end{array}
\right.   \\
\sigma_{G\mino{1}f,\mu}^{-1}(e)
& = &
\left\{
\begin{array}{ll}
\sigma_{G,\mu}^{-1}(f) ,  &  \hbox{if $e=\sigma_{G,\mu}(f)$},\caseSpace  \\  
\sigma_{G,\mu}^{-1}(e) ,  &  \hbox{otherwise}.\caseSpace
\end{array}
\right.
\end{eqnarray*}
Thus, in any event, each of these two values of $g$ actually falls into
one of the first two cases.  This completes the treatment of the exceptional
values of $g$ (apparently four in number, but really just two).

We have now proved
(\ref{eq:contractions-commute-mu-g}) for $\mu\in\{\omega,\omega^2\}$.
But it then follows immediately for $\mu=1$ too, since
$\sigma_{H,1}=\sigma_{H,\omega^2}^{-1}\circ\sigma_{H,\omega}^{-1}$ for any $H$.
So (\ref{eq:contractions-commute-mu-g}) holds for all $\mu$ and all $g$,
which establishes (\ref{eq:contractions-commute}).

Now that we know contractions commute, we can use triality to show that
any two reductions of the same type commute.
For any $\mu\in\{1,\omega,\omega^2\}$,
\[
G\mino{\mu}e\mino{\mu}f = 
(G^{\mu}\mino{1}e^{\mu}\mino{1}f^{\mu})^{\mu^{-1}} = 
(G^{\mu}\mino{1}f^{\mu}\mino{1}e^{\mu})^{\mu^{-1}} = 
G\mino{\mu}f\mino{\mu}e . 
\]
\eopf   \\

\begin{figure}
\begin{center}
\begin{tikzpicture}
\begin{scope}[>=triangle 45]
%  vertex  u
\coordinate (u) at (2,3) {};
\draw[fill] (u) circle (0.23cm);
%%  \draw (1.9,2.75) node {$u$};
%  other edges at  u
\coordinate (a) at (1.06,3.42) {};
\coordinate (b) at (2.71,3.71) {};
\draw[->>] (a) -- (u);
\draw[->>] (b) -- (u);
\draw[dotted,thick] (2.35,3.35) arc (45:150:0.5);
%  vertex  v
\coordinate (v) at (1.5,1.5) {};
\draw[fill] (v) circle (0.23cm);
%%  \draw (1.3,1.65) node {$v$};
%  edge  f(u,v)
\draw[->>] (u) -- (v);
\draw (2,2.25) node {$f$};
%  vertices  w1, w2
\coordinate (w1) at (0.5,1.5) {};
\coordinate (w2) at (1,0.65) {};
%  other edges at  v
\draw[->] (v) -- (w1);
\draw[->>] (w2) -- (v);
\draw[dotted,thick] (0.75,1.5) arc (180:240:0.75);
%  vertex  x
\coordinate (x) at (2.95,1) {};
%%  \draw (3,1.35) node {$x$};
\draw[fill=greenish,opaque] (x) circle (0.23cm);
%   ... repeated at end to ensure colour is on top!
\coordinate (x1) at (2.95,0.25) {};
\coordinate (x2) at (3.5,0.45) {};
\draw[->>] (v) -- (x);
\draw (2.2,0.9) node {$e$};
\draw[greenish,->] (x) -- (x1);
\draw[greenish,->>] (x2) -- (x);
\draw[greenish,dotted,thick] (2.95,0.5) arc (270:315:0.5);
%   vertex  y  and its edges
\coordinate (y) at (4.375,1.5) {};
%%  \draw (4.35,1.15) node {$y$};
\draw[fill] (y) circle (0.23cm);
\coordinate (y1) at (5,0.875) {};
\coordinate (y2) at (5,2.125) {};
\draw[->>] (x) -- (y);
\draw[->] (y) -- (y1);
\draw[->] (y) -- (y2);
\draw[dotted,thick] (4.875,1) arc (-45:45:0.7071);
\draw[fill=greenish,opaque] (x) circle (0.23cm);
%   ... repeated from earlier to ensure colour is on top!
\end{scope}
\end{tikzpicture}
\end{center}
\begin{center}
$G$
\end{center}
\begin{tikzpicture}
\begin{scope}[>=triangle 45]
%  vertex  u
\coordinate (u) at (2,3) {};
\draw[fill] (u) circle (0.23cm);
%%  \draw (1.9,2.75) node {$u$};
%  other edges at  u
\coordinate (a) at (1.06,3.42) {};
\coordinate (b) at (2.71,3.71) {};
\draw[->>] (a) -- (u);
\draw[->>] (b) -- (u);
\draw[dotted,thick] (2.35,3.35) arc (45:150:0.5);
%  vertex  v
\coordinate (v) at (1.5,1.5) {};
\draw[fill] (v) circle (0.23cm);
%%  \draw (1.3,1.65) node {$v$};
%  edge  f(u,v)
%  \draw[->>] (u) -- (v);
%  \draw (2,2.25) node {$f$};
%  vertices  w1, w2
\coordinate (w1) at (0.5,1.5) {};
\coordinate (w2) at (1,0.65) {};
%  other edges at  v
\draw[->] (v) -- (w1);
\draw[->>] (w2) -- (v);
\draw[dotted,thick] (0.75,1.5) arc (180:240:0.75);
%  vertex  x:  now merged with  u
\coordinate (x) at (2.1,2.9) {};
%%  \draw (3,1.35) node {$x$};
\draw[fill=greenish,opaque] (x) circle (0.23cm);
\coordinate (x1) at (2.1,2.15) {};
\coordinate (x2) at (2.7,2.35) {};
%  \draw[->>] (v) -- (x);
%  \draw (2.3,0.9) node {$e$};
\draw[greenish,->] (x) -- (x1);
\draw[greenish,->>] (x2) -- (x);
\draw[greenish,dotted,thick] (2.1,2.4) arc (270:315:0.5);
%   vertex  y  and its edges
\coordinate (y) at (4.375,1.5) {};
%%  \draw (4.35,1.15) node {$y$};
\draw[fill] (y) circle (0.23cm);
\coordinate (y1) at (5,0.875) {};
\coordinate (y2) at (5,2.125) {};
\draw[->>] (x) -- (y);
\draw[->] (y) -- (y1);
\draw[->] (y) -- (y2);
\draw[dotted,thick] (4.875,1) arc (-45:45:0.7071);
\draw[fill=greenish,opaque] (x) circle (0.23cm);
%   ... repeated from earlier to ensure colour is on top!
\end{scope}
\end{tikzpicture}
\hfill
\begin{tikzpicture}
\begin{scope}[>=triangle 45]
%  vertex  u
\coordinate (u) at (2,3) {};
\draw[fill] (u) circle (0.23cm);
%%  \draw (1.9,2.75) node {$u$};
%  other edges at  u
\coordinate (a) at (1.06,3.42) {};
\coordinate (b) at (2.71,3.71) {};
\draw[->>] (a) -- (u);
\draw[->>] (b) -- (u);
\draw[dotted,thick] (2.35,3.35) arc (45:150:0.5);
%  vertex  v
\coordinate (v) at (1.5,1.5) {};
\draw[fill] (v) circle (0.23cm);
%%  \draw (1.3,1.65) node {$v$};
%  edge  f(u,v)
%  \draw[->>] (u) -- (v);
%  \draw (2,2.25) node {$f$};
%  vertices  w1, w2
\coordinate (w1) at (0.5,1.5) {};
\coordinate (w2) at (1,0.65) {};
%  other edges at  v
\draw[->] (v) -- (w1);
\draw[->>] (w2) -- (v);
\draw[dotted,thick] (0.75,1.5) arc (180:240:0.75);
%  vertex  x:  now merged with  v
\coordinate (x) at (1.7,1.5) {};
%%  \draw (3,1.35) node {$x$};
\draw[fill=greenish,opaque] (x) circle (0.23cm);
\coordinate (x1) at (1.7,0.75) {};
\coordinate (x2) at (2.25,0.95) {};
%  \draw[->>] (v) -- (x);
%  \draw (2.3,0.9) node {$e$};
\draw[greenish,->] (x) -- (x1);
\draw[greenish,->>] (x2) -- (x);
\draw[greenish,dotted,thick] (1.7,1) arc (270:315:0.5);
%   vertex  y  and its edges
\coordinate (y) at (4.375,1.5) {};
%%  \draw (4.35,1.15) node {$y$};
\draw[fill] (y) circle (0.23cm);
\coordinate (y1) at (5,0.875) {};
\coordinate (y2) at (5,2.125) {};
\draw[->>] (u) -- (y);    %  changed from  (x) -- (y)
\draw[->] (y) -- (y1);
\draw[->] (y) -- (y2);
\draw[dotted,thick] (4.875,1) arc (-45:45:0.7071);
\draw[fill=greenish,opaque] (x) circle (0.23cm);
%   ... repeated from earlier to ensure colour is on top!
\end{scope}
\end{tikzpicture}
\begin{center}
$G\mino{\omega}f\mino{1}e$
\hspace{8cm}
$G\mino{1}e\mino{\omega}f$
\end{center}
\caption{Non-commutativity of minor operations}
\label{fig:non-commut}
\end{figure}
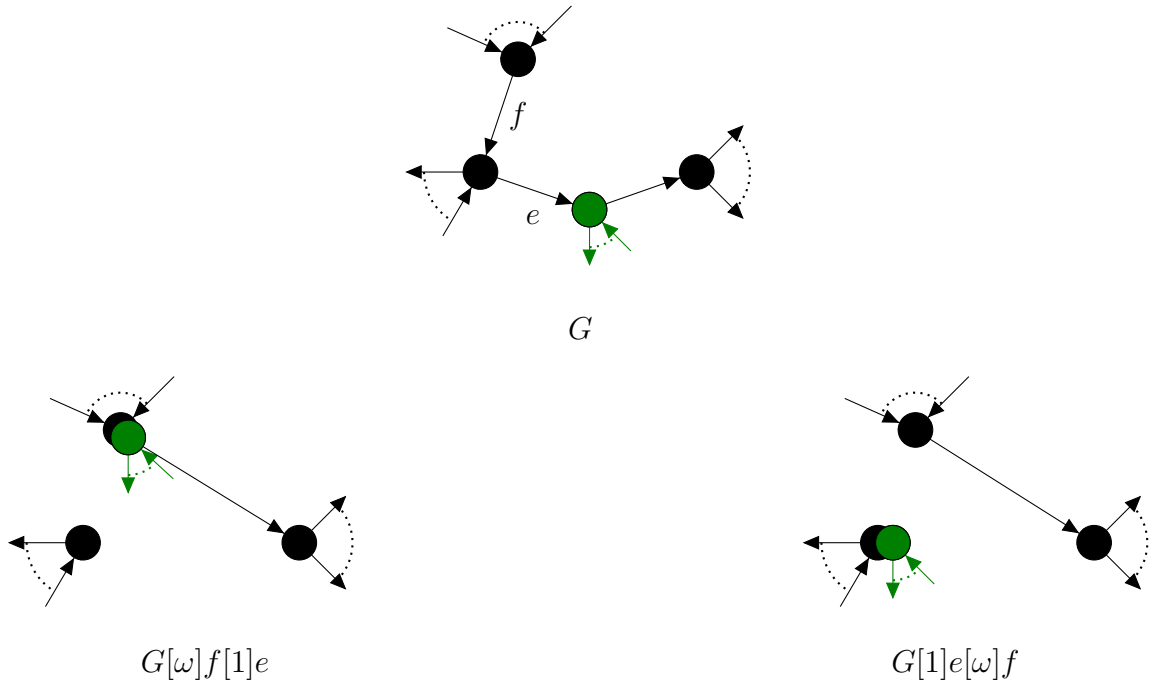

However, it is not always the case that two reductions commute.
Figure \ref{fig:non-commut} illustrates the fact that, in general,
if $f=\sigma_{G,\omega}(e)$ then
$G\mino{1}e\mino{\omega}f\not=G\mino{\omega}f\mino{1}e$.
By triality, it follows that if $f=\sigma_{G,\omega^2}(e)$ then in general
$G\mino{\omega^2}e\mino{1}f\not=G\mino{1}f\mino{\omega^2}e$,
and if $f=\sigma_{G,1}(e)$ then in general
$G\mino{\omega}e\mino{\omega^2}f\not=G\mino{\omega^2}f\mino{\omega}e$.

Most of the remainder of this section is devoted to showing that these
exceptional cases are the only situations where reductions do not commute.

We will need some lemmas.

\begin{lemma}
\label{lemma:minors-perms-inv}
(a)
\[
\sigma_{G\mino{1}e,\omega}^{-1}(h) = \left\{
\begin{array}{ll}
\sigma_{G,\omega}^{-1}(e),  &  \hbox{if $h=\sigma_{G,\omega}(e)$;}   \\
\sigma_{G,\omega}^{-1}(h),  &  \hbox{otherwise.}\caseSpace
\end{array}
\right.
\]
(b)
\[
\sigma_{G\mino{\omega}f,\omega}^{-1}(h) = \left\{
\begin{array}{ll}
\sigma_{G,\omega}^{-1}(f),  &  \hbox{if $h=\sigma_{G,\omega}(f)$;}   \\
\sigma_{G,\omega}^{-1}(h),  &  \hbox{otherwise.}\caseSpace
\end{array}
\right.
\]
(c)
\[
\sigma_{G\mino{1}e,1}^{-1}(h) = \left\{
\begin{array}{ll}
\sigma_{G,\omega}(e),  &  \hbox{if $h=\sigma_{G,1}(e)$;}   \\
\sigma_{G,1}^{-1}(e),  &  \hbox{if $h=\sigma_{G,\omega^2}^{-1}(e)$;}\caseSpace   \\
\sigma_{G,1}^{-1}(h),  &  \hbox{otherwise.}\caseSpace
\end{array}
\right.
\]
(d)
\[
\sigma_{G\mino{\omega}f,1}^{-1}(h) = \left\{
\begin{array}{ll}
\sigma_{G,1}^{-1}(f),  &  \hbox{if $h=\sigma_{G,1}(f)$;}   \\
\sigma_{G,1}^{-1}(h),  &  \hbox{otherwise.}\caseSpace
\end{array}
\right.
\]
(e)
\[
\sigma_{G\mino{1}e,\omega^2}^{-1}(h) = \left\{
\begin{array}{ll}
\sigma_{G,\omega^2}^{-1}(e),  &  \hbox{if $h=\sigma_{G,\omega^2}(e)$;}   \\
\sigma_{G,\omega^2}^{-1}(h),  &  \hbox{otherwise.}\caseSpace
\end{array}
\right.
\]
(f)
\[
\sigma_{G\mino{\omega}f,\omega^2}^{-1}(h) = \left\{
\begin{array}{ll}
\sigma_{G,\omega^2}^{-1}(f),  &  \hbox{if $h=\sigma_{G,\omega}^{-1}(f)$;}   \\
\sigma_{G,1}(f),  &  \hbox{if $h=\sigma_{G,\omega^2}(f)$;}\caseSpace   \\
\sigma_{G,\omega^2}^{-1}(h),  &  \hbox{otherwise.}\caseSpace
\end{array}
\right.
\]
\end{lemma}

\pf
Immediate from :
(a) Theorem \ref{thm:minors-perms}(2,1),(5,1);
(b) Theorem \ref{thm:minors-perms}(2,2),(5,2);
(c) Theorem \ref{thm:minors-perms}(1,1),(4,1),(5,1);
(d) Theorem \ref{thm:minors-perms}(1,2),(5,2);
(e) Theorem \ref{thm:minors-perms}(3,1),(5,1);
(f) Theorem \ref{thm:minors-perms}(3,2),(4,2),(5,2).
\eopf

\begin{lemma}
\label{lemma:reductions-commute-sigma-omega}
If $f\not=\sigma_{G,\omega}(e)$ then
\[
\sigma_{G\mino{1}e\mino{\omega}f,\omega} =
\sigma_{G\mino{\omega}f\mino{1}e,\omega}.
\]
\end{lemma}

\pf
The proof has some similarities to
that of Theorem \ref{thm:same-reductions-commute},
but is significantly more complicated.

We prove that
\[
\sigma_{G\mino{1}e\mino{\omega}f,\omega}(g) =
\sigma_{G\mino{\omega}f\mino{1}e,\omega}(g)
\]
for all $g\in E(G)\setminus\{e,f\}$.

If $g\not\in\{
\sigma_{G,\omega}^{-1}(e),
\sigma_{G,\omega}^{-1}(f),
\sigma_{G\mino{1}e,\omega}^{-1}(f),
\sigma_{G\mino{\omega}f,\omega}^{-1}(e)
\}$ then
\[
\sigma_{G\mino{1}e\mino{\omega}f,\omega}(g) =
\sigma_{G\mino{1}e,\omega}(g) =
\sigma_{G,\omega}(g) =
\sigma_{G\mino{\omega}f,\omega}(g) =
\sigma_{G\mino{\omega}f\mino{1}e,\omega}(g) ,
\]
by Theorem \ref{thm:minors-perms}(5,1),(5,2).

This leaves four special cases for $g$, which we will consider in turn,
after noting some facts which we will use repeatedly.

Observe that the Lemma's condition, $f\not=\sigma_{G,\omega}(e)$, implies
\begin{equation}
\label{eq:implied-by-forbidden-e-f}
\sigma_{G\mino{1}e,\omega}^{-1}(f) \not= \sigma_{G,\omega}^{-1}(e) ,
\end{equation}
by Lemma \ref{lemma:minors-perms-inv}(a).

Also, by Lemma \ref{lemma:minors-perms-inv}(b),
\begin{equation}
\label{eq:implied-by-e-left-f}
e=\sigma_{G,\omega}(f)    ~~~~ \Longleftrightarrow ~~~~
\sigma_{G\mino{\omega}f,\omega}^{-1}(e) = \sigma_{G,\omega}^{-1}(f) .
\end{equation}

Case 1:
$g = \sigma_{G,\omega}^{-1}(e)$.

\begin{eqnarray*}
\sigma_{G\mino{1}e\mino{\omega}f,\omega}(\sigma_{G,\omega}^{-1}(e))
& = &
\sigma_{G\mino{1}e,\omega}(\sigma_{G,\omega}^{-1}(e))
~~~~~~~
\hbox{(by Theorem \ref{thm:minors-perms}(5,2) applied to $G\mino{1}e$,
using (\ref{eq:implied-by-forbidden-e-f}))}   \\
& = &
\sigma_{G,\omega}(e) ~~~~~~~~~~~~~~~~~~~~~~~~~~~
\hbox{(by Theorem \ref{thm:minors-perms}(2,1))} .
\end{eqnarray*}

On the other hand,
if $e\not=\sigma_{G,\omega}(f)$ then 
\begin{eqnarray*}
\sigma_{G\mino{\omega}f\mino{1}e,\omega}(\sigma_{G,\omega}^{-1}(e))
& = &
\sigma_{G\mino{\omega}f\mino{1}e,\omega}(\sigma_{G\mino{\omega}f,\omega}^{-1}(e))
~~~~~~~~~~
\hbox{(by Lemma \ref{lemma:minors-perms-inv}(b))}   \\
& = &
\sigma_{G\mino{\omega}f,\omega}(e)
~~~~~~~~~~~
\hbox{(by Theorem \ref{thm:minors-perms}(2,1))
applied to $G\mino{\omega}f$)}   \\
& = &
\sigma_{G,\omega}(e)
~~~~~~~~~~~
\hbox{(by Theorem \ref{thm:minors-perms}(5,2),
using (\ref{eq:implied-by-forbidden-e-f}))} .
\end{eqnarray*}
while if $e=\sigma_{G,\omega}(f)$ then 
$\sigma_{G,\omega}^{-1}(e)=f$ which is not in the domain of
$\sigma_{G\mino{\omega}f\mino{1}e,\omega}$
so this situation does not arise.

So, in any event,
$\sigma_{G\mino{1}e\mino{\omega}f,\omega}(g) =
\sigma_{G\mino{\omega}f\mino{1}e,\omega}(g)$ in this case.   \\

Case 2:
$g = \sigma_{G,\omega}^{-1}(f)$.

\begin{eqnarray*}
\sigma_{G\mino{1}e\mino{\omega}f,\omega}(\sigma_{G,\omega}^{-1}(f))
& = &
\sigma_{G\mino{1}e\mino{\omega}f,\omega}(\sigma_{G\mino{1}e,\omega}^{-1}(f))
~~~~~~~
\hbox{(by Lemma \ref{lemma:minors-perms-inv}(a),
using $f\not=\sigma_{G,\omega}(e)$)}   \\
& = &
\sigma_{G\mino{1}e,\omega}(f)
~~~~~~~
\hbox{(by Theorem \ref{thm:minors-perms}(2,2))}   \\
& = &
\left\{
\begin{array}{ll}
\sigma_{G,\omega}(f)
~~~~~~~~
\hbox{(by Theorem \ref{thm:minors-perms}(5,1))},   &
\hbox{if $f\not=\sigma_{G,\omega}^{-1}(e)$,}   \\
\sigma_{G\mino{1}e,\omega}(\sigma_{G,\omega}^{-1}(e)),   &
\hbox{if $f=\sigma_{G,\omega}^{-1}(e)$,}\caseSpace   \\
\end{array}
\right.   \\
& = &
\left\{
\begin{array}{ll}
\sigma_{G,\omega}(f),  &
\hbox{if $f\not=\sigma_{G,\omega}^{-1}(e)$,}   \\
\sigma_{G,\omega}(e)),   &
\hbox{if $f=\sigma_{G,\omega}^{-1}(e)$.}\caseSpace   \\
\end{array}
\right.   \\
\end{eqnarray*}

Now consider
$\sigma_{G\mino{\omega}f\mino{1}e,\omega}(\sigma_{G,\omega}^{-1}(f))$.

If $e\not=\sigma_{G,\omega}(f)$, we have
\begin{eqnarray*}
\sigma_{G\mino{\omega}f\mino{1}e,\omega}(\sigma_{G,\omega}^{-1}(f))
& = &
\sigma_{G\mino{\omega}f,\omega}(\sigma_{G,\omega}^{-1}(f))
~~~~~~~
\hbox{(by Theorem \ref{thm:minors-perms}(5,1) and
(\ref{eq:implied-by-e-left-f}))}   \\
& = &
\sigma_{G,\omega}(f)
~~~~~~~
\hbox{(by Theorem \ref{thm:minors-perms}(2,2)).}
\end{eqnarray*}
If $e=\sigma_{G,\omega}(f)$ then
\begin{eqnarray*}
\sigma_{G\mino{\omega}f\mino{1}e,\omega}(\sigma_{G,\omega}^{-1}(f))
& = &
\sigma_{G\mino{\omega}f\mino{1}e,\omega}(\sigma_{G\mino{\omega}f,\omega}^{-1}(e))
~~~~~~~
\hbox{(by Lemma \ref{lemma:minors-perms-inv}(b) and
(\ref{eq:implied-by-e-left-f}))}   \\
& = &
\sigma_{G\mino{\omega}f,\omega}(e)
~~~~~~~
\hbox{(by Theorem \ref{thm:minors-perms}(2,1))}   \\
& = &
\sigma_{G,\omega}(e)
~~~~~~~
\hbox{(by Theorem \ref{thm:minors-perms}(5,2),
using $e\not=\sigma_{G,\omega}^{-1}(f)$).}
\end{eqnarray*}

Case 3:
$g = \sigma_{G\mino{1}e,\omega}^{-1}(f)$

\begin{eqnarray*}
\sigma_{G\mino{1}e\mino{\omega}f,\omega}(\sigma_{G\mino{1}e,\omega}^{-1}(f))
& = &
\sigma_{G\mino{1}e,\omega}(f)
~~~~~~~~
\hbox{(by Theorem \ref{thm:minors-perms}(2,2))}   \\
& = &
\left\{
\begin{array}{ll}
\sigma_{G,\omega}(f)
~~~~~~~
\hbox{(by Theorem \ref{thm:minors-perms}(5,1)),}  &
\hbox{if $e\not=\sigma_{G,\omega}(f)$,}  \\
\sigma_{G\mino{1}e,\omega}(\sigma_{G,\omega}^{-1}(e)),  &
\hbox{if $e=\sigma_{G,\omega}(f)$,}\caseSpace  \\
\end{array}
\right.   \\
& = &
\left\{
\begin{array}{ll}
\sigma_{G,\omega}(f),  &
\hbox{if $e\not=\sigma_{G,\omega}(f)$,}  \\
\sigma_{G,\omega}(e),  &
\hbox{if $e=\sigma_{G,\omega}(f)$}
~~~~~~~
\hbox{(by Theorem \ref{thm:minors-perms}(2,1)).}\caseSpace  \\
\end{array}
\right.
\end{eqnarray*}

If $e\not=\sigma_{G,\omega}(f)$ then
\begin{eqnarray*}
\sigma_{G\mino{\omega}f\mino{1}e,\omega}(\sigma_{G\mino{1}e,\omega}^{-1}(f))
& = &
\sigma_{G\mino{\omega}f\mino{1}e,\omega}(\sigma_{G,\omega}^{-1}(f))
~~~~~~~
\hbox{(by Lemma \ref{lemma:minors-perms-inv}(a))}   \\
& = &
\sigma_{G\mino{\omega}f,\omega}(\sigma_{G,\omega}^{-1}(f))
~~~~~~~
\hbox{(by Theorem \ref{thm:minors-perms}(5,1)) and (\ref{eq:implied-by-e-left-f}))}   \\
& = &
\sigma_{G,\omega}(f) ~~~~~~~~
\hbox{(by Theorem \ref{thm:minors-perms}(2,2))}.
\end{eqnarray*}

If $e=\sigma_{G,\omega}(f)$ then
\begin{eqnarray*}
\sigma_{G\mino{\omega}f\mino{1}e,\omega}(\sigma_{G\mino{1}e,\omega}^{-1}(f))
& = &
\sigma_{G\mino{\omega}f\mino{1}e,\omega}(\sigma_{G,\omega}^{-1}(f))
~~~~~~~~
\hbox{(by Lemma \ref{lemma:minors-perms-inv}(a), using $f\not=\sigma_{G,\omega}(e)$)}   \\
& = &
\sigma_{G\mino{\omega}f\mino{1}e,\omega}(\sigma_{G\mino{\omega}f,\omega}^{-1}(e))
~~~~~~~~
\hbox{(by Lemma \ref{lemma:minors-perms-inv}(b) and (\ref{eq:implied-by-e-left-f}))}   \\
& = &
\sigma_{G\mino{\omega}f,\omega}(e)
~~~~~~~~
\hbox{(by Theorem \ref{thm:minors-perms}(2,1))}   \\
& = &
\sigma_{G,\omega}(e)
~~~~~~~~
\hbox{(by Theorem \ref{thm:minors-perms}(5,2),
using $e\not=\sigma_{G,\omega}^{-1}(f)$)}.
\end{eqnarray*}

Case 4:
$g = \sigma_{G\mino{\omega}f,\omega}^{-1}(e)$

This case can be proved in a manner similar to the previous cases.
But in fact this is not necessary, since we have shown the permutations
$\sigma_{G\mino{1}e\mino{\omega}f,\omega}$ and
$\sigma_{G\mino{\omega}f\mino{1}e,\omega}$ agree on every element of
their common domain except one, so they must agree on this last element too.
\eopf

\begin{lemma}
\label{lemma:reductions-commute-sigma-1}
If $f\not=\sigma_{G,\omega}(e)$ then
\[
\sigma_{G\mino{1}e\mino{\omega}f,1} =
\sigma_{G\mino{\omega}f\mino{1}e,1}.
\]
\end{lemma}

\pf
This proof is more complicated again than that of
Lemma \ref{lemma:reductions-commute-sigma-omega}.

We may suppose that neither $e$ nor $f$ is
a triloop, since we have already established commutativity in such cases,
in Lemma \ref{lemma:reductions-commute-loops}.  So
$\sigma_{G,\mu}(e)\not=e$ and $\sigma_{G,\mu}(f)\not=f$ ,
for $\mu\in\{1,\omega,\omega^2\}$.

We prove that
\begin{equation}
\label{eq:sigmaG1eof1invg-sigmaGof1einvg}
\sigma_{G\mino{1}e\mino{\omega}f,1}^{-1}(g) =
\sigma_{G\mino{\omega}f\mino{1}e,1}^{-1}(g)
\end{equation}
for all $g\in E(G)\setminus\{e,f\}$.

Observe that
\begin{eqnarray*}
\sigma_{G\mino{1}e\mino{\omega}f,\omega^2}(g)
& = &
\sigma_{G\mino{1}e,\omega^2}(g) ~~~~~~~~
\hbox{if
$g\not\in\{\sigma_{G\mino{1}e,\omega^2}^{-1}(f),\sigma_{G\mino{1}e,1}(f)\}$,
by Theorem \ref{thm:minors-perms}(5,2)}   \\
& = &
\sigma_{G,\omega^2}(g) ~~~~~~~~
\hbox{if $g\not=\sigma_{G,\omega^2}^{-1}(e)$,
by Theorem \ref{thm:minors-perms}(5,1)}   \\
& = &
\sigma_{G\mino{\omega}f,\omega^2}(g) ~~~~~~~~
\hbox{if
$g\not\in\{\sigma_{G,\omega^2}^{-1}(f),\sigma_{G,1}(f)\}$,
by Theorem \ref{thm:minors-perms}(5,2)}   \\
& = &
\sigma_{G\mino{\omega}f\mino{1}e,\omega^2}(g) ~~~~~~~~
\hbox{if $g\not=\sigma_{G\mino{\omega}f,\omega^2}^{-1}(e)$,
by Theorem \ref{thm:minors-perms}(5,1)} .
\end{eqnarray*}

It follows that if $g\not\in\{\sigma_{G,1}(f),\sigma_{G,\omega^2}^{-1}(f),\sigma_{G,\omega^2}^{-1}(e),\sigma_{G\mino{\omega}f,\omega^2}^{-1}(e),\sigma_{G\mino{1}e,\omega^2}^{-1}(f),\sigma_{G\mino{1}e,1}(f)\}$ then
\begin{eqnarray*}
\sigma_{G\mino{1}e\mino{\omega}f,1}^{-1}(g)
& = &
\sigma_{G\mino{1}e\mino{\omega}f,\omega}(\sigma_{G\mino{1}e\mino{\omega}f,\omega^2}(g))
~~~~~~~~ \hbox{(using $\sigma_{1}\circ\sigma_{\omega}\circ\sigma_{\omega^2}=\hbox{identity}$)}   \\
& = &
\sigma_{G\mino{1}e\mino{\omega}f,\omega}(\sigma_{G\mino{\omega}f\mino{1}e,\omega^2}(g))
~~~~~~~~ \hbox{(by the previous paragraph)}   \\
& = &
\sigma_{G\mino{\omega}f\mino{1}e,\omega}(\sigma_{G\mino{\omega}f\mino{1}e,\omega^2}(g))
~~~~~~~~ \hbox{(by Lemma \ref{lemma:reductions-commute-sigma-omega})}   \\
& = &
\sigma_{G\mino{\omega}f\mino{1}e,1}^{-1}(g)
~~~~~~~~ \hbox{(using $\sigma_{1}\circ\sigma_{\omega}\circ\sigma_{\omega^2}=\hbox{identity}$, again)} .
\end{eqnarray*}

We now consider in turn how to deal with the exceptional values of $g$,
apparently six in number.

Case 1:  $g = \sigma_{G,1}(f)$.

We must have $\sigma_{G,1}(f)\not=e$, else $g=e$ which is forbidden.

Theorem \ref{thm:minors-perms}(5,1) tells us that
\begin{equation}
\label{eq:G1e1f-eq-G1f}
\sigma_{G\mino{1}e,1}(f) = 
\sigma_{G,1}(f) ,
\end{equation}
since $f\not=\sigma_{G,1}^{-1}(e)$ (by the previous paragraph)
and $f\not=\sigma_{G,\omega}(e)$ (by hypothesis), so
Theorem \ref{thm:minors-perms}(1,1) and (4,1) do not apply.

We have
\begin{eqnarray*}
\sigma_{G\mino{1}e\mino{\omega}f,1}^{-1}(\sigma_{G,1}(f))
& = &
\sigma_{G\mino{1}e\mino{\omega}f,1}^{-1}(\sigma_{G\mino{1}e,1}(f))
~~~~~~~~ \hbox{(by (\ref{eq:G1e1f-eq-G1f}))}   \\
& = &
\sigma_{G\mino{1}e,1}^{-1}(f)
~~~~~~~~ \hbox{(by Lemma \ref{lemma:minors-perms-inv}(d), first case)}   \\
& = &
\left\{
\begin{array}{ll}
\sigma_{G,\omega}(e),  &  \hbox{if $f = \sigma_{G,1}(e)$},  \\
\sigma_{G,1}^{-1}(e),  &  \hbox{if $f = \sigma_{G,\omega^2}^{-1}(e)$},\caseSpace  \\
\sigma_{G,1}^{-1}(f),  &  \hbox{otherwise,}\caseSpace
\end{array}
\right.   \\
\end{eqnarray*}
by Lemma \ref{lemma:minors-perms-inv}(c).

Now consider $\sigma_{G\mino{\omega}f\mino{1}e,1}^{-1}(\sigma_{G,1}(f))$.

If $f = \sigma_{G,1}(e)$ then
$\sigma_{G\mino{\omega}f,1}(e) =
\sigma_{G\mino{\omega}f,1}(\sigma_{G,1}^{-1}(f)) =
\sigma_{G,1}(f)$,
with the second equality following from 
Theorem \ref{thm:minors-perms}(1,2).
This justifies the first step of the following.
\begin{eqnarray*}
\sigma_{G\mino{\omega}f\mino{1}e,1}^{-1}(\sigma_{G,1}(f))
& = &
\sigma_{G\mino{\omega}f\mino{1}e,1}^{-1}(\sigma_{G\mino{\omega}f,1}(e))   \\
& = &
\sigma_{G\mino{\omega}f,\omega}(e)  ~~~~~~~~
\hbox{(by Lemma \ref{lemma:minors-perms-inv}(c))}   \\
& = &
\sigma_{G,\omega}(e)  ~~~~~~~~
\hbox{(using our hypothesis, $e\not=\sigma_{G,\omega}^{-1}(f)$).}
\end{eqnarray*}

If $f = \sigma_{G,\omega^2}^{-1}(e)$, i.e.,
$e = \sigma_{G,\omega^2}(f)$,
then $\sigma_{G\mino{\omega}f,\omega^2}^{-1}(e) = \sigma_{G,1}(f)$,
by Lemma \ref{lemma:minors-perms-inv}(f) (second case).
This justifies the first step of the following.
\begin{eqnarray*}
\sigma_{G\mino{\omega}f\mino{1}e,1}^{-1}(\sigma_{G,1}(f))
& = &
\sigma_{G\mino{\omega}f\mino{1}e,1}^{-1}(\sigma_{G\mino{\omega}f,\omega^2}^{-1}(e))   \\
& = &
\sigma_{G\mino{\omega}f,1}^{-1}(e)  ~~~~~~~~
\hbox{(by Lemma \ref{lemma:minors-perms-inv}(c))}   \\
& = &
\sigma_{G,1}^{-1}(e)  ~~~~~~~~
\hbox{(by Lemma \ref{lemma:minors-perms-inv}(d),
using $e\not=\sigma_{G,1}(f)$).}
\end{eqnarray*}

Suppose, then, that $f \not= \sigma_{G,1}(e)$ and
$f \not= \sigma_{G,\omega^2}^{-1}(e)$.

From $e \not= \sigma_{G,1}^{-1}(f)$
we deduce that $\sigma_{G\mino{\omega}f,1}(e) = \sigma_{G,1}(e)$,
by Theorem \ref{thm:minors-perms}(5,2).  Also, since $e\not=f$ and
$\sigma_{G,1}$ is a bijection, we have $\sigma_{G,1}(e)\not=\sigma_{G,1}(f)$.
So $\sigma_{G\mino{\omega}f,1}(e) \not= \sigma_{G,1}(f)$.

From $e \not= \sigma_{G,\omega^2}(f)$, and our hypothesis
$e\not=\sigma_{G,\omega}^{-1}(f)$, we deduce from
Lemma \ref{lemma:minors-perms-inv}(f) that
$\sigma_{G\mino{\omega}f,\omega^2}^{-1}(e) = \sigma_{G,\omega^2}^{-1}(e)$.
Our hypothesis $e\not=\sigma_{G,\omega}^{-1}(f)$ implies
$e\not=\sigma_{G,\omega^2}(\sigma_{G,1}(f))$, which in turn implies
$\sigma_{G,\omega^2}^{-1}(e) \not= \sigma_{G,1}(f)$.  Combining the
conclusions of the two previous sentences, we obtain
$\sigma_{G\mino{\omega}f,\omega^2}^{-1}(e) \not= \sigma_{G,1}(f)$.

The conclusions of the previous two paragraphs, together with
Lemma \ref{lemma:minors-perms-inv}(c), justify the first step in the following.
\begin{eqnarray*}
\sigma_{G\mino{\omega}f\mino{1}e,1}^{-1}(\sigma_{G,1}(f))
& = &
\sigma_{G\mino{\omega}f,1}^{-1}(\sigma_{G,1}(f))   \\
& = &
\sigma_{G,1}^{-1}(f)  ~~~~~~~~
\hbox{(by Lemma \ref{lemma:minors-perms-inv}(d)).}
\end{eqnarray*}

We have shown, then, that
$\sigma_{G\mino{1}e\mino{\omega}f,1}^{-1}$ and
$\sigma_{G\mino{\omega}f\mino{1}e,1}^{-1}$ agree on
$g = \sigma_{G,1}(f) $, in all circumstances.
This deals with the first of our exceptional values of $g$.

Case 2:  $g = \sigma_{G,\omega^2}^{-1}(f)$.

We must have $\sigma_{G,\omega^2}^{-1}(f)\not=e$, else $g=e$ which is forbidden.

Firstly, observe that
$\sigma_{G\mino{1}e,1}(f)\in\{\sigma_{G,\omega^2}^{-1}(e),\sigma_{G,1}(f)\}$,
by Theorem \ref{thm:minors-perms}(1,1),(5,1), using the hypothesis
$f\not=\sigma_{G,\omega}(e)$.
Now, $\sigma_{G,\omega^2}^{-1}(e) \not= \sigma_{G,\omega^2}^{-1}(f)$, 
since $e\not=f$.
Furthermore, $\sigma_{G,1}(f) \not= \sigma_{G,\omega^2}^{-1}(f)$,
since if $\sigma_{G,1}(f) = \sigma_{G,\omega^2}^{-1}(f)$
then $f = \sigma_{G,\omega^2}( \sigma_{G,1}(f) ) = \sigma_{G,\omega}^{-1}(f)$,
which means that $f$ is a triloop, which we excluded at the start.
So, whatever its value, we have
$\sigma_{G\mino{1}e,1}(f) \not= \sigma_{G,\omega^2}^{-1}(f)$. 

Secondly, observe that
$\sigma_{G\mino{1}e,\omega^2}^{-1}(f) = \sigma_{G,\omega^2}^{-1}(f)$,
by Lemma \ref{lemma:minors-perms-inv}(e), using $f\not=\sigma_{G,\omega^2}(e)$
(see the start of this Case).

The conclusions of these two previous paragraphs justify the first two steps
in the following.
\begin{eqnarray*}
\sigma_{G\mino{1}e\mino{\omega}f,1}^{-1}(\sigma_{G,\omega^2}^{-1}(f))
& = &
\sigma_{G\mino{1}e,1}^{-1}(\sigma_{G,\omega^2}^{-1}(f))
~~~~~~~~
\hbox{(by Lemma \ref{lemma:minors-perms-inv}(d))}  \\
& = &
\sigma_{G\mino{1}e,1}^{-1}(\sigma_{G\mino{1}e,\omega^2}^{-1}(f))   \\
& = &
\sigma_{G\mino{1}e,\omega}(f)   \\
& = &
\left\{
\begin{array}{ll}
\sigma_{G,\omega}(e),  &  \hbox{if $f=\sigma_{G,\omega}^{-1}(e)$,}  \\
\sigma_{G,\omega}(f),  &  \hbox{otherwise,}\caseSpace
\end{array}
\right.
\end{eqnarray*}
by Theorem \ref{thm:minors-perms}(2,1),(5,2).

Now consider
$\sigma_{G\mino{\omega}f\mino{1}e,1}^{-1}(\sigma_{G,\omega^2}^{-1}(f))$.

If $f=\sigma_{G,\omega}^{-1}(e)$, then
$\sigma_{G,\omega^2}^{-1}(f)=\sigma_{G,1}(e)$.
Also, $e\not=\sigma_{G,1}^{-1}(f)$, since
if $e=\sigma_{G,1}^{-1}(f)$ then
$f=\sigma_{G,\omega}^{-1}( \sigma_{G,1}^{-1}(f) ) =
\sigma_{G,\omega^2}(f)$, so that $f$ is a triloop, which we have excluded.
So $\sigma_{G\mino{\omega}f,1}(e)=\sigma_{G,1}(e)=\sigma_{G,\omega^2}^{-1}(f)$,
with the first equality holding by Theorem \ref{thm:minors-perms}(5,2).
This justifies the first step in the following.
\begin{eqnarray*}
\sigma_{G\mino{\omega}f\mino{1}e,1}^{-1}(\sigma_{G,\omega^2}^{-1}(f))
& = &
\sigma_{G\mino{\omega}f\mino{1}e,1}^{-1}(\sigma_{G\mino{\omega}f,1}(e))  \\
& = &
\sigma_{G\mino{\omega}f,\omega}(e)  ~~~~~~~~
\hbox{(by Lemma \ref{lemma:minors-perms-inv}(c))}   \\
& = &
\sigma_{G,\omega}(e)  ~~~~~~~~
\hbox{(by Theorem \ref{thm:minors-perms}(5,2), since $e\not=\sigma_{G,\omega}^{-1}(f)$ by hypothesis).}
\end{eqnarray*}

If $f\not=\sigma_{G,\omega}^{-1}(e)$, then
$\sigma_{G,\omega^2}^{-1}(f)\not=\sigma_{G,1}(e)$.
Also, $\sigma_{G,\omega^2}^{-1}(f)\not=\sigma_{G,1}(f)$, else $f$ is a triloop,
as we saw early in this Case.
So $\sigma_{G,\omega^2}^{-1}(f)\not\in\{\sigma_{G,1}(e),\sigma_{G,1}(f)\}$.
But $\sigma_{G\mino{\omega}f,1}(e)\in\{\sigma_{G,1}(e),\sigma_{G,1}(f)\}$,
by Theorem \ref{thm:minors-perms}(1,2),(5,2). 
So $\sigma_{G,\omega^2}^{-1}(f)\not=\sigma_{G\mino{\omega}f,1}(e)$.

Since $e\not=\sigma_{G,\omega}^{-1}(f)$ by hypothesis,
$\sigma_{G\mino{\omega}f,\omega^2}^{-1}(e)\in
\{\sigma_{G,1}(f),\sigma_{G,\omega^2}^{-1}(e)\}$.
Now, as we have seen, $\sigma_{G,\omega^2}^{-1}(f)\not=\sigma_{G,1}(f)$,
else $f$ is a triloop; also,
$\sigma_{G,\omega^2}^{-1}(f)\not=\sigma_{G,\omega^2}^{-1}(e)$, since $e\not=f$.
So $\sigma_{G,\omega^2}^{-1}(f)\not=\sigma_{G\mino{\omega}f,\omega^2}^{-1}(e)$.

The conclusions of the previous two paragraphs, together with
Lemma \ref{lemma:minors-perms-inv}(c), justify the first step of the following.
\begin{eqnarray*}
\sigma_{G\mino{\omega}f\mino{1}e,1}^{-1}(\sigma_{G,\omega^2}^{-1}(f))
& = &
\sigma_{G\mino{\omega}f,1}^{-1}(\sigma_{G,\omega^2}^{-1}(f))   \\
& = &
\sigma_{G,1}^{-1}(\sigma_{G,\omega^2}^{-1}(f))  ~~~~~~~~
\hbox{(by Lemma \ref{lemma:minors-perms-inv}(d), using
$\sigma_{G,\omega^2}^{-1}(f)\not=\sigma_{G,1}(f)$)}   \\
& = &
\sigma_{G,\omega}(f) .
\end{eqnarray*}

So $\sigma_{G\mino{1}e\mino{\omega}f,1}^{-1}$ and
$\sigma_{G\mino{\omega}f\mino{1}e,1}^{-1}$ agree on
$g = \sigma_{G,\omega^2}^{-1}(f) $, in all circumstances.
This deals with the second of our exceptional values of $g$.   \\

Case 3:  $g = \sigma_{G,\omega^2}^{-1}(e)$.

We must have $\sigma_{G,\omega^2}^{-1}(e)\not=f$, else $g=f$ which is forbidden.

Observe that
$\sigma_{G\mino{1}e,1}(\sigma_{G,1}^{-1}(e)) = \sigma_{G,\omega^2}^{-1}(e)$,
by Theorem \ref{thm:minors-perms}(1,1).

So, if $f=\sigma_{G,1}^{-1}(e)$ then
\begin{eqnarray*}
\sigma_{G\mino{1}e\mino{\omega}f,1}^{-1}( \sigma_{G,\omega^2}^{-1}(e) )
& = &
\sigma_{G\mino{1}e\mino{\omega}f,1}^{-1}(
\sigma_{G\mino{1}e,1}(f) )   \\
& = &
\sigma_{G\mino{1}e,1}^{-1}(f)  ~~~~~~~~
\hbox{(by Lemma \ref{lemma:minors-perms-inv}(d))}   \\
& = &
\left\{
\begin{array}{ll}
\sigma_{G,\omega}(e),  &  \hbox{if also $f=\sigma_{G,1}(e)$,}  \\
\sigma_{G,1}^{-1}(f),  &  \hbox{if also $f\not=\sigma_{G,1}(e)$,}\caseSpace
\end{array}
\right.
\end{eqnarray*}
by Lemma \ref{lemma:minors-perms-inv}(c) with $f\not=\sigma_{G,\omega^2}^{-1}(e)$.

On the other hand, if $f\not=\sigma_{G,1}^{-1}(e)$ then
$\sigma_{G\mino{1}e,1}(f) \not= \sigma_{G\mino{1}e,1}( \sigma_{G,1}^{-1}(e) )$,
since $\sigma_{G\mino{1}e,1}$ is a bijection.  So
$\sigma_{G,\omega^2}^{-1}(e) \not= \sigma_{G\mino{1}e,1}(f)$.
We therefore have
\begin{eqnarray*}
\sigma_{G\mino{1}e\mino{\omega}f,1}^{-1}( \sigma_{G,\omega^2}^{-1}(e) )
& = &
\sigma_{G\mino{1}e,1}^{-1}( \sigma_{G,\omega^2}^{-1}(e) ) ~~~~~~~~
\hbox{(by Lemma \ref{lemma:minors-perms-inv}(d))}   \\
& = &
\sigma_{G\mino{1}e,1}^{-1}(
\sigma_{G\mino{1}e,1}(\sigma_{G,1}^{-1}(e)) )   \\
& = &
\sigma_{G\mino{1}e,1}^{-1}(
\sigma_{G\mino{1}e,1}(\sigma_{G,1}^{-1}(e)) )  ~~~~~~~~
\hbox{(by Lemma \ref{lemma:minors-perms-inv}(d))}   \\
& = &
\sigma_{G,1}^{-1}(e) .
\end{eqnarray*}

So, in summary,
\begin{eqnarray*}
\sigma_{G\mino{1}e\mino{\omega}f,1}^{-1}( \sigma_{G,\omega^2}^{-1}(e) )
& = &
\left\{
\begin{array}{ll}
\sigma_{G,\omega}(e),  &  \hbox{if $f=\sigma_{G,1}^{-1}(e)$
and $f=\sigma_{G,1}(e)$,}  \\
\sigma_{G,1}^{-1}(f),  &  \hbox{if $f=\sigma_{G,1}^{-1}(e)$
and $f\not=\sigma_{G,1}(e)$,}  \\
\sigma_{G,1}^{-1}(e),  &  \hbox{if $f\not=\sigma_{G,1}^{-1}(e)$}.
\end{array}
\right.
\end{eqnarray*}

Now consider
$\sigma_{G\mino{\omega}f\mino{1}e,1}^{-1}(\sigma_{G,\omega^2}^{-1}(e))$.

Since $e\not=\sigma_{G,\omega}^{-1}(f)$, by hypothesis, and
$e\not=\sigma_{G,\omega^2}(f)$ (see start of this Case),
Lemma \ref{lemma:minors-perms-inv}(f) gives
$\sigma_{G,\omega^2}^{-1}(e)=\sigma_{G\mino{\omega}f,\omega^2}^{-1}(e)$.

Consider, for a moment, the circumstances under which
$\sigma_{G\mino{\omega}f,1}^{-1}(e)=e$.  Lemma \ref{lemma:minors-perms-inv}(d)
tells us that
\[
\sigma_{G\mino{\omega}f,1}^{-1}(e)=
\left\{
\begin{array}{ll}
\sigma_{G,1}^{-1}(f),  &  \hbox{if $e=\sigma_{G,1}(f)$},  \\
\sigma_{G,1}^{-1}(e),  &  \hbox{if $e\not=\sigma_{G,1}(f)$}.\caseSpace
\end{array}
\right.
\]
If $e\not=\sigma_{G,1}(f)$ then $\sigma_{G\mino{\omega}f,1}^{-1}(e)\not=e$,
since otherwise $e=\sigma_{G,1}^{-1}(e)$, so that $e$
is a triloop, which we have excluded.
Also, if $e\not=\sigma_{G,1}^{-1}(f)$ then $e$ cannot equal either of
the two possible expressions just given for
$\sigma_{G\mino{\omega}f,1}^{-1}(e)$ (using the triloop exclusion, again,
for the second of these).  So, again,
$\sigma_{G\mino{\omega}f,1}^{-1}(e)\not=e$.  On the other hand,
if $e=\sigma_{G,1}(f)$ and $e=\sigma_{G,1}^{-1}(f)$ then the first case
above gives
$\sigma_{G\mino{\omega}f,1}^{-1}(e) = \sigma_{G,1}^{-1}(f) = e$.
In this situation, in applying Theorem \ref{thm:minors-perms} to
$G\mino{\omega}f$, we cannot use case (1,1), since that would require
$\sigma_{G\mino{\omega}f,1}^{-1}(e)\not=e$.  Similarly,
we cannot use the inverse
of case (1,1) to find
$\sigma_{G\mino{\omega}f\mino{1}e,1}^{-1}(
\sigma_{G\mino{\omega}f,\omega^2}^{-1}(e))$; instead, we must use case (5,1).

If $e=\sigma_{G,1}(f)$ and $e=\sigma_{G,1}^{-1}(f)$, then, we have
\begin{eqnarray*}
\sigma_{G\mino{\omega}f\mino{1}e,1}^{-1}(\sigma_{G,\omega^2}^{-1}(e))
& = &
\sigma_{G\mino{\omega}f\mino{1}e,1}^{-1}(
\sigma_{G\mino{\omega}f,\omega^2}^{-1}(e))   \\
& = &
\sigma_{G\mino{\omega}f,1}^{-1}(
\sigma_{G\mino{\omega}f,\omega^2}^{-1}(e))  ~~~~~~~~
\hbox{(by Theorem \ref{thm:minors-perms}(5,1))}   \\
& = &
\sigma_{G\mino{\omega}f,\omega}(e)   \\
& = &
\sigma_{G,\omega}(e)  ~~~~~~~~
\hbox{(by Theorem \ref{thm:minors-perms}(5,2), using our hypothesis
$e\not=\sigma_{G,\omega}^{-1}(f)$)} .
\end{eqnarray*}

Otherwise, we have
\begin{eqnarray*}
\sigma_{G\mino{\omega}f\mino{1}e,1}^{-1}(\sigma_{G,\omega^2}^{-1}(e))
& = &
\sigma_{G\mino{\omega}f\mino{1}e,1}^{-1}(
\sigma_{G\mino{\omega}f,\omega^2}^{-1}(e))   \\
& = &
\sigma_{G\mino{\omega}f,1}^{-1}(e)  ~~~~~~~~
\hbox{(by Lemma \ref{lemma:minors-perms-inv}(c))}   \\
& = &
\left\{
\begin{array}{ll}
\sigma_{G,1}^{-1}(f),  &  \hbox{if $e=\sigma_{G,1}(f)$ (and so $e\not=\sigma_{G,1}^{-1}(f)$ too,}  \\
&  \hbox{\ \ \ \ else we are in the previous paragraph),}  \\
\sigma_{G,1}^{-1}(e),  &  \hbox{if $e\not=\sigma_{G,1}(f)$},
\end{array}
\right.
\end{eqnarray*}
by Lemma \ref{lemma:minors-perms-inv}(d).

So $\sigma_{G\mino{1}e\mino{\omega}f,1}^{-1}$ and
$\sigma_{G\mino{\omega}f\mino{1}e,1}^{-1}$ agree on
$g = \sigma_{G,\omega^2}^{-1}(e)$, in all circumstances.
This deals with the third of our exceptional values of $g$.

Cases 4--6:
$g\in\{\sigma_{G\mino{\omega}f,\omega^2}^{-1}(e),\sigma_{G\mino{1}e,\omega^2}^{-1}(f),\sigma_{G\mino{1}e,1}(f)\}$.

Theorem \ref{thm:minors-perms} and Lemma \ref{lemma:minors-perms-inv} tell us that
\begin{eqnarray*}
\sigma_{G\mino{\omega}f,\omega^2}^{-1}(e)
& \in &
\{ \sigma_{G,1}(f), \sigma_{G,\omega^2}^{-1}(e) \} ,   \\
\sigma_{G\mino{1}e,\omega^2}^{-1}(f)
& \in &
\{ \sigma_{G,\omega^2}^{-1}(e), \sigma_{G,\omega^2}^{-1}(f) \} ,   \\
\sigma_{G\mino{1}e,1}(f)
& \in &
\{ \sigma_{G,1}(f), \sigma_{G,\omega^2}^{-1}(e) \} .
\end{eqnarray*}

So these are not really new cases at all; they each take us back
into one of Cases 1--3.

This completes our proof of (\ref{eq:sigmaG1eof1invg-sigmaGof1einvg}),
and hence of the Lemma.
\eopf

\begin{thm}
If $f\not=\sigma_{G,\omega}(e)$ then
\[
G\mino{1}e\mino{\omega}f = G\mino{\omega}f\mino{1}e.
\]
\end{thm}

\pf
In view of Lemmas \ref{lemma:reductions-commute-sigma-omega}
and \ref{lemma:reductions-commute-sigma-1}, we know that
$\sigma_{G\mino{1}e\mino{\omega}f,\mu}=\sigma_{G\mino{\omega}f\mino{1}e,\mu}$
for $\mu\in\{1,\omega\}$.  But then it follows for $\mu=\omega^2$ too,
since $\sigma_{\omega^2}=\sigma_{\omega}^{-1}\circ\sigma_{1}^{-1}$.
\eopf   \\

Triality gives the following two corollaries.

\begin{cor}
If $f\not=\sigma_{G,1}(e)$ then
\[
G\mino{\omega}e\mino{\omega^2}f = G\mino{\omega^2}f\mino{\omega}e.
\]
\eopf
\end{cor}

\begin{cor}
If $f\not=\sigma_{G,\omega^2}(e)$ then
\[
G\mino{\omega^2}e\mino{1}f = G\mino{1}f\mino{\omega^2}e.
\]
\eopf
\end{cor}

The results so far in this section
(together with the fact of non-commutativity in general
for the excluded cases for the previous three results)
give us a complete description of when
the $\mu$-reductions do, or do not, commute, in general.

But some interesting questions remain.
Given that the excluded (generally non-commutative) cases are so specific,
it is natural to ask for a characterisation of those alternating dimaps
for which all reductions always commute.

Consider $f=\sigma_{G,\omega}(e)$, illustrated in Figure \ref{fig:non-commut}.
In this case, $\bullet\mino{1}e$ and $\bullet\mino{\omega}f$ do
not commute in general, but we can still investigate when they do.

\begin{propn}
\label{propn:f-eq-sigma-G-omega-e-commut}
If $f=\sigma_{G,\omega}(e)$ then
$\bullet\mino{1}e$ and $\bullet\mino{\omega}f$ commute
if and only if at least one of $e$, $f$ is a triloop.
\end{propn}

\pf
If either $e$ or $f$ is a triloop, then they commute by 
Lemma \ref{lemma:reductions-commute-loops}.
Suppose then that neither $e$ nor $f$ is a triloop.
If $e$ and $f$ form an a-face of size 2, then it is routine to
show that these reductions do not commute unless the head of $f$
meets no other edge except $f$, but that would make $f$ a 1-loop.
If $e$ and $f$ do not form such an a-face, then the endpoints of
$e$ and $f$ --- three in number --- are all distinct.  The situation
is then exactly as in Figure \ref{fig:non-commut} (except that the
right-hand vertex might coincide with the tail of $f$ or the head of $e$,
but that is immaterial).  It is evident from the Figure that the only
way the reductions can commute in this case is if the head of $e$ has
in-degree 1 (i.e., if the edges shown in green do not exist), which
would make $e$ a 1-loop.
\eopf

\begin{thm}
\label{thm:2-reduction-commutativity}
Every pair of reductions on $G$ commutes
if and only if the set of triloops of $G$ includes at least one
of each pair of edges that are consecutive in any in-star, a-face
or c-face.
\end{thm}

\pf
Use Proposition \ref{propn:f-eq-sigma-G-omega-e-commut}
and triality.
\eopf   \\

% ***** relate to trimedial graph

We pause now to introduce a graph derived from $G$ which gives
an alternative way of framing Theorem \ref{thm:2-reduction-commutativity}.
%  and whose acyclic orientations turn out to be related to the ordering
%  of reductions.

The \textit{trimedial graph} $\tri(G)$ of the alternating dimap $G$
has vertex set $E(G)$ with two vertices of $\tri(G)$ being adjacent
if their corresponding edges in $G$
are consecutive in an a-face, a c-face, or an in-star
of $G$.  The trimedial graph is always undirected and 6-regular, and
may have loops and/or multiple edges.  Its 6-regularity implies that,
if it has no loops or multiple edges, then it
is nonplanar even if $G$ is plane (in contrast to the usual medial graph).

% **** trimedial graph comes with a natural embedding.  Do I mention that?

With this definition, we may rewrite
Theorem \ref{thm:2-reduction-commutativity}.

\begin{cor}
\label{cor:2-reduction-commutativity-trimedial}
Every pair of reductions on $G$ commutes
if and only if the set of triloops of $G$ form a vertex cover of $\tri(G)$.
\eopf
\end{cor}

So far, we have considered the usual kind of commutativity,
where the order in which two operations are applied does not matter.
We can also ask about stronger forms of commutativity.
If a set of $k$ reductions (each of the form $\bullet\mino{\mu}e$,
where each $\mu\in\{1,\omega,\omega^2\}$ and all the $e$ are distinct)
has the property that applying them
in any order always gives the same result, then we say that it is
$k$\textit{-commutative} on $G$.
We say that $G$ is
\textit{$k$-reduction-commutative} if every set
of $k$ reductions is $k$-commutative on $G$.
It is \textit{totally reduction-commutative} if
it is $k$-reduction-commutative for every $k$.

In this terminology, ordinary commutativity is 2-commutativity, in the sense that,
if two particular reductions $\bullet\mino{\mu}e$ and $\bullet\mino{\nu}f$ commute,
then the set $\{\bullet\mino{\mu}e,\bullet\mino{\nu}f\}$ is 2-commutative.
Theorem \ref{thm:2-reduction-commutativity}
characterises alternating dimaps that are 2-reduction-commutative.

While total reduction-commutativity implies
$k$-reduction-commutativity for any fixed $k$, which in turn implies
$l$-reduction-commutativity for any $l<k$, the converses do not hold.

Consider how taking minors affects these properties.

\begin{propn}
\label{propn:minors-totally-reduction-commutative}
If $G$ is totally reduction-commutative, then so is any minor of $G$.
\eopf
\end{propn}

By contrast, 2-reduction-commutativity is not in general preserved
by taking minors.  To see this, let $H$ be any alternating dimap with no triloops, and form $G$
from it by inserting an $\omega^2$-loop at each vertex of each anticlockwise face and an
$\omega$-loop at each vertex of each clockwise face.  Then $H$ is a minor of $G$, yet
Theorem \ref{thm:2-reduction-commutativity} tells us that $G$ is 2-reduction-commutative
yet $H$ is not.

We now characterise alternating dimaps that are totally reduction-commutative.

A \textit{1-circuit} is an alternating dimap consisting of a single directed
circuit, in which every edge is a 1-loop.  An $\omega$\textit{-circuit}
(respectively, $\omega^2$\textit{-circuit}) consists of a single
vertex together with a number of $\omega$-loops (resp., $\omega^2$-loops)
at it.  A \textit{tricircuit} is an alternating dimap that can be
constructed from a 1-circuit, an $\omega$-circuit and an
$\omega^2$-circuit (any of which may have no edges), taking a single
vertex in each, and identifying these three vertices in the natural way.
This is done
so as to preserve the alternating dimap property, and will entail
having the $\omega$-circuit and $\omega^2$-circuit on opposite sides
of the 1-circuit.

\begin{thm}
\label{thm:totally-reduction-commut-tricircuits}
An alternating dimap $G$ is totally reduction-commutative if and only if
each of its components is a tricircuit.
\end{thm}

\pf
Suppose $G$ is totally reduction-commutative.  Then it is certainly
2-reduction-commutative, so by Theorem \ref{thm:2-reduction-commutativity}
the set of triloops of $G$ includes at least one
of each pair of edges that are consecutive in any in-star, a-face
or c-face.
%  Suppose some a-face $A$ has two edges $e$ and $f$ that are not loops.
%  Let $g_1, \ldots, g_l$ be the sequence of edges from $e$ to $f$ in
%  anticlockwise order around $A$.  Assume that $e$ and $f$ are chosen
%  from $A$ so as to minimise $l$.  Then all the $g_i$ are all loops,
%  otherwise the minimality of $l$ is violated.  The alternating dimap
%  $G'=G\mino{*}g_1\cdots\mino{*}g_l$ has $f=\sigma_{G',\omega}^{-1}(e)$,
%  so by Proposition \ref{propn:f-eq-sigma-G-omega-e-commut} at least one
%  of $e$ and $f$ is a triloop.  Since it is not a loop (by choice of $e$
%  and $f$), it must be a 1-loop.  It follows that each a-face has at most
%  one edge that is not a triloop.  Similarly, by triality, each c-face
%  and each in-star has at most one edge that is not a triloop.

Consider those edges of $G$ which     % are not triloops and
have distinct endpoints (i.e., the non-loops).

%  This para only assumes  e, f  are non-loops.  But the assumed
%  configuration (both pointing into  v)  means not a 1-loop either,
%  hence not a triloop.
Suppose two non-loop
edges $e$ and $f$ share an endpoint $v$, so $e,f\in I(v)$.
Since $e$ and $f$ are not loops, they do not come out of $v$.
The number of half-edges going out of $v$ must be two greater than the
number of half-edges other than $e$ and $f$ going into $v$.
So there must be two half-edges going out of $v$
that do not match (i.e., are not part of the same edge as) any half-edge
going into $v$.  Let $g$ be an edge to which one of these half-edges belongs.
Without loss of generality, suppose that $e,g,f$ occur in that order,
going clockwise around $v$.  Let the sequence of edges of $I(v)$
which are between $g$ and $e$ going anticlockwise be $h_1,\ldots,h_a$,
and let the sequence of edges of $I(v)$ which are between $g$ and $f$
going clockwise be $i_1,\ldots,i_c$.  Then the alternating dimap
$G':=G\mino{\omega^2}(h_1,\ldots,h_a)\mino{\omega}(i_1,\ldots,i_c)$
is left with $e,g,f$ intact, still in this same order around $v$,
and with no edges intervening between them any more.  Then $e$ and
$f$ are consecutive (clockwise) in the in-star at $v$ in $G'$.
By Theorem \ref{thm:2-reduction-commutativity}, this implies
non-commutativity of some reductions on $G'$, which in turn implies
that $G$ is not totally reduction-commutative.

Now suppose two non-triloop non-loops $e$ and $f$ are head-to-tail: say, with
$v$ = head of $e$ = tail of $f$.  Since $e$ is not a 1-loop, there must
be other edges at $v$.  If all of those edges lie between $e$ and $f$ going
clockwise, then $e$ and $f$ are consecutive around the clockwise face
containing $e$, so Theorem \ref{thm:2-reduction-commutativity} gives
non-commutativity of some reductions, so
$G$ is not totally reduction-commutative.  Similarly, if those extra
edges at $v$ all lie on the other side --- between $f$ and $e$ going
clockwise --- then, again, $G$ is not totally reduction-commutative.
So there are some edges on each side.  Let the edges of $I(v)$ between
$f$ and $e$ going anticlockwise be $h_1,\ldots,h_k$.  Then
$G':=G\mino{\omega^2}(h_1,\ldots,h_k)$ has $e$ and $f$ as consecutive
edges in the anticlockwise face containing $e$.  This gives some
non-commutative reductions in $G'$, so $G$ is not totally
reduction-commutative.

If non-triloop non-loops $e$ and $f$ belong to the same component
of $G$, then let $P$ be the shortest path, in the underlying undirected
graph, from one to the other.  (Note, $e,f\not\in E(P)$, and $P$ meets $e$ and $f$
only at the endpoints of $P$, by its minimality.)
%  This path must consist only of non-loops.
%  If there is any variation in the directions of the edges along $e,P,f$,
%  then there will be some two non-loops that share an endpoint, and as
%  we showed above, $G$ is then not totally reduction-commutative.
%  So all the edges of $e,P,f$ must go in the same
%  direction along $P$.  If any 
If all the edges of $P$ are contracted, to give $G\mino{1}E(P)$,
then we have $e$ and $f$ sharing an endpoint and we are in one of
the previous two paragraphs, so $G\mino{1}E(P)$ is not totally
reduction-commutative, so neither is $G$.

So each component of $G$ has at most one edge that is neither a triloop
nor a loop.

All the 1-loops in a component of $G$ must lie in a single
directed circuit in that component.  To see this, take any 1-loop $e=uv$.
It has a unique successor, which cannot be a loop or $e$ would not be
a 1-loop.  So it must either be a 1-loop or the sole edge which is neither
a triloop nor a loop.  Now let us go back the other way.  Consider the
edges in $I(u)$.  At least one of them must be a non-loop.
But if $I(u)$ has two non-loops, then both of them are not 1-loops,
and so this component has at least two edges that are neither a triloop
nor a loop, which is a contradiction.  So $I(u)$ has only one non-loop,
which must either be a 1-loop or the sole non-triloop non-loop.
We can follow 1-loops forwards and backwards in this way until we are
forced to stop.  This happens when we complete a circuit, which will either
be a circuit consisting entirely of 1-loops --- in which case it is
an entire component of $G$ --- or consisting of 1-loops except for the
sole non-triloop non-loop, which we call $f=wx$.
In the latter case, other edges may meet
the head $x$ of that special edge, but cannot meet any other vertex on
the circuit.  The other edges at $x$ must all be loops, since if any
is an outgoing non-loop then another must be an incoming non-loop which
is then not a 1-loop either, a contradiction with the uniqueness of $f$.
Furthermore, if any edge $g$ at $x$ is a proper 1-semiloop, then we can form a minor,
by reduction of any $\omega$-loops or $\omega^2$-loops that get in the way,
in which $f$ and $g$ form a configuration that allows non-commutativity.
So those other edges at $x$ must all be $\omega$-loops or $\omega^2$-loops.

This description of the component of $G$, as a circuit whose edges are
1-loops with possibly one exception, and with the head of that exception
holding $\omega$-loops and $\omega^2$-loops, identifies the component
as a tricircuit.  So every component of $G$ is in fact a tricircuit.

Conversely, if every component of $G$ is a tricircuit, then each component
has at most one edge that is not a triloop, so any two reductions
on $G$ commute, by Theorem \ref{thm:reductions-commute-loops}.
Therefore $G$ is totally reduction-commututative.
\eopf   \\

% **** above:  need to allow for the possibility that the 1-circuit
% **** consists of a single 1-semiloop.  And the argument about two
% **** edges pointing into the same vertex should, if possible, be
% **** extended to cope with the possibility that one or even both
% **** are 1-semiloops.

% **** check that Propn 14 and Thm 15 cope with 1-semiloops ...

So far, we have considered commutativity (or otherwise) with
respect to \textit{identity}: reductions commute if and only if
carrying them out in each possible order gives alternating dimaps
that are \textit{identical}.  We could also define commutativity
with respect to \textit{isomorphism}.  \\

\noindent\textbf{Problem}   \\
Characterise alternating dimaps $G$ for which,
for all $\mu_1,\mu_2\in\{1,\omega,\omega^2\}$ and all $e,f\in E(G)$,
\[
G\mino{\mu_1}e\mino{\mu_2}f \cong G\mino{\mu_2}f\mino{\mu_1}e .
\]

%  Raise the qn of commutativity in the "isomorphic" sense.  Give the example I came up with.
%  The thing about acyclic orientations in the trimedial graph.
%  Also, connection of trimedial gph, and vertex covers or stable sets in it,
%      with the alt dimaps where all reductions commute (in identity sense)

\section{Connections with binary functions}
\label{sec:bin-fns}

The relationship between triality and minor operations for
alternating dimaps is reminiscent of properties of binary
functions found by the author in \cite{farr2013}.  In this
section we briefly summarise those properties and compare
the relationships found there with those found here.
We will determine those alternating dimaps which can be
represented, in a certain faithful manner, by binary functions.

Let $E$ be a finite set, with $m=|E|$.
A \textit{binary function} with \textit{ground set} $E$
and \textit{dimension} $m$
is a function $f:2^E\rightarrow\mathbb{C}$
such that $f(\emptyset)=1$.  Equivalently, we regard it as
a $2^{m}$-element complex vector $\mathbf{f}$
whose elements are indexed
by the subsets of $E$ and whose first element (indexed by $\emptyset$) is 1.
(The restriction $f(\emptyset)=1$ was not imposed as part of the definition
in earlier work \cite{farr93,farr04,farr07a,farr07b}, but all scalar
multiples of a binary function are equivalent for our purposes, and
we have always been most interested in the cases where $f(\emptyset)\not=0$.)
We often represent a subset $X\subseteq E$ by its characteristic
vector $\mathbf{x}\in\{0,1\}^E$, with $x_e=1$ if $e\in X$ and $x_e=0$ otherwise.
Since $\mathbf{x}$ may be thought of as a binary string, it may also
be taken to be the binary representation of a number $x$ such that
$0\le x\le 2^{m}-1$.  With this notation, $f(X)$ may also be written
$f_{\mathbf{x}}$ or $f_x$.
In particular, $f(\emptyset)=f_{(0,\ldots,0)}=f_0=1$.
We write $0_k$ for the sequence of $k$ 0s,
and sometimes drop the subscript $k$ when it is clear from the context.

The definition was motivated by indicator functions of linear spaces
over GF(2), especially of cutset spaces of graphs: if $N$ is a matrix
over GF(2) whose columns are indexed by $E$ (such as the incidence
matrix of a graph, or the matrix representation of a binary matroid),
then the indicator function of the rowspace of $N$ takes value 1
on a set $X\subseteq E$ if the characteristic vector of $X$ belongs
to the rowspace of $N$, and takes value 0 otherwise.

If $f, g:2^E\rightarrow\mathbb{C}$ and there exists a constant
$c\in\mathbb{C}\setminus\{0\}$ such that $f(X) = c g(X)$
for all $X\subseteq E$, then we write $f\simeq g$.

Define
\[
M(\mu) := \frac{1}{2\sqrt{2}}
\left(
\begin{array}{cc}
\sqrt{2}+1+(\sqrt{2}-1)\mu  &  1-\mu   \\
1-\mu  &  \sqrt{2}-1+(\sqrt{2}+1)\mu
\end{array}
\right) .
\]
The $\mu$\textit{-transform} of $\mathbf{f}$, denoted by
$L^{[\mu]}\mathbf{f}$, is given by
\[
L^{[\mu]}\mathbf{f} := M(\mu)^{\otimes m}\mathbf{f} ,
\]
where the $2^m\times2^m$ matrix on the right is the $m$-th Kronecker
power of $M(\mu)$.

When $\mu=1$, we have the identity transform, while when $\mu=-1$,
we have a scalar multiple of the Hadamard transform.  It is well known that
the Hadamard transform takes the indicator function of a linear space
to a scalar multiple of the indicator function of its dual, from which
it follows that the indicator functions of the cutset and circuit spaces
of a graph are related by the Hadamard transform in the same way.
It was shown in \cite{farr93} that general matroid duality is also
described by the Hadamard transform.

It is easy to show that $M(\mu_1\mu_2)=M(\mu_1)M(\mu_2)$, see \cite{farr2013}.
It follows (using the mixed-product property for the Kronecker product)
that composition of the $L^{[\mu]}$ transforms corresponds to multiplication
of their $\mu$ paramters: $L^{[\mu_1]} L^{[\mu_2]} = L^{[\mu_1\mu_2]}$,
from \cite[Theorem 2]{farr2013}.  At this point, readers may ask: what
happens when $\mu=\omega$?  We look at this shortly.

Suppose $E=\{e_0,\ldots,e_{m-1}\}$.

We use $\mathbf{f}_{b\bullet}$ as shorthand for the
vector of length $2^{m-1}$, with elements indexed by
subsets of $E\setminus e_0$, whose $X$-element
is $f(X)$, if $b=0$, or $f(X\cup\{e_0\})$, if $b=1$
(for $X\subseteq E\setminus\{e_0\}$).
We define $\mathbf{f}_{\bullet b}$ in the same way,
except that we use $e_{m-1}$ instead of $e_0$ throughout.
The vectors $\mathbf{f}_{0\bullet}$ and $\mathbf{f}_{1\bullet}$
give the top and bottom halves, respectively, of $\mathbf{f}$,
while $\mathbf{f}_{\bullet 0}$ and $\mathbf{f}_{\bullet 1}$
give the elements in even and odd positions, respectively,
of $\mathbf{f}$.

Let $I_l$ denote
the $l\times l$ identity matrix.
If $e\in E$, then the $[\mu]$\textit{-minor} of $\mathbf{f}$ by $e$ is
the $2^{m-1}$-element
vector $\mathbf{f}\minor{[\mu]}e$, with entries indexed by
subsets of $E\setminus\{e\}$, given by
\begin{equation}
\label{eq:bin-fn-minor}
\mathbf{f}\minor{[\mu]}e_i ~ := ~ c \cdot  ~(~ I_2^{\otimes i} \otimes
( ~~ 1 ~~~ \frac{1+\mu}{\sqrt{2}+1-(\sqrt{2}-1)\mu} ~~ )
\otimes I_2^{\otimes(m-i-1)} ~)~ \mathbf{f} ,
\end{equation}
where $c$ is such that the $\emptyset$-element
of $\mathbf{f}\minor{[\mu]}e_i$ is 1.

Put $(\mu_0,\mu_1,\mu_2)=(1,\omega,\omega^2)$ and, for each $j\in\{0,1,2\}$,
\begin{equation}
\label{eq:lambda-mu}
\lambda_j := \frac{1+\mu_j}{\sqrt{2}+1-(\sqrt{2}-1)\mu_j} .
\end{equation}
Then $\mathbf{f}\minor{[\mu_j]}e_i$ is a scalar multiple of
$(~ I_2^{\otimes i} \otimes
( ~~ 1 ~~~ \lambda_j ~~ )
\otimes I_2^{\otimes(m-i-1)} ~)~ \mathbf{f}$.

When $f$ is the indicator function of the cutset space of a graph,
the minor $f\minor{[\mu]}e$ amounts to deletion when $\mu=1$ and
contraction when $\mu=-1$.  See \cite[\S2,\S6]{farr2013}, and also
\cite{farr04} for the first definition of generalised minor operations
interpolating between deletion and contraction (albeit with a different
parameterisation to that used here and in \cite{farr2013}).
This work has its roots in \cite{farr93}, where deletion and contraction
are expressed in terms of indicator
functions of cutset spaces, and these operations are extended to
general binary functions.

It is shown in \cite[Theorem 9]{farr2013} that,
for all $\mu_1\in\mathbb{C}\setminus\{0\}$ and $\mu_2\in\mathbb{C}$,
\[
(L^{[\mu_1]}f)\minor{[\mu_2/\mu_1]}e  \simeq
L^{[\mu_1]}(f\minor{[\mu_2]}e) .
\]
In particular, we have
\begin{eqnarray*}
(L^{[\omega]}f)\minor{[1]}e  &  \simeq  &
L^{[\omega]}(f\minor{[\omega]}e) ,   \\
(L^{[\omega]}f)\minor{[\omega]}e  &  \simeq  &
L^{[\omega]}(f\minor{[\omega^2]}e) ,   \\
(L^{[\omega]}f)\minor{[\omega^2]}e  &  \simeq  &
L^{[\omega]}(f\minor{[1]}e) .
\end{eqnarray*}
This relationship between the transform $L^{[\omega]}$
(called the \textit{trinity transform} \cite{farr2013} or
\textit{triality transform}) and
the minor operations 
for binary functions follows the same pattern as the relationships
between triality and minors for alternating dimaps,
given in Theorem \ref{thm:trial-minor}.
It is natural to ask what connection there may be between the two.

For binary functions, the minor operations always
commute \cite[Lemma 4]{farr04}.  In fact, that result implies that
every binary function is totally reduction-commutative (using the
natural analogue of that definition for binary functions).
But, as we saw in \S\ref{sec:non-commut},
the minor operations for alternating dimaps do not always commute.
It follows that alternating dimaps, along with triality and minor
operations, cannot be represented faithfully by binary functions with
their triality transform and minor operations described above. 

Nonetheless, we can ask if there is a subclass of alternating dimaps
which can be represented faithfully by binary functions in this way.
For this to occur, this subclass must consist only of alternating
dimaps that are totally reduction-commutative.  Such alternating
dimaps are disjoint unions of tricircuits, by
Theorem \ref{thm:totally-reduction-commut-tricircuits}.

Later we will give a definition of faithful representation by
binary functions, and determine when such a representation is possible.
To do the latter, it will help to characterise those binary functions
for which any reduction, on any element of the ground set, gives the
same given binary function.

To do this, we need some more notation.

Throughout, we write
\[
\mathbf{i}=\vr{1}{0}, ~~~~
\mathbf{j}=\vr{0}{1}, ~~~~
H = (\mathbf{h}_0,\ldots,\mathbf{h}_{k-1})
\in \{\mathbf{i},\mathbf{j}\}^{\{0,\ldots,k-1\}} .
\]
For each $i$,
\[
H^{(i)} = (\mathbf{h}_0,\ldots,\mathbf{h}_{i-1},
\mathbf{h}_{i+1},\ldots,\mathbf{h}_{k-1})
\]
is the sequence obtained from $H$ by omitting the term indexed by $i$.

For each $H$, define the sequence $G=G(H)=(g_0,\ldots,g_{k-1})$ by
\[
g_i = \left\{
\begin{array}{ll}
0,  &  \hbox{if $\mathbf{h}_i = \mathbf{i}$};   \\
1,  &  \hbox{if $\mathbf{h}_i = \mathbf{j}$}.
\end{array}
\right.
\]
The sequence obtained from this by omitting the term indexed by $i$ is
\[
G^{(i)} = (g_0,\ldots,g_{i-1},g_{i+1},\ldots,g_{k-1}) .
\]
The subsequence $(g_{i_1},\ldots,g_{i_2})$ of $G$ is denoted by
$G[i_1\dd i_2]$.

If $b\in\{0,1\}$, then $G:i\leftarrow b$ denotes the sequence obtained by inserting $b$
between the $i$-th and $(i+1)$-th elements of $G$:
\[
G:i\leftarrow b ~=~ (g_0, \ldots, g_{i-1}, b, g_i, \ldots, g_{k-1}) .
\]
The two-element vector $\mathbf{f}_{G:i}$ is defined by
\[
\mathbf{f}_{G:i} = \vr{f_{G:i\leftarrow0}}{f_{G:i\leftarrow1}} .
\]

Write $\mathbf{u}$ for a $2^k$-element vector indexed by the numbers
$0,\ldots,2^{k-1}$ --- or, equivalently, by vectors of $k$ bits, or by subsets of
$\{0,\ldots,k-1\}$.

For a given $G$, we write $u_{G}$ for the entry of $\mathbf{u}$ whose
index has binary representation given by $G$, i.e., whose index is
$\sum_{i=0}^{k-1}g_i2^{k-1-i}$.
%  Write $\mathbf{u}_{G}^{(i)}$ for the vector
%  \[
%  \vr{u_{G[0\dd i-1],0,G[i+1\dd k-1]}}{u_{G[0\dd i-1],1,G[i+1\dd k-1]}} ,
%  \]
%  where each subscript is a sequence of $k$ bits, interpreted for indexing
%  purposes as a binary number.

It is routine to show that, if $m\ge1$ and
$u$ is a (vector representation of a)
binary function with ground set of size $m-1$, then
\begin{equation}
\label{eq:u-sum}
\mathbf{u} = \sum_{H}
(~ \mathbf{h}_0 \otimes \cdots \otimes
\mathbf{h}_{m-2} ~) ~ u_G .
\end{equation}
If $m=1$ then there is a single $H$ to sum over, consisting of the
empty sequence, and the empty product
$\mathbf{h}_0 \otimes \cdots \otimes
\mathbf{h}_{m-2}$ is the trivial single-element vector $~(~ 1 ~)$.
Also $G$ is the empty bit-sequence, representing the number 0, and
$u_G=1$, so $\mathbf{u} = (~1~)$, as expected.

%  Also, though I'm not sure yet if I'll need it, we have ********
%  \[
%  \mathbf{u} = \sum_{H^{(i)}}
%  \mathbf{h}_0 \otimes \cdots \otimes
%  \mathbf{h}_{i-1} \otimes \mathbf{u}_{G}^{(i)} \otimes
%  \mathbf{h}_{i+1} \otimes \cdots \otimes
%  \mathbf{h}_{m-2} .
%  \]

\begin{lemma}
\label{lemma:reduce-f-to-u-soln}
Suppose $\mathbf{f}$ and $\mathbf{u}$ are binary functions with
\[
\mathbf{f}\minor{[\mu]}e_i = \mathbf{u}, ~~~~
\hbox{for all $\mu\in\{1,\omega,\omega^2\}$} .
\]
Then
for all $G\in\{0,1\}^{\{0,\ldots,m-2\}}$ and all $b\in\{0,1\}$,
\begin{equation}
\label{eq:reduce-f-to-u-soln}
f_{G:i\leftarrow b}  =
f_{0:i\leftarrow b} \, u_{G} \,.
\end{equation}
\end{lemma}

\pf
Let us write the hypothesis as a set of equations, using
(\ref{eq:bin-fn-minor}) and (\ref{eq:lambda-mu}).
If $\mathbf{f}\minor{[\mu_j]}e_i=\mathbf{u}$ for all $j$,
then, for each $j$, there exists $c_{ij}$ such that
\begin{equation}
\label{eq:reduce-f-to-u}
(~ I_2^{\otimes i} \otimes
( ~~ 1 ~~~ \lambda_j ~~ )
\otimes I_2^{\otimes(m-i-1)} ~)~ \mathbf{f} = c_{ij} \, \mathbf{u} \,.
%  ~~~~~~~~
%  i=0,\ldots,m-1, ~~~ j=1,\omega,\omega^2,
\end{equation}

Put
\[
R = \left(
\begin{array}{cc}
1 & \lambda_0   \\
1 & \lambda_1   \\
1 & \lambda_2
\end{array}
\right)
~~~~~~ \hbox{and} ~~~~~~
\mathbf{c}_i =
\left(
\begin{array}{c}
c_{i0}  \\  c_{i1}  \\  c_{i2}
\end{array}
\right) .
\]
The equations (\ref{eq:reduce-f-to-u}) may be written (using (\ref{eq:u-sum})),
\begin{equation}
\label{eq:reduce-f-to-u-new}
(~ I_2^{\otimes i} \otimes  R
\otimes I_2^{\otimes(m-i-1)} ~)~ \mathbf{f}
= \sum_{H}
(~ \mathbf{h}_0 \otimes \cdots \otimes \mathbf{h}_{i-1} \otimes
\mathbf{c}_i \otimes
\mathbf{h}_{i} \otimes \cdots \otimes \mathbf{h}_{m-2} ~) ~ u_G .
\end{equation}
We show by induction on $m$ that the solutions to this equation satisfy
\begin{equation}
\label{eq:reduce-f-to-u-soln1}
R \, \mathbf{f}_{G:i}  =  \mathbf{c}_i \, u_G  \,,
\end{equation}
for all $G\in\{0,1\}^{\{0,\ldots,m-2\}}$.

For the inductive basis, let $m=1$, so $i=0$.
Then we may write
\[
\mathbf{f} = \vr{1}{f_1} , ~~~~~~ \mathbf{u} = ~(~ 1 ~) 
\]
(where $u_0=1$ since $\mathbf{u}$ is a binary function),
and our equations (\ref{eq:reduce-f-to-u-new}) are
\[
R \, \mathbf{f} = 
\mathbf{c}_0 \,.
\]
In this case, any $\mathbf{f}$ will do,
with appropriate choice of $\mathbf{c}_0$.
Equations (\ref{eq:reduce-f-to-u-soln1})
are satisfied in this case (with there being just a single
sequence $G$, which is the empty sequence and represents the number 0).
Observe that $\mathbf{f}_{G:i}=\mathbf{f}_{G:0}=\mathbf{f}$.

Now suppose that the claim is true for $m=k\ge2$.  We show
that it is true for $m=k+1$.

We wish to solve (\ref{eq:reduce-f-to-u-new}) when $m=k+1$.
Since $k\ge1$, either $i>0$ or $k-i>0$ (or both).  We treat these
two cases in turn.

Suppose $i>0$.

We can write
\[
\mathbf{f} = \mathbf{i} \otimes \mathbf{f}_{0\bullet} +
\mathbf{j} \otimes \mathbf{f}_{1\bullet} .
\]

This allows us to rewrite the left-hand side of (\ref{eq:reduce-f-to-u-new}):
\begin{eqnarray*}
\lefteqn{
(~ I_2^{\otimes i} \otimes  R
\otimes I_2^{\otimes(m-i-1)} ~)~ \mathbf{f} }   \\
& = &
(~ I_2 \otimes I_2^{\otimes(i-1)} \otimes  R
\otimes I_2^{\otimes(m-i-1)} ~)~
( \mathbf{i} \otimes \mathbf{f}_{0\bullet} +
\mathbf{j} \otimes \mathbf{f}_{1\bullet} )   \\
& = &
I_2 \mathbf{i} \otimes
~(~ I_2^{\otimes(i-1)} \otimes  R
\otimes I_2^{\otimes(m-i-1)} ~)~ \mathbf{f}_{0\bullet} +
I_2 \mathbf{j} \otimes
~(~ I_2^{\otimes(i-1)} \otimes  R
\otimes I_2^{\otimes(m-i-1)} ~)~ \mathbf{f}_{1\bullet}   \\
& = &
\mathbf{i} \otimes
~(~ I_2^{\otimes(i-1)} \otimes  R
\otimes I_2^{\otimes((m-1)-(i-1)-1)} ~)~ \mathbf{f}_{0\bullet} +
\mathbf{j} \otimes
~(~ I_2^{\otimes(i-1)} \otimes  R
\otimes I_2^{\otimes((m-1)-(i-1)-1)} ~)~ \mathbf{f}_{1\bullet}   \\
\end{eqnarray*}

On the other hand, the right-hand side of (\ref{eq:reduce-f-to-u-new})
may be rewritten: 
\begin{eqnarray*}
\lefteqn{\sum_{H}
(~ \mathbf{h}_0 \otimes \cdots \otimes \mathbf{h}_{i-1} \otimes
\mathbf{c}_i \otimes
\mathbf{h}_{i} \otimes \cdots \otimes \mathbf{h}_{m-2} ~) ~ u_G }   \\
& = &
\sum_{\mathbf{h}_0, H^{(0)}}
(~ \mathbf{h}_0 \otimes \mathbf{h}_1 \otimes
\cdots \otimes \mathbf{h}_{i-1} \otimes
\mathbf{c}_i \otimes
\mathbf{h}_{i} \otimes \cdots \otimes \mathbf{h}_{m-2} ~) ~ u_G   \\
& = &
\sum_{H^{(0)}}
(~ \mathbf{i} \otimes \mathbf{h}_1 \otimes
\cdots \otimes \mathbf{h}_{i-1} \otimes
\mathbf{c}_i \otimes
\mathbf{h}_{i} \otimes \cdots \otimes \mathbf{h}_{m-2} ~) ~ u_{0G^{(0)}}  +   \\
&&
\sum_{H^{(0)}}
(~ \mathbf{j} \otimes \mathbf{h}_1 \otimes
\cdots \otimes \mathbf{h}_{i-1} \otimes
\mathbf{c}_i \otimes
\mathbf{h}_{i} \otimes \cdots \otimes \mathbf{h}_{m-2} ~) ~ u_{1G^{(0)}}    \\
& = &
\mathbf{i} \otimes
\sum_{H^{(0)}}
(~ \mathbf{h}_1 \otimes
\cdots \otimes \mathbf{h}_{i-1} \otimes
\mathbf{c}_i \otimes
\mathbf{h}_{i} \otimes \cdots \otimes \mathbf{h}_{m-2} ~) ~ u_{0G^{(0)}}  +   \\
&&
\mathbf{j} \otimes
\sum_{H^{(0)}}
(~ \mathbf{h}_1 \otimes
\cdots \otimes \mathbf{h}_{i-1} \otimes
\mathbf{c}_i \otimes
\mathbf{h}_{i} \otimes \cdots \otimes \mathbf{h}_{m-2} ~) ~ u_{1G^{(0)}} .    \\
&&
\end{eqnarray*}
These rewritten forms of each side give a top half and a bottom half for each.
Equating these tells us that (\ref{eq:reduce-f-to-u-new}) is equivalent
to the following two simultaneous equations.
\begin{eqnarray*}
~(~ I_2^{\otimes(i-1)} \otimes  R
\otimes I_2^{\otimes((m-1)-(i-1)-1)} ~)~ \mathbf{f}_{0\bullet}
& = &
\sum_{H^{(0)}}
(~ \mathbf{h}_1 \otimes
\cdots \otimes \mathbf{h}_{i-1} \otimes
\mathbf{c}_i \otimes
\mathbf{h}_{i} \otimes \cdots \otimes \mathbf{h}_{m-2} ~) ~ u_{0G^{(0)}}  \,, \\
~(~ I_2^{\otimes(i-1)} \otimes  R
\otimes I_2^{\otimes((m-1)-(i-1)-1)} ~)~ \mathbf{f}_{1\bullet}
& = &
\sum_{H^{(0)}}
(~ \mathbf{h}_1 \otimes
\cdots \otimes \mathbf{h}_{i-1} \otimes
\mathbf{c}_i \otimes
\mathbf{h}_{i} \otimes \cdots \otimes \mathbf{h}_{m-2} ~) ~ u_{1G^{(0)}}  \,.
\end{eqnarray*}
Each of these equations is an instance of the same type of equation
as (\ref{eq:reduce-f-to-u-new}), with dimension $m$ and position $i$
each reduced by one.  So, by
the inductive hypothesis, their solutions are
\begin{eqnarray*}
R \, \mathbf{f}_{0G^{(0)}:i}  & = &  \mathbf{c}_i \, u_{0G^{(0)}}  \,,   \\
R \, \mathbf{f}_{1G^{(0)}:i}  & = &  \mathbf{c}_i \, u_{1G^{(0)}}  \,.
\end{eqnarray*}
Combining these gives (\ref{eq:reduce-f-to-u-soln1}), for all $G$.
This deals with the case $i>0$.

If $k-i>0$ then a similar argument can be used, peeling off
the identity matrix from the right, rather than the left, of
$I_2^{\otimes i} \otimes R \otimes I_2^{\otimes(m-i-1)}$
in (\ref{eq:reduce-f-to-u-new}), and using
$\mathbf{f} = \mathbf{f}_{\bullet0} \otimes \mathbf{i} +
\mathbf{f}_{\bullet1} \otimes \mathbf{j}$, and so on.

It follows by induction that the solutions to (\ref{eq:reduce-f-to-u-new})
satisfy (\ref{eq:reduce-f-to-u-soln1}).

When $G=0$, (\ref{eq:reduce-f-to-u-soln1}) and $u_0=1$ give
\[
R \, \mathbf{f}_{0:i}  =  \mathbf{c}_i .
\]
Using this with (\ref{eq:reduce-f-to-u-soln1}) gives, for any $G$,
\[
R \, \mathbf{f}_{G:i} = R \, \mathbf{f}_{0:i} \, u_G .
\]
Since $R$ has rank 2 (because the $\lambda_i$ are distinct),
\[
\mathbf{f}_{G:i} = \mathbf{f}_{0:i} \, u_G .
\]
Hence, for all $b\in\{0,1\}$ and all $G$, 
\[
f_{G:i\leftarrow b}  =
f_{0:i\leftarrow b} \, u_{G} \,.
\]
\eopf   \\

We now define our notion of faithful representation, and then
determine when it is possible.   \\

\noindent\textbf{Definition}

A \textit{strict binary representation} of a minor-closed set
$\mathcal{A}$ of alternating dimaps is a triple $(F,\varepsilon,\nu)$
such that
\begin{itemize}
\item[(a)]
$F : \mathcal{A} \rightarrow \{\hbox{binary functions}\}$
\item[(b)]
$\varepsilon=(\varepsilon_G\mid G\in\mathcal{A})$ is a family
of bijections $\varepsilon_G:E(G)\rightarrow E(F(G))$;
\item[(c)]
$\nu\in\mathbb{C}$ with $|\nu|=1$;
%  \item[(d)]    %  Not needed: covered by (b).
%   $\varepsilon_G(e)\in E(F(G))$
%   for all $G\in\mathcal{A}$ and all $e\in E(G)$;
\item[(d)]
$F(G^{(\omega)}) \simeq L^{[\omega]} F(G)$ for all $G\in\mathcal{A}$;
\item[(e)]
$F(G\mino{\mu}e) \simeq F(G)\minor{[\nu\mu]}\varepsilon_G(e)$ for all
$G\in\mathcal{A}$, $e\in E(G)$ and $\mu\in\{1,\omega,\omega^2\}$.
\end{itemize}

\vspace{0.1in}

Let $C_1$ denote the ultraloop.
We write $\mathcal{U}_k=\{iC_1\mid i=0,\ldots,k\}$ and
$\mathcal{U}_{\infty}=\{iC_1\mid i\in\mathbb{N}\cup\{0\}\}$,
where $0C_1$ is the empty alternating dimap.

\begin{thm}
If $\mathcal{A}$ is a minor-closed class of alternating dimaps
which has a strict binary representation then $\mathcal{A}=\emptyset$,
or $\mathcal{A}=\mathcal{U}_k$
for some $k$, or $\mathcal{A}=\mathcal{U}_{\infty}$.
\end{thm}

\pf
Suppose $(F,\varepsilon,\nu)$ is
a strict binary representation of $\mathcal{A}$.

The theorem is immediately true if $\mathcal{A}=\emptyset$.
So suppose $\mathcal{A}\not=\emptyset$.

If $|\mathcal{A}|\ge1$ then, since it is minor-closed, it must
contain the empty alternating dimap $C_0$, and the image $F(C_0)$,
representing $C_0$ as a binary function, must be the binary function
$f:2^{\emptyset}\rightarrow\mathcal{C}$ defined by $f(\emptyset)=1$.

So, if $|\mathcal{A}|=1$ then $\mathcal{A}=\mathcal{U}_0$, and the previous paragraph gives
a strict binary representation of $\mathcal{A}$.

Similarly, if $|\mathcal{A}|\ge2$, then it must contain the
ultraloop $C_1$, since that is the only alternating dimap on one edge.   \\

Claim 1:  The image $F(C_1)$ of the ultraloop $C_1$ is given by
\[
F(C_1)  =  \vr{1}{\sqrt{2}-1} .
\]

Proof:

$F(C_1)$ must be some binary
function $f$ on a singleton ground set, $E=\{e\}$ say, with
$f(\emptyset)=1$ and $f(\{e\})=u$ for some $u\in\mathbb{C}$.
Since $C_1$ is self trial, so must $f$ be (by (d) above).  This means
that its vector form $\mathbf{f}=\textvr{1}{u}$ must
%  satisfy
%  \[
%  M(\omega) \mathbf{f} = \mathbf{f} ,
%  \]
%  i.e., it must
be an eigenvector for eigenvalue 1
of the matrix $M(\omega)$.  Now this matrix has eigenvalues 1 and $\omega$,
and the eigenvectors for the former are the scalar multiples of
\[
\vr{1}{\sqrt{2}-1} .
\]
So this is $F(C_1)$, and $u=\sqrt{2}-1$.  So Claim 1 is proved.   \\

If $|\mathcal{A}|=2$ then $\mathcal{A}$ consists of just the empty alternating
dimap and the ultraloop.  The $F$ given by Claim 1, together with
appropriate identity maps $\varepsilon$ (and, in fact, any $\nu$),
gives a strict binary representation.  So we are done in this case.   \\

It remains to deal with $|\mathcal{A}|\ge3$, when $\mathcal{A}$ contains
at least one alternating dimap on two edges.   \\

Claim 2:
The only binary function $\mathbf{f}$ with the property that 
every reduction, on any of the elements of its ground set,
gives $\mathbf{u}=F(C_1)^{\otimes k}$, is 
$\mathbf{f}=F(C_1)^{\otimes(k+1)}$.   \\

Proof:

Observe that, by Claim 1, $\mathbf{u}=(u_G\mid G\in\{0,1\}^E)$ where
$u_G=(\sqrt{2}-1)^{|G|}$, where $|E|=k$ and $|G|$ is the number of 1s in $G$.

Applying Lemma \ref{lemma:reduce-f-to-u-soln},
for all $i\in\{0,\ldots,k\}$,
to $\mathbf{u} = F(C_1)^{\otimes k}$ gives
\[
f_{G:i\leftarrow b}  =
f_{0:i\leftarrow b} \, u_{G} \,.
\]
Hence, for each $i$ and each $G$,
\[
f_{G:i\leftarrow 0}  =
f_{0_{k+1}} \, u_{G}  =
u_{G}  =  (\sqrt{2}-1)^{|G|} .
\]
Now consider $f_{G:i\leftarrow 1}$.  Put $j:=0$ if $i\not=0$ and $j:=1$
otherwise (so $j\not=i$).  Then
\begin{eqnarray*}
f_{G:i\leftarrow 1}  & = &
f_{0_{k}:i\leftarrow 1} \, u_{G}  =
f_{(0_{k-1}:i\leftarrow 1):j\leftarrow 0} \, u_{G}  =
f_{0_{k}:j\leftarrow 0} \, u_{0_{k-1}:i\leftarrow 1} \, u_{G}  =
f_{0_{k+1}} \, u_{0_{k-1}:i\leftarrow 1} \, u_{G}  \\
&  = &
1 \cdot  (\sqrt{2}-1) \cdot (\sqrt{2}-1)^{|G|}  =
(\sqrt{2}-1)^{|G|+1} .
\end{eqnarray*}
It follows that, for all $G'\in\{0,1\}^{k+1}$,
\[
f_{G'} = (\sqrt{2}-1)^{|G'|} .
\]
Therefore 
\[
\mathbf{f} = F(C_1)^{\otimes(k+1)} ,
\]
proving the Claim.   \\

Claim 3:
If $k\ge2$ and every reduction of $G$ is $kC_1$, then $G=(k+1)C_1$.   \\

Before proving the claim, consider the case $k=1$, which it does not cover.
Then every alternating dimap on two edges (of which there are four) has
the claimed property.  Of these, the only self-trial one is $2C_1$.   \\

Proof:

Suppose $k=2$.  If $G$ is connected, then there must be some
$e\in E(G)$ and some $\mu\in\{1,\omega,\omega^2\}$ such that
$G\mino{\mu}e=2C_1$ and is therefore disconnected.  The only way in which
$\mu$-reducing a single edge can disconnect a connected alternating dimap is if the edge
is a proper $\mu^{-1}$-semiloop.  It is easily determined that the
only alternating dimaps on three edges which have this property are those
consisting of two triloops and a semiloop.  These do not have three proper semiloops.
So, although they have the specified property for
\textit{one} of their edges, they do not have it for \textit{all} of their
edges.  So $G$ must be disconnected.  Since $G$ has only three edges, some
component of $G$ must be an ultraloop.  But this disappears when reduced,
so the rest of $G$ must be $2C_1$, so $G=3C_1$.

Now suppose $k\ge3$.  It is impossible for $G$ to be connected, because
no reduction of any edge of any connected alternating dimap can possibly
break it up into three or more components.  So consider the components of $G$.
If any of these is not an ultraloop,
then it has at least two edges, and also is left unchanged by reduction of
any edge in any other component (of which there must be at least one), so we
would have a reduction of $G$ that does not give $kC_1$,
which is a contradiction.  So every component of $G$ must be an ultraloop.
Each of these just disappears on reduction, giving $kC_1$, as desired.

So Claim 3 is proved.   \\

Claim 4:
For all $k\ge0$, either $\mathcal{A}$ has no
members with $k$ edges, or it has just one such member which is $kC_1$.   \\

Proof:

We prove the claim by induction on $k$.

We have seen that this is true already for $k\le1$.

Suppose $k=2$.
Every alternating dimap $G_2$
on two edges has the property that every
reduction of it gives the ultraloop.  Therefore, if $G_2\in\mathcal{A}$
then $F(G_2)=F(C_1)^{\otimes2}$.  But $F(C_1)^{\otimes2}$ is self-trial,
since $F(C_1)$ is.  Therefore $G_2$ must be self-trial too.  But the only
self-trial alternating dimap on two edges is $2C_1$.
So the only member of $\mathcal{A}$ with two edges
is $2C_1$.

Now suppose it is true regarding members of $\mathcal{A}$ with $k-1$
edges, where $k\ge3$.  We show that it is true for $k$ edges.

If $\mathcal{A}$ has no members with $k-1$ edges, then it can have
no members with $k$ edges either, since it is minor-closed.

If $\mathcal{A}$ has at least one member with $k-1$ edges, then by the
inductive hypothesis it can have only one such member, and
this must be $(k-1)C_1$.  We must show that, if $\mathcal{A}$ has at
least one member with $k$ edges, then it can have only one, and it is $kC_1$.

Let $G$ be a member of $\mathcal{A}$ with $k$ edges.  Since $\mathcal{A}$
is minor-closed and has $(k-1)C_1$ as its only member with $k-1$ edges,
all reductions of $G$ must give $(k-1)C_1$.  So, by the requirements of
a strict binary representation, all reductions of $F(G)$ must give
$F(C_1)^{\otimes(k-1)}$.  This implies that $F(G)=F(C_1)^{\otimes k}$,
by Claim 3.  This completes the proof of Claim 4.   \\

It follows from Claim 4 that $\mathcal{A}$ can only be one of the classes
given in the statement of the theorem.  It remains to establish that
a strict binary representation is possible for each of those classes.
This is routine, using
\[
F(kC_1)  =  \vr{1}{\sqrt{2}-1}^{\otimes k}
\]
for every $k$ for which $kC_1\in\mathcal{A}$.
Let $\varepsilon$ consist just of identity maps.
To show that this does indeed enable a strict binary representation,
use Claims 1--3.  The details are a routine exercise.
\eopf   \\

It is possible to develop broader definitions of binary
representations of classes of alternating dimaps.
For example, we could allow the edges of $G$ to be
represented by disjoint subsets of elements of the ground set of $F(G)$
instead of just by distinct single elements.   \\

\noindent\textbf{Problem}

Characterise those minor-closed classes of alternating dimaps that
have binary function representations of a more general type, such
as that suggested above.

\section{Excluded minors for fixed genus}
\label{sec:excl-minors-fixed-genus}

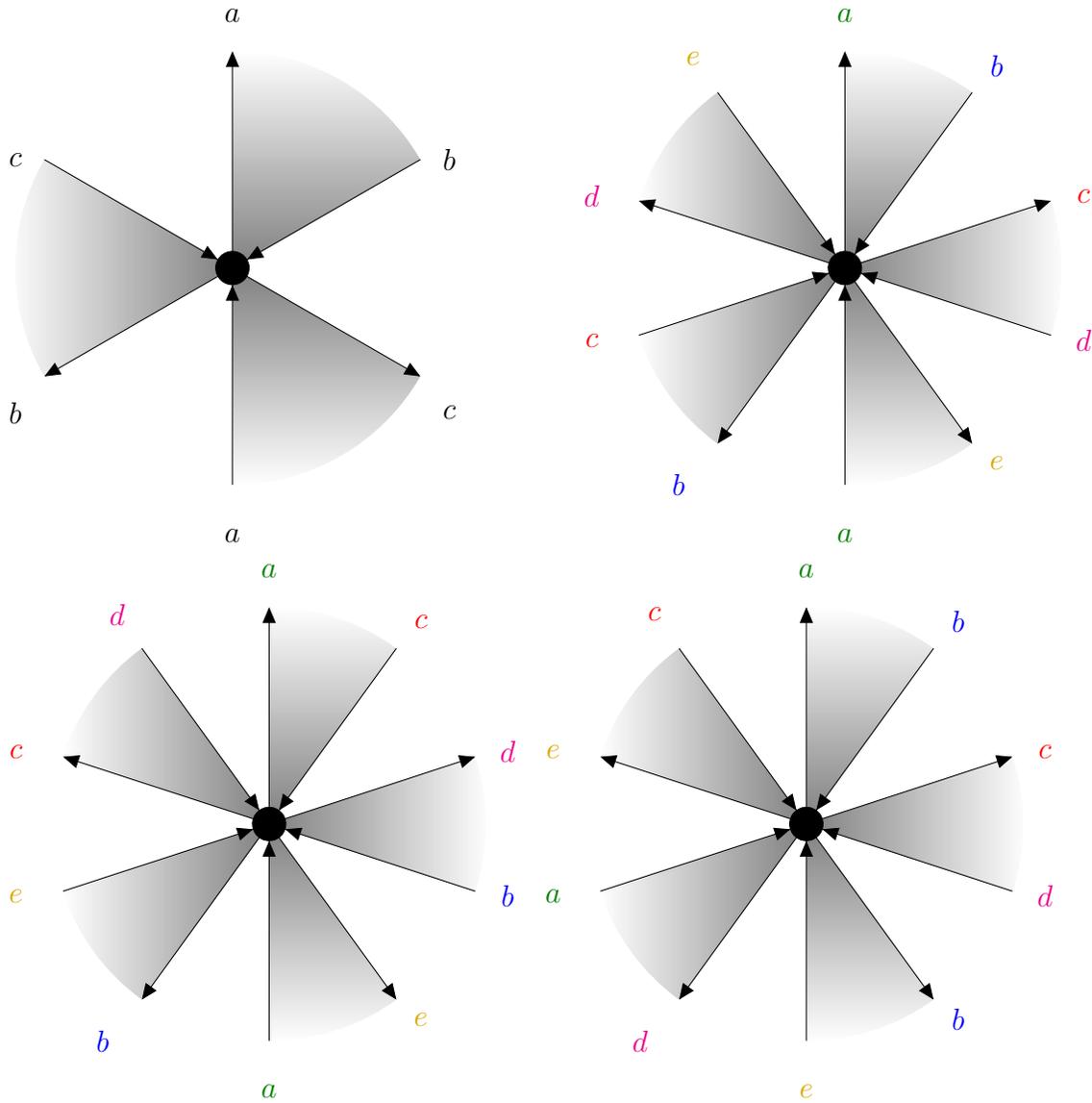
\begin{figure}[ht]
\newcommand{\alabel}{\textcolor{greenish}{$a$}}
\newcommand{\blabel}{\textcolor{blue}{$b$}}
\newcommand{\clabel}{\textcolor{red}{$c$}}
\newcommand{\dlabel}{\textcolor{magenta}{$d$}}
\newcommand{\elabel}{\textcolor{gold}{$e$}}
\begin{tikzpicture}
\begin{scope}[>=triangle 45]
\coordinate (v) at (4,4);
\coordinate (a1) at (4,7);
\coordinate (b1) at (6.6,5.5);   %  6.6 := 4 + 1.5*sqrt(3)
\coordinate (c1) at (6.6,2.5);
\coordinate (a2) at (4,1);
\coordinate (b2) at (1.4,2.5);   %  1.4 := 4 - 1.5*sqrt(3)
\coordinate (c2) at (1.4,5.5);
\coordinate (a1label) at (4,7.5);
\coordinate (b1label) at (7,5.5);
\coordinate (c1label) at (7,2);
\coordinate (a2label) at (4,0.3);
\coordinate (b2label) at (1,2);
\coordinate (c2label) at (1,5.5);
\shade[bottom color=gray, top color=white] (v) -- (a1) arc (90:30:3) -- cycle;
\shade[top color=gray, bottom color=white] (v) -- (a2) arc (270:330:3) -- cycle;
\shade[right color=gray, left color=white] (v) -- (b2) arc (210:150:3) -- cycle;
\draw[fill] (v) circle (0.23cm);
\draw[->] (v) -- (a1);
\draw[->>] (b1) -- (v);
\draw[->] (v) -- (c1);
\draw[->>] (a2) -- (v);
\draw[->] (v) -- (b2);
\draw[->>] (c2) -- (v);
\draw (a1label) node {$a$};
\draw (b1label) node {$b$};
\draw (c1label) node {$c$};
\draw (a2label) node {$a$};
\draw (b2label) node {$b$};
\draw (c2label) node {$c$};
\end{scope}
\end{tikzpicture}
\begin{tikzpicture}
\begin{scope}[>=triangle 45]
\coordinate (v) at (4,4);
\coordinate (a1) at (4,7);
\coordinate (b1) at (5.76,6.43);
%   5.76 := 4 + 3*cos(54 deg), 6.43 := 4 + 3*sin(54 deg)
\coordinate (c1) at (6.85,4.93);
%   6.85 := 4 + 3*cos(18 deg), 4.93 := 4 + 3*sin(18 deg)
\coordinate (d1) at (6.85,3.07);
\coordinate (e1) at (5.76,1.57);
\coordinate (a2) at (4,1);
\coordinate (b2) at (2.24,1.57);
\coordinate (c2) at (1.15,3.07);
\coordinate (d2) at (1.15,4.93);
\coordinate (e2) at (2.24,6.43);
\coordinate (a1label) at (4,7.5);
\coordinate (b1label) at (6.1,6.8);
\coordinate (c1label) at (7.3,5);
\coordinate (d1label) at (7.3,3);
\coordinate (e1label) at (6.1,1.3);
\coordinate (a2label) at (4,0.3);
\coordinate (b2label) at (1.7,1);
\coordinate (c2label) at (0.5,3);
\coordinate (d2label) at (0.5,5);
\coordinate (e2label) at (1.9,6.9);
\shade[bottom color=gray, top color=white] (v) -- (a1) arc (90:54:3) -- cycle;
\shade[left color=gray, right color=white] (v) -- (d1) arc (-18:18:3) -- cycle;
\shade[top color=gray, bottom color=white] (v) -- (a2) arc (270:306:3) -- cycle;
\shade[right color=gray, left color=white] (v) -- (b2) arc (234:198:3) -- cycle;
\shade[right color=gray, left color=white] (v) -- (e2) arc (126:162:3) -- cycle;
\draw[fill] (v) circle (0.23cm);
\draw[->] (v) -- (a1);
\draw[->>] (b1) -- (v);
\draw[->] (v) -- (c1);
\draw[->>] (d1) -- (v);
\draw[->] (v) -- (e1);
\draw[->>] (a2) -- (v);
\draw[->] (v) -- (b2);
\draw[->>] (c2) -- (v);
\draw[->] (v) -- (d2);
\draw[->>] (e2) -- (v);
\draw (a1label) node {\alabel};
\draw (b1label) node {\blabel};
\draw (c1label) node {\clabel};
\draw (d1label) node {\dlabel};
\draw (e1label) node {\elabel};
\draw (a2label) node {\alabel};
\draw (b2label) node {\blabel};
\draw (c2label) node {\clabel};
\draw (d2label) node {\dlabel};
\draw (e2label) node {\elabel};
\end{scope}
\end{tikzpicture}
\begin{tikzpicture}
\begin{scope}[>=triangle 45]
\coordinate (v) at (4,4);
\coordinate (a1) at (4,7);
\coordinate (b1) at (5.76,6.43);
%   5.76 := 4 + 3*cos(54 deg), 6.43 := 4 + 3*sin(54 deg)
\coordinate (c1) at (6.85,4.93);
%   6.85 := 4 + 3*cos(18 deg), 4.93 := 4 + 3*sin(18 deg)
\coordinate (d1) at (6.85,3.07);
\coordinate (e1) at (5.76,1.57);
\coordinate (a2) at (4,1);
\coordinate (b2) at (2.24,1.57);
\coordinate (c2) at (1.15,3.07);
\coordinate (d2) at (1.15,4.93);
\coordinate (e2) at (2.24,6.43);
\coordinate (a1label) at (4,7.5);
\coordinate (b1label) at (6.1,6.8);
\coordinate (c1label) at (7.3,5);
\coordinate (d1label) at (7.3,3);
\coordinate (e1label) at (6.1,1.3);
\coordinate (a2label) at (4,0.3);
\coordinate (b2label) at (1.7,1);
\coordinate (c2label) at (0.5,3);
\coordinate (d2label) at (0.5,5);
\coordinate (e2label) at (1.9,6.9);
\shade[bottom color=gray, top color=white] (v) -- (a1) arc (90:54:3) -- cycle;
\shade[left color=gray, right color=white] (v) -- (d1) arc (-18:18:3) -- cycle;
\shade[top color=gray, bottom color=white] (v) -- (a2) arc (270:306:3) -- cycle;
\shade[right color=gray, left color=white] (v) -- (b2) arc (234:198:3) -- cycle;
\shade[right color=gray, left color=white] (v) -- (e2) arc (126:162:3) -- cycle;
\draw[fill] (v) circle (0.23cm);
\draw[->] (v) -- (a1);
\draw[->>] (b1) -- (v);
\draw[->] (v) -- (c1);
\draw[->>] (d1) -- (v);
\draw[->] (v) -- (e1);
\draw[->>] (a2) -- (v);
\draw[->] (v) -- (b2);
\draw[->>] (c2) -- (v);
\draw[->] (v) -- (d2);
\draw[->>] (e2) -- (v);
\draw (a1label) node {\alabel};
\draw (b1label) node {\clabel};
\draw (c1label) node {\dlabel};
\draw (d1label) node {\blabel};
\draw (e1label) node {\elabel};
\draw (a2label) node {\alabel};
\draw (b2label) node {\blabel};
\draw (c2label) node {\elabel};
\draw (d2label) node {\clabel};
\draw (e2label) node {\dlabel};
\end{scope}
\end{tikzpicture}
\begin{tikzpicture}
\begin{scope}[>=triangle 45]
\coordinate (v) at (4,4);
\coordinate (a1) at (4,7);
\coordinate (b1) at (5.76,6.43);
%   5.76 := 4 + 3*cos(54 deg), 6.43 := 4 + 3*sin(54 deg)
\coordinate (c1) at (6.85,4.93);
%   6.85 := 4 + 3*cos(18 deg), 4.93 := 4 + 3*sin(18 deg)
\coordinate (d1) at (6.85,3.07);
\coordinate (e1) at (5.76,1.57);
\coordinate (a2) at (4,1);
\coordinate (b2) at (2.24,1.57);
\coordinate (c2) at (1.15,3.07);
\coordinate (d2) at (1.15,4.93);
\coordinate (e2) at (2.24,6.43);
\coordinate (a1label) at (4,7.5);
\coordinate (b1label) at (6.1,6.8);
\coordinate (c1label) at (7.3,5);
\coordinate (d1label) at (7.3,3);
\coordinate (e1label) at (6.1,1.3);
\coordinate (a2label) at (4,0.3);
\coordinate (b2label) at (1.7,1);
\coordinate (c2label) at (0.5,3);
\coordinate (d2label) at (0.5,5);
\coordinate (e2label) at (1.9,6.9);
\shade[bottom color=gray, top color=white] (v) -- (a1) arc (90:54:3) -- cycle;
\shade[left color=gray, right color=white] (v) -- (d1) arc (-18:18:3) -- cycle;
\shade[top color=gray, bottom color=white] (v) -- (a2) arc (270:306:3) -- cycle;
\shade[right color=gray, left color=white] (v) -- (b2) arc (234:198:3) -- cycle;
\shade[right color=gray, left color=white] (v) -- (e2) arc (126:162:3) -- cycle;
\draw[fill] (v) circle (0.23cm);
\draw[->] (v) -- (a1);
\draw[->>] (b1) -- (v);
\draw[->] (v) -- (c1);
\draw[->>] (d1) -- (v);
\draw[->] (v) -- (e1);
\draw[->>] (a2) -- (v);
\draw[->] (v) -- (b2);
\draw[->>] (c2) -- (v);
\draw[->] (v) -- (d2);
\draw[->>] (e2) -- (v);
\draw (a1label) node {\alabel};
\draw (b1label) node {\blabel};
\draw (c1label) node {\clabel};
\draw (d1label) node {\dlabel};
\draw (e1label) node {\blabel};
\draw (a2label) node {\elabel};
\draw (b2label) node {\dlabel};
\draw (c2label) node {\alabel};
\draw (d2label) node {\elabel};
\draw (e2label) node {\clabel};
\end{scope}
\end{tikzpicture}
\caption{The 1-posy and the three 2-posies.  In each posy,
the two faces are coloured grey (clockwise) and white (anticlockwise).}
\label{fig:posies}
\end{figure}

A \textit{posy}, or $k$-\textit{posy},
is an alternating dimap with one vertex, $2k+1$ edges
(all loops), and two faces.
Its genus is $k$.  Up to isomorphism, there is a single 0-posy,
a single 1-posy and
three 2-posies.  The 0-posy is just a single ultraloop.
The 1-posy and the three 2-posies are shown in
Figure \ref{fig:posies}.

\begin{thm}
\label{thm:genus-excl-posy-minors}
A nonempty alternating dimap $G$ has genus $< k$
if and only if none of its minors is a disjoint union of posies
of total genus $k$.
\end{thm}

\pf
The forward implication is clear, since every such union of
posies has genus $k$.

For the reverse implication, we prove by induction on $|E(G)|$ that,
if $G$ is nonempty and has genus $\ge k$, then it has, as a minor,
a disjoint union of posies of total genus $k$.

This is true for $|E(G)|=1$, since a nonempty planar alternating dimap must
have a directed cycle, hence a submap contractible to a loop, hence a 0-posy
minor.

Now suppose it is true for all alternating dimaps of $< m$ edges,
where $m>1$.  Let $G$ be any alternating dimap with $m$ edges and genus $\ge k$.
Let $e\in E(G)$.  Now, $G\mino{1}e$, $G\mino{\omega}e$ and
$G\mino{\omega^2}e$ each have $m-1$ edges, so by the inductive hypothesis,
$G\mino{\mu}e$ has as a minor a disjoint union of posies
of total genus $\gamma(G\mino{\mu}e)$,
for each $\mu\in\{1,\omega,\omega^2\}$.
Such a minor of $G\mino{\mu}e$ is also a minor of $G$, so we see that
$G$ has such a minor, for each such $\mu$.
If $\gamma(G\mino{\mu}e)=\gamma(G)$ for any such $\mu$, we are done.
So it remains to consider the case where $\gamma(G\mino{\mu}e)<\gamma(G)$
(in which case $\gamma(G\mino{\mu}e)=\gamma(G)-1$)
for each $\mu$ and each $e\in E(G)$.

\begin{figure}[ht]
%  copied from 1-posy fig, then modified
\begin{center}
\begin{tikzpicture}
\begin{scope}[>=triangle 45]
\coordinate (v) at (4,4);
\coordinate (a1) at (4,7);
\coordinate (b1) at (6.6,5.5);   %  6.6 := 4 + 1.5*sqrt(3)
\coordinate (c1) at (6.6,2.5);
\coordinate (a2) at (4,1);
\coordinate (b2) at (1.4,2.5);   %  1.4 := 4 - 1.5*sqrt(3)
\coordinate (c2) at (1.4,5.5);
\coordinate (e) at (3.5,1.5);
\coordinate (e1) at (2.2,2.5);
\coordinate (f) at (2.2,5.6);
\coordinate (F) at (5,2.2);
\coordinate (F1) at (2,4);
\coordinate (F2) at (5,5.5);
\draw[fill] (v) circle (0.23cm);
\draw[->] (v) -- (a1);
\draw[->>] (b1) -- (v);
\draw[->] (v) -- (c1);
\draw[->>] (a2) -- (v);
\draw[->] (v) -- (b2);
\draw[->>] (c2) -- (v);
\draw (e) node {$e$};
\draw (e1) node {$e'$};
\draw (f) node {$f$};
\draw (F) node {$F$};
\draw (F1) node {$F'$};
\draw (F2) node {$F''$};
\end{scope}
\end{tikzpicture}
\end{center}
\caption{Proof of Theorem \ref{thm:genus-excl-posy-minors}:
faces around $v$ when $e$ is a 1-semiloop and an $\omega$-semiloop.}
\label{fig:loops-and-faces-around-v}
\end{figure}
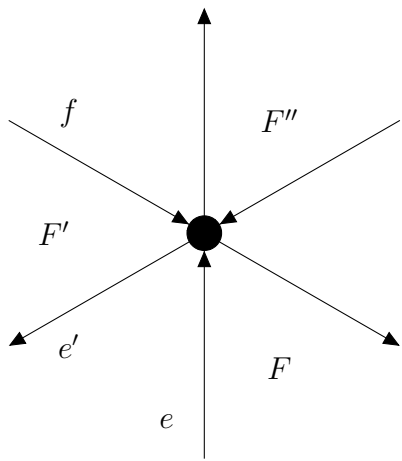

The condition $\gamma(G\mino{1}e)<\gamma(G)$ implies that $e$ is a
proper 1-semiloop, so already we know that every edge of $G$ is a loop that
does not enclose its own face.
Let $e$ be any edge and let $v$ be the vertex at which $e$ is a loop.
The condition $\gamma(G\mino{\omega^2}e)<\gamma(G)$
implies that $e$ is also a proper
$\omega$-semiloop.  Let $F$ be the face on the right side of $e$,
and let $F'$ be the face on the right side of the left
successor $e'$ of $e$
(see Figure \ref{fig:loops-and-faces-around-v}).
If faces $F$ and $F'$ were distinct, then $\omega$-reduction of $e$
would not reduce the genus and $e$ would not be an $\omega$-semiloop.
So $F=F'$.  Applying this same reasoning to the next edge (clockwise
from $e$) in the in-star at $v$ (denoted by $f$ in
Figure \ref{fig:loops-and-faces-around-v})
shows that the face $F'$ is, in turn,
identical to the second face beyond it (denoted by $F''$),
continuing to go clockwise
around $v$.  Continuing in this manner we find that every second ``face''
around $v$ is really just part of one single face.

In a similar manner, the condition $\gamma(G\mino{\omega}e)<\gamma(G)$
implies that $e$ is also a proper $\omega^2$-semiloop, and we find that all the
``faces'' at $v$ which were not accounted for in the previous paragraph
(being every \textit{other} second face around $v$) are, again, just one
single face (but necessarily distinct from the face $F$ we found there).
So the component of $G$ that consists of loops at $v$ is a posy.
So $G$ itself is just a
disjoint union of all these posies and has genus $k$, so we are done.
\eopf

\section{Tutte invariants}
\label{sec:tutte-invariants}

We now extend the notion of a Tutte invariant to alternating dimaps,
and investigate what invariants of this type exist.    \\

\noindent\textbf{Definition}

A \textit{simple Tutte invariant}
for alternating dimaps is a function $F$ defined
on every alternating dimap such that
$F$ is invariant under isomorphism, $F(\hbox{empty alternating dimap})=1$,
and there exist
$w, x, y, z$ such that, for any alternating dimap $G$,
\begin{enumerate}
\item
for any ultraloop $e$ of $G$,
\begin{equation}
\label{eq:tutte-invariant-ultraloop}
F(G) = w\,F(G-e)     % ***** ensure  G-e  defined early on
\end{equation}
\item
for any proper 1-loop $e$ of $G$,
\begin{equation}
\label{eq:tutte-invariant-1-loop}
F(G) = x\,F(G\mino{1}e)
\end{equation}
\item
for any proper $\omega$-loop $e$ of $G$,
\begin{equation}
\label{eq:tutte-invariant-omega-loop}
F(G) = y\,F(G\mino{\omega}e)
\end{equation}
\item
for any proper $\omega^2$-loop $e$ of $G$,
\begin{equation}
\label{eq:tutte-invariant-omega2-loop}
F(G) = z\,F(G\mino{\omega^2}e)
\end{equation}
\item
for any edge $e$ of $G$ that is not an ultraloop or a triloop,
\begin{equation}
\label{eq:tutte-invariant-other}
F(G) = F(G\mino{1}e) + F(G\mino{\omega}e) + F(G\mino{\omega^2}e) .
\end{equation}
\end{enumerate}

\begin{lemma}
\label{lemma:tutte-w-eq-0}
Let $F$ be a simple Tutte invariant of alternating dimaps.
If $w=0$ then $F(G)=0$ for any nonempty $G$.
\end{lemma}

\pf
Induction on $|E(G)|$.
\eopf

\begin{thm}
\label{thm:tutte-recipe1}
The only simple Tutte invariants of alternating dimaps are:
\begin{itemize}
\item[(a)] $F(G)=0$ for nonempty $G$, with $w=0$;
\item[(b)] $F(G)=3^{|E(G)|}$, with $w=x=y=z=3$;
\item[(c)] $F(G)=(-1)^{|V(G)|}$, with $y=z=1$ and $x=w=-1$;
\item[(d)] $F(G)=(-1)^{\cfaces(G)}$, with $x=z=1$ and $y=w=-1$;
\item[(e)] $F(G)=(-1)^{\afaces(G)}$, with $x=y=1$ and $z=w=-1$.
\end{itemize}
\end{thm}

\pf
Let $F$ be a simple Tutte invariant of alternating dimaps, with $w, x, y, z$ as
in the definition.

If $w=0$ then $F(G)=0$ for nonempty $G$, by Lemma \ref{lemma:tutte-w-eq-0}.
So suppose $w\not=0$.

For any $k\ge1$, let $U_k$ be a disjoint union of $k$ ultraloops.
For any $k\ge2$, let $L_{k,1}$ be the directed cycle of $k$ vertices
and $k$ edges, which has one a-face, one c-face and $k$ in-stars.
For $k\ge2$ and $\mu\in\{\omega,\omega^2\}$, let $L_{k,\mu}$ be the
alternating dimap consisting of a single vertex with $k$ $\mu$-loops.
% **** do I want this defn earlier in a defns and notn section? *****

It is clear that
\begin{eqnarray}
F(U_k)  & = &  w^k,   \\
F(L_{2,1})  & = &  xw,   \\
F(L_{2,\omega})  & = &  yw,   \\
F(L_{2,\omega^2})  & = &  zw .
\end{eqnarray}

Consider the alternating dimap $L_{2,1}+e_{\omega^2}$
obtained by adding an $\omega^2$-loop $e$ to
a vertex $v$ of $L_{2,1}$.  Let $f$ (respectively, $g$) be the edge
of $L_{2,1}$ going out of (resp., into) $v$.  Observe that
$f$ is a 1-loop and $g$ is not a triloop.
We can calculate $F(L_{2,1}+e_{\omega^2})$ by
applying (\ref{eq:tutte-invariant-1-loop}) at $f$
(or (\ref{eq:tutte-invariant-omega2-loop}) at $e$), obtaining $xzw$.
Alternatively, we can apply
(\ref{eq:tutte-invariant-other}) at $g$, obtaining $(x+z+w)w$.
Equating the results
and using $w\not=0$,
we obtain
\begin{equation}
\label{eq:xzw-xz}
x + z + w = xz .
\end{equation}
Similar reasoning for the trials $(L_{2,1}+e_{\omega^2})^{\omega}$
and $(L_{2,1}+e_{\omega^2})^{\omega^2}$ gives
\begin{eqnarray}
F((L_{2,1}+e_{\omega^2})^{\omega})  & = &  xyw,   \label{eq:tutte-xyw}  \\
F((L_{2,1}+e_{\omega^2})^{\omega^2})  & = &  yzw ,   \label{eq:tutte-yzw}   \\
x + y + w  & = & xy ,
\label{eq:xyw-xy}   \\
y + z + w  & = & yz .
\label{eq:yzw-yz}
\end{eqnarray}

Now consider the alternating dimap obtained from $L_{2,1}$ (with edges
$g,h$) by adding, to the endpoint $v$ of $h$, a clockwise loop $e$ within
the anticlockwise face and an anticlockwise loop $f$ within the clockwise face.
Call it $A$.  Using the $\omega^2$-loop $e$, or the $\omega$-loop $f$,
or the 1-loop $g$, we find that
\begin{equation}
\label{eq:tutte-xyzw}
F(A)=xyzw.
\end{equation}
But using $h$, which is not
a triloop, we have
\begin{eqnarray*}
F(A)  & = &
F(L_{2,1}+e_{\omega^2}) + F((L_{2,1}+e_{\omega^2})^{\omega})
+ F((L_{2,1}+e_{\omega^2})^{\omega^2})   \\
& = &  xzw + xyw + yzw   \\
& = &  (xz + xy + yz)w .
\end{eqnarray*}
Equating with (\ref{eq:tutte-xyzw}), and using $w\not=0$, we obtain
\begin{equation}
\label{eq:xz-xy-yz-xyz}
xz + xy + yz = xyz .
\end{equation}

From (\ref{eq:xzw-xz}), (\ref{eq:xyw-xy}), (\ref{eq:yzw-yz}) we obtain
\[
w = xz - x - z = xy - x - y = yz - y - z .
\]
The second equality here gives $(x-1)z = (x-1)y$, so either $x=1$ or $y=z$.
Similarly, either $y=1$ or $x=z$, and either $z=1$ or $x=y$.
Combining these, we have one of
\begin{enumerate}
\item[(i)]  $x=y=z$,
\item[(ii)]  $x=y=1$ and $z\not=1$,
\item[(iii)]  $x=z=1$ and $y\not=1$,
\item[(iv)]  $y=z=1$ and $x\not=1$.
\end{enumerate}
If (i)    % **** make the numbering styles here, and in the list, match ****
holds, then any of
(\ref{eq:xzw-xz}), (\ref{eq:xyw-xy}), (\ref{eq:yzw-yz}) gives
\begin{equation}
\label{eq:w-x-x2}
w=x(x-2) .
\end{equation}
Also (\ref{eq:xz-xy-yz-xyz}) gives $3x^2=x^3$, whence $x=3$
(since $x=0$ would imply $w=0$, by (\ref{eq:w-x-x2})) and $w=3$ (by (\ref{eq:w-x-x2})).

If (ii) holds, then (\ref{eq:xyw-xy}) gives $w=-1$.
Similarly, cases
(iii) and (iv) give $w=-1$ too.

Also, (\ref{eq:xz-xy-yz-xyz}) implies $z=-1$ in case (ii), $y=-1$ in case (iii)
and $z=-1$ in case (iv).

We now establish the form of $F$ for each of cases (i)--(iv) in turn.
The numbering of the claims indicates the case to which each applies.

Claim (i):
$F(G) = 3^{|E(G)|}$.

Proof of Claim (i):  we use induction on $|E(G)|$.
If $|E(G)| = 0$ then the claim is true by the definition of $F$.
Suppose $|E(G)|=m>1$.  Let $e\in E(G)$.  If $e$ is an ultraloop,
then $F(G)=wF(G-e)=3F(G-e)=3\cdot3^{m-1}$, by the inductive hypothesis,
which equals $3^m$.  If $e$ is a proper $\mu$-loop, with $\mu\in\{1,\omega,\omega^2\}$,
then $F(G)=xF(G\mino{\mu}e)=3\cdot3^{m-1}=3^m$.
If $e$ is neither an ultraloop or a triloop, then
$F(G)=F(G\mino{1}e)+F(G\mino{\omega}e)+F(G\mino{\omega^2}e)=3^{m-1}+3^{m-1}+3^{m-1}=3^m$.

Claim (ii):
$F(G) = (-1)^{\afaces(G)}$.

Proof of Claim (ii):  we use induction on $|E(G)|$.
If $|E(G)| = 0$ then the claim is true by the definition of $F$.
Suppose $|E(G)|=m>1$.  Let $e\in E(G)$.  If $e$ is an ultraloop,
then $F(G)=wF(G-e)=-F(G-e)=-(-1)^{\afaces(G)-1}=(-1)^{\afaces(G)}$.
Observe that the number of anticlockwise faces in an alternating dimap
is unchanged by 1- or $\omega$-reduction, but may be altered by
$\omega^2$-reduction.
If $e$ is a proper $\omega$-loop,
then $F(G)=zF(G\mino{\omega^2}e)=-F(G\mino{\omega^2}e)=-(-1)^{\afaces(G\mino{\omega^2}e)}=-(-1)^{\afaces(G)-1}=(-1)^{\afaces(G)}$.
If $e$ is a proper $\mu$-loop with $\mu\in\{1,\omega^2\}$, then
$F(G)=F(G\mino{\mu}e)=(-1)^{\afaces(G\mino{\mu}e)}
=(-1)^{\afaces(G)}$.
If $e$ is a proper $\omega$-semiloop, then $\hbox{af}(G\mino{\omega^2}e)=\hbox{af}(G)+1$,
while if $e$ is neither an ultraloop nor a triloop nor an $\omega$-semiloop, then
$\hbox{af}(G\mino{\omega^2}e)=\hbox{af}(G)-1$.  In any event, if $e$ is not a triloop, then
$F(G)=F(G\mino{1}e)+F(G\mino{\omega}e)+F(G\mino{\omega^2}e)
=(-1)^{\afaces(G\mino{1}e)}
+(-1)^{\afaces(G\mino{\omega}e)}
+(-1)^{\afaces(G\mino{\omega^2}e)}
=(-1)^{\afaces(G)}
+(-1)^{\afaces(G)}
+(-1)^{\afaces(G)\pm1}
=(-1)^{\afaces(G)}$.

Claim (iii):
$F(G) = (-1)^{\cfaces(G)}$.

Claim (iv):
$F(G) = (-1)^{|V(G)|}$.

The proofs of Claims (iii) and (iv) are similar to that of Claim (ii), and
are left as an exercise.  For Claim (iv), bear in mind that $|V(G)|$ is the
number of in-stars of $G$.
\eopf   \\

A \textit{Tutte invariant} for alternating dimaps is defined as for
a simple Tutte invariant, except that condition 5, with
(\ref{eq:tutte-invariant-other}), is replaced by a requirement that
there exist nonzero $a,b,c$ such that, for any alternating dimap $G$
and any $e\in E(G)$,
\[
F(G) = aF(G\mino{1}e)+bF(G\mino{\omega}e)+cF(G\mino{\omega^2}e) ,
\]

\begin{thm}
\label{thm:tutte-recipe2}
The only Tutte invariants of alternating dimaps are:
\begin{itemize}
\item[(a)] $F(G)=0$ for nonempty $G$, with $w=x=y=z=0$
\item[(b)]  $F(G)=3^{|E(G)|}a^{|V(G)|}b^{\cfaces(G)}c^{\afaces(G)}$,
$\frac{x}{a}=\frac{y}{b}=\frac{z}{c}=3$, $\frac{w}{abc}=3$.
\item[(c)]  $F(G)=
a^{|V(G)|}b^{\cfaces(G)}(-c)^{\afaces(G)}$,
$\frac{x}{a}=\frac{y}{b}=1$, $\frac{z}{c}=-1$, $\frac{w}{abc}=-1$.
\item[(d)]  $F(G)=
a^{|V(G)|}(-b)^{\cfaces(G)}c^{\afaces(G)}$,
$\frac{x}{a}=\frac{z}{c}=1$, $\frac{y}{b}=-1$, $\frac{w}{abc}=-1$.
\item[(e)]  $F(G)=
(-a)^{|V(G)|}b^{\cfaces(G)}c^{\afaces(G)}$,
$\frac{y}{b}=\frac{z}{c}=1$, $\frac{x}{a}=-1$, $\frac{w}{abc}=-1$.
\end{itemize}
\end{thm}

\pf
If
\[
F(G) = aF(G\mino{1}e)+bF(G\mino{\omega}e)+cF(G\mino{\omega^2}e) ,
\]
then define
\begin{equation}
\label{eq:F-dash}
F'(G) = \frac{F(G)}{a^{|V(G)|}b^{\cfaces(G)}c^{\afaces(G)}} .
\end{equation}

If $e$ is an ultraloop,
then
\begin{eqnarray*}
F'(G) & = & 
\frac{wF(G-e)}{a^{|V(G)|}b^{\cfaces(G)}c^{\afaces(G)}}   \\
& = &
\frac{wF(G-e)}{a^{|V(G-e)|+1}b^{\cfaces(G-e)+1}c^{\afaces(G)+1}}   \\
& = &
\frac{w}{abc} \cdot
\frac{F(G-e)}{a^{|V(G-e)|}b^{\cfaces(G-e)}c^{\afaces(G-e)}}   \\
& = &
\frac{w}{abc} \cdot F'(G-e) .
\end{eqnarray*}

If $e$ is a 1-loop, then
\begin{eqnarray*}
F'(G) & = & 
\frac{xF(G\mino{1}e)}{a^{|V(G)|}b^{\cfaces(G)}c^{\afaces(G)}}   \\
& = &
\frac{xF(G\mino{1}e)}{a^{|V(G\mino{1}e)|+1}b^{\cfaces(G\mino{1}e)}c^{\afaces(G\mino{1}e)}}   \\
& = &
\frac{x}{a} \cdot
\frac{F(G\mino{1}e)}{a^{|V(G\mino{1}e)|}b^{\cfaces(G\mino{1}e)}c^{\afaces(G\mino{1}e)}}   \\
& = &
\frac{x}{a} \cdot F'(G\mino{1}e) .
\end{eqnarray*}

Similarly, if $e$ is an $\omega$-loop, then
\[
F'(G) =
\frac{y}{b} \cdot F'(G\mino{\omega}e) ,
\]
and if $e$ is an $\omega^2$-loop, then
\[
F'(G) =
\frac{z}{c} \cdot F'(G\mino{\omega^2}e) .
\]

If $e\in E(G)$ is neither an ultraloop nor a triloop,
\begin{eqnarray*}
F'(G) & = &
\frac{aF(G\mino{1}e)+bF(G\mino{\omega}e)+cF(G\mino{\omega^2}e)}{a^{|V(G)|}b^{\cfaces(G)}c^{\afaces(G)}}   \\
& = &
\frac{F(G\mino{1}e)}{a^{|V(G)|-1}b^{\cfaces(G)}c^{\afaces(G)}}
+ \frac{F(G\mino{\omega}e)}{a^{|V(G)|}b^{\cfaces(G)-1}c^{\afaces(G)}} +  \\
&&
\frac{F(G\mino{\omega^2}e)}{a^{|V(G)|}b^{\cfaces(G)}c^{\afaces(G)-1}}   \\
& = &
\frac{F(G\mino{1}e)}{a^{|V(G\mino{1}e)|}b^{\cfaces(G\mino{1}e)}c^{\afaces(G\mino{1}e)}}
+ \frac{F(G\mino{\omega}e)}{a^{|V(G\mino{\omega}e)|}b^{\cfaces(G\mino{\omega}e)}c^{\afaces(G\mino{\omega}e)}} +   \\
&&
\frac{F(G\mino{\omega^2}e)}{a^{|V(G\mino{\omega^2}e)|}b^{\cfaces(G\mino{\omega^2}e)}c^{\afaces(G\mino{\omega^2}e)}}   \\
& = &
F'(G\mino{1}e) + F'(G\mino{\omega}e) + F'(G\mino{\omega^2}e) .
\end{eqnarray*}

Hence $F'$ is a simple Tutte invariant.

By Theorem \ref{thm:tutte-recipe1}, we must have one of
\begin{itemize}
\item $F'(G)=0$, $w=x=y=z=0$
\item $F'(G)=3^{|E(G)|}$,
$\frac{x}{a}=\frac{y}{b}=\frac{z}{c}=3$, $\frac{w}{abc}=3$.
\item $F'(G)=(-1)^{\afaces(G)}$,
$\frac{x}{a}=\frac{y}{b}=1$, $\frac{z}{c}=-1$, $\frac{w}{abc}=-1$.
\item $F'(G)=(-1)^{\cfaces(G)}$,
$\frac{x}{a}=\frac{z}{c}=1$, $\frac{y}{b}=-1$, $\frac{w}{abc}=-1$.
\item $F'(G)=(-1)^{|V(G)|}$,
$\frac{y}{b}=\frac{z}{c}=1$, $\frac{x}{a}=-1$, $\frac{w}{abc}=-1$.
\end{itemize}
Use (\ref{eq:F-dash}) to complete the proof.
\eopf   \\

%  The following just here as a convenient note, can be deleted later.  ***
%   F'(G) =
%   \frac{F(G)}{a^{|V(G)|}b^{\cfaces(G)}c^{\afaces(G)}} .
%  **** Next few lines need changing to make in terms of F somehow.
%  \begin{itemize}
%  \item $F(G)=0$, $w=x=y=z=0$
%  \item $F(G)=3^{|E(G)|}a^{|V(G)|}b^{\cfaces(G)}c^{\afaces(G)}$,
%  $\frac{x}{a}=\frac{y}{b}=\frac{z}{c}=3$, $\frac{w}{abc}=3$.
%  \item $F(G)=
%  a^{|V(G)|}b^{\cfaces(G)}(-c)^{\afaces(G)}$,
%  $\frac{x}{a}=\frac{y}{b}=1$, $\frac{z}{c}=-1$, $\frac{w}{abc}=-1$.
%  \item $F(G)=
%  a^{|V(G)|}(-b)^{\cfaces(G)}c^{\afaces(G)}$,
%  $\frac{x}{a}=\frac{z}{c}=1$, $\frac{y}{b}=-1$, $\frac{w}{abc}=-1$.
%  \item $F(G)=
%  (-a)^{|V(G)|}b^{\cfaces(G)}c^{\afaces(G)}$,
%  $\frac{y}{b}=\frac{z}{c}=1$, $\frac{x}{a}=-1$, $\frac{w}{abc}=-1$.
%  \end{itemize}

Other definitions of Tutte invariants for alternating dimaps are possible.   \\

\noindent\textbf{Definition}

An \textit{extended Tutte invariant}
for alternating dimaps is a function $F$ defined
on every alternating dimap such that
$F$ is invariant under isomorphism, $F(\hbox{empty alternating dimap})=1$,
and there exist
$w, x, y, z, a, b, c, d, e, f, g, h, i, j, k, l$
such that, for any alternating dimap $G$,
\begin{enumerate}
\item
for any ultraloop $e$ of $G$,
\begin{equation}
F(G) = w\,F(G-e)     % ***** ensure  G-e  defined early on
\end{equation}
\item
for any proper 1-loop $e$ of $G$,
\begin{equation}
F(G) = x\,F(G\mino{1}e)
\end{equation}
\item
for any proper $\omega$-loop $e$ of $G$,
\begin{equation}
F(G) = y\,F(G\mino{\omega}e)
\end{equation}
\item
for any proper $\omega^2$-loop $e$ of $G$,
\begin{equation}
F(G) = z\,F(G\mino{\omega^2}e)
\end{equation}
\item
for any proper 1-semiloop $e$ of $G$,
\begin{equation}
F(G) = aF(G\mino{1}e) + bF(G\mino{\omega}e) + cF(G\mino{\omega^2}e) .
\end{equation}
\item
for any proper $\omega$-semiloop $e$ of $G$,
\begin{equation}
F(G) = dF(G\mino{1}e) + eF(G\mino{\omega}e) + fF(G\mino{\omega^2}e) .
\end{equation}
\item
for any proper $\omega^2$-semiloop $e$ of $G$,
\begin{equation}
F(G) = gF(G\mino{1}e) + hF(G\mino{\omega}e) + iF(G\mino{\omega^2}e) .
\end{equation}
\item
for any edge $e$ of $G$ that is not an ultraloop or a triloop
or a semiloop,
\begin{equation}
F(G) = jF(G\mino{1}e) + kF(G\mino{\omega}e) + lF(G\mino{\omega^2}e) .
\end{equation}
\end{enumerate}

\noindent\textbf{Problem}

Characterise all extended Tutte invariants of alternating dimaps.   \\

One basic extended Tutte invariant is
\[
F(G) = \alpha^{|E(G)|}\beta^{|V(G)|}\gamma^{\afaces(G)}\delta^{\cfaces(G)} ,
\]
with $\alpha,\beta,\gamma,\delta\not=0$.
This satisfies the definition with $w=\alpha\beta\gamma\delta$,
$x=\alpha\beta$, $y=\alpha\gamma$, $z=\alpha\delta$, $a=\alpha/\beta$, $f=\alpha/\delta$,
$h=\alpha/\gamma$, $b=c=d=e=g=i=0$, $j=\alpha\beta/3$, $k=\alpha\gamma/3$, and
$l=\alpha\delta/3$.

Extended Tutte invariants are much richer than Tutte invariants, since they
include the Tutte polynomial for planar graphs,
in a sense we now explain.

The \textit{Tutte polynomial} $T(G;x,y)$ of a graph $G$ has the following inductive definition.
If $E(G)=\emptyset$ then $T(G;x,y)=1$.  Otherwise, for any $e\in E(G)$,
\begin{eqnarray*}
T(G;x,y) = \left\{\begin{array}{ll}
x \, T(G\setminus e;x,y), & \hbox{if $e$ is a coloop} ;   \\
y \,T(G/e;x,y), & \hbox{if $e$ is a loop} ;   \\
T(G/e;x,y) + T(G\setminus e;x,y), & \hbox{otherwise}.
\end{array}\right.
\end{eqnarray*}

To any orientably 2-cell-embedded (undirected) graph $G$,
we can associate two alternating dimaps
$\altc(G)$ and $\alta(G)$ as follows.
For $\altc(G)$ (respectively, $\alta(G)$),
replace each edge $e=uv\in E(G)$ by a pair of oppositely directed
edges $(u,v)$ and $(v,u)$, forming a clockwise (resp., anticlockwise) face of size two.
The faces of $G$ now all correspond to anticlockwise (resp., clockwise) faces
in $\altc(G)$ (resp., $\alta(G)$).

%  We write $G_{\hbox{\scriptsize c},i}$ for any alternating dimap which can be obtained from some
%  $\altc(G)$ by adding $i$ $\omega^2$-loops to vertices so that they sit in a-faces,
%  and $G_{\hbox{\scriptsize a},i}$ is defined similarly,
%  using $\alta(G)$, $\omega$-loops, and c-faces.

For any alternating dimap $G$, define $T_c(G;x,y)$ and $T_a(G;x,y)$ as follows.
If $E(G)=\emptyset$, then $T_c(G;x,y)=T_a(G;x,y)=1$.  Otherwise, for any $e\in E(G)$,
\begin{eqnarray*}
T_c(G;x,y) & = & \left\{\begin{array}{lll}
T_c(G\mino{*}e;x,y), &
\multicolumn{2}{l}{\hbox{if $e$ is an $\omega^2$-loop (including an ultraloop);}}   \\
x \, T_c(G\mino{\omega^2}e;x,y), ~~~~~~ &
\multicolumn{2}{l}{\hbox{if $e$ is an $\omega$-semiloop;}}    \\
y \, T_c(G\mino{1}e;x,y), &
\multicolumn{2}{l}{\hbox{if $e$ is a proper 1-semiloop or an $\omega$-loop;}}    \\
\multicolumn{2}{l}{T_c(G\mino{1}e;x,y) + T_c(G\mino{\omega^2}e;x,y),} &
\hbox{if $e$ is not a semiloop.}
\end{array}\right.   \\
&&   \\
T_a(G;x,y) & = & \left\{\begin{array}{lll}
T_a(G\mino{*}e;x,y), &
\multicolumn{2}{l}{\hbox{if $e$ is an $\omega$-loop (including an ultraloop);}}    \\
x \, T_a(G\mino{\omega}e;x,y), ~~~~~~ &
\multicolumn{2}{l}{\hbox{if $e$ is an $\omega^2$-semiloop;}}    \\
y \, T_a(G\mino{1}e;x,y), &
\multicolumn{2}{l}{\hbox{if $e$ is a proper 1-semiloop or an $\omega^2$-loop;}}    \\
\multicolumn{2}{l}{T_a(G\mino{1}e;x,y) + T_a(G\mino{\omega}e;x,y),}  &
\hbox{if $e$ is not a semiloop.}
\end{array}\right.
\end{eqnarray*}

\begin{thm}
For any plane graph $G$,
\[
T(G;x,y) = T_c(\altc(G);x,y) = T_a(\alta(G);x,y).
\]
\end{thm}

\pf
For any vertex $v$,
write $L^{(\omega)}(v)$ and $L^{(\omega^2)}(v)$ for an $\omega$-loop and an $\omega^2$-loop,
respectively, at $v$.  If such a loop is added to an alternating dimap, it must be
placed within a c-face or an a-face, respectively.

Consider $\altc(G)$.  Observe that, for any $uv\in E(G)$,
\begin{eqnarray*}
\altc(G)\mino{1}(u,v)  & = &  \altc(G/uv) + L^{(\omega^2)}(u') ,   \\
\altc(G)\mino{\omega^2}(u,v)  & = &  \altc(G\setminus uv) + L^{(\omega^2)}(u') ,
\end{eqnarray*}
for some $u'$.  (Mostly $u'=u$, except that a little more detail is needed if $(u,v)$ is a
proper 1-semiloop, but the exact location of these extra triloops is not important.)

These observations can be used to prove, by induction on $|E(G)|$, that
$T(G;x,y) = T_c(\altc(G);x,y)$ for any plane graph $G$.

It is clear from the definitions that they are identical when $G$ is empty.

Suppose then that $T(G;x,y) = T_c(\altc(G);x,y)$ when $|E(G)|<m$, where $m\ge1$.

Let $G$ be any orientably 2-cell-embedded graph on $m$ edges, and let $e=uv\in E(G)$.

If $e$ is a coloop, then $(u,v)$ and $(v,u)$ are both $\omega$-semiloops in $\altc(G)$.
(Conversely, if $(u,v)$ and $(v,u)$ are both \textit{proper}
$\omega$-semiloops in $\altc(G)$, then $uv$ is a coloop in $G$.  This does not hold in
general if $G$ is not plane, however.)  We have
\begin{eqnarray*}
T_c(\altc(G);x,y)  & = &  x \, T_c(\altc(G)\mino{\omega^2}(u,v);x,y)
= x \, T_c(\altc(G/uv)+ L^{(\omega^2)}(u');x,y)    \\
& = &  x \, T_c(\altc(G/uv);x,y) = x \, T(G/uv;x,y) = T(G;x,y) ,
\end{eqnarray*}
where the penultimate equality uses the inductive hypothesis.

if $e$ is a loop, then in $\altc(G)$ the directed versions $(u,v)$ and $(v,u)$
are both 1-semiloops.  (This time, the converse holds even if $G$ is not plane.)
One of them may also be
an $\omega$-loop, but neither is an $\omega^2$-loop.  In any case, we find that
$T_c(\altc(G);x,y)=T(G;x,y)$ by a similar argument to that just used for coloops.

If $e$ is neither a coloop nor a loop, we have
\begin{eqnarray*}
T_c(\altc(G);x,y)  & = &  T_c(\altc(G)\mino{1}(u,v);x,y) + T_c(\altc(G)\mino{\omega^2}(u,v);x,y)   \\
& = &  T_c(\altc(G/uv)+ L^{(\omega^2)}(u');x,y) +
T_c(\altc(G\setminus uv)+ L^{(\omega^2)}(u');x,y)    \\
& = &  T_c(\altc(G/uv);x,y) + T_c(\altc(G\setminus uv);x,y)    \\
& = &  T(G/uv;x,y) + T(G\setminus uv;x,y)   ~~~~~ \hbox{(by the inductive hypothesis)}    \\
& = &  T(G;x,y) .
\end{eqnarray*}

We conclude by induction that $T(G;x,y) = T_c(\altc(G);x,y)$ holds for all $G$.

The proof that $T(G;x,y) = T_a(\alta(G);x,y)$ follows the same line, with appropriate adjustments.
\eopf   \\

Having constructed alternating dimaps from embedded graphs, by replacing edges by c-faces,
or by a-faces, of size 2, it is natural to ask about replacing edges by in-stars of size 2.
To do this, for an embedded
graph $G$, first construct its medial graph, $\med(G)$,
then turn it into an alternating dimap by directing the edges so as to ensure the
alternating property.  For each component of $G$, there are two such ways of directing the
edges in that component.  There are therefore $2^{k(G)}$ different alternating dimaps
constructible from $G$ in this way, all with $\med(G)$ as the underlying embedded graph.
We refer to any one of them as $\alti(G)$.

Write $T_i(G;x)$ for any invariant of alternating dimaps that satisfies the following.
\[
T_i(G;x)  =  \left\{\begin{array}{lll}
1,  &  \hbox{if $G$ is empty,}   \\
T_i(G\mino{*}e;x), &
\multicolumn{2}{l}{\hbox{if $e$ is a 1-loop (including an ultraloop);}}   \\
x \, T_i(G\mino{\omega^2}e;x), ~~~~~~ &
\multicolumn{2}{l}{\hbox{if $e$ is a proper $\omega$-semiloop or an $\omega^2$-loop;}}    \\
x \, T_i(G\mino{\omega}e;x), &
\multicolumn{2}{l}{\hbox{if $e$ is a proper $\omega^2$-semiloop or an $\omega$-loop;}}    \\
\multicolumn{2}{l}{T_i(G\mino{\omega}e;x) + T_i(G\mino{\omega^2}e;x),} &
\hbox{if $e$ is not a semiloop.}
\end{array}\right. 
\]
This is not a full definition of a unique $T_i(G;x)$, since we have not specified what happens
if $e$ is a proper 1-semiloop.  But it will turn out that, in the alternating dimaps of interest here,
proper 1-semiloops do not arise.

\begin{thm}
For any plane graph $G$,
\[
T(G;x,x) = T_i(\alti(G);x) .
\]
\end{thm}

\pf
For any alternating dimap $H$, write $H^{(k)}$ for any alternating dimap obtained from $H$
by performing, $k$ times, a subdivision of an edge (by insertion in it of a new vertex of
indegree = outdegree = 1, with the edge going into the new vertex being a proper 1-loop).

Let $G$ be a plane graph and fix any specific $\alti(G)$.
If $e\in E(G)$, write $e^{\downarrow}$ for either of the edges
of $\alti(G)$ that are directed into the vertex representing $e$.

Observe that, for any $e\in E(G)$,
\[
\{ \alti(G)\mino{\omega}e^{\downarrow}, \alti(G)\mino{\omega^2}e^{\downarrow} \}
=  \{ \alti(G/e)^{(1)}, \alti(G\setminus e)^{(1)} \} .
\]
Therefore
\[
T_i(\alti(G)\mino{\omega}e^{\downarrow}; x) + T_i(\alti(G)\mino{\omega^2}e^{\downarrow};x)
=  T_i(\alti(G/e)^{(1)}; x) + T_i(\alti(G\setminus e)^{(1)};x) .
\]
If $e$ is either a coloop or a loop, then
%**** edit next line to suit ****
\begin{eqnarray*}
\alti(G)\mino{\omega^2}e^{\downarrow}
& \in &  \{ \alti(G/e)^{(1)}, \alti(G\setminus e)^{(1)} \} ,   \\
\alti(G)\mino{\omega}e^{\downarrow}
& \in &  \{ \alti(G/e)^{(1)}, \alti(G\setminus e)^{(1)} \} .
\end{eqnarray*}

We now prove the theorem by induction on $|E(G)|$.  The base case is immediate from the definition.  So suppose $G$ is an orientably 2-cell-embedded graph
with $m$ edges, where $m\ge1$.

If $e$ is a coloop or a loop, then $e^{\downarrow}$ is an $\omega$- or an $\omega^2$-semiloop
in $\alti(G)$,
except that it is not a 1-loop.  From our (partial) definition of $T_i(\alti(G);x)$, we have
\begin{eqnarray*}
T_i(\alti(G);x)  & \in &
\{ x T_i(\alti(G)\mino{\omega}e^{\downarrow};x), \,
x T_i(\alti(G)\mino{\omega^2}e^{\downarrow};x) \}    \\
& = & \{  x T_i(\alti(G/e)^{(1)};x), \,
x T_i(\alti(G\setminus e)^{(1)};x)  \}   \\
& = & \{  x T_i(\alti(G/e);x), \,
x T_i(\alti(G\setminus e);x)  \}   \\
& = & \{  x T(G/e;x,x), \,
x T(G\setminus e;x,x)  \},
\end{eqnarray*}
by the inductive hypothesis.  But, for such an $e$, the graphs $G/e$ and $G\setminus e$ have
isomorphic cycle matroids, so their Tutte polynomials are identical.  Therefore
\[
T_i(\alti(G);x)  =  x \, T(G/e;x,x)  =  x \, T(G\setminus e;x,x) .
\]
But these two quantities each equal $T(G;x,x)$, for such an $e$, so we are done in this case.

If $e$ is neither a loop nor a coloop, then 
\begin{eqnarray*}
T_i(\alti(G);x)
& = &  T_i(\alti(G)\mino{\omega}e^{\downarrow};x) +
T_i(\alti(G)\mino{\omega^2}e^{\downarrow};x)   \\
& = &  T_i(\alti(G/e)^{(1)}; x) + T_i(\alti(G\setminus e)^{(1)};x)   \\
& = &  T_i(\alti(G/e); x) + T_i(\alti(G\setminus e);x)   \\
& = &  T(G/e; x, x) + T(G\setminus e; x, x)   ~~~~~~ \hbox{(by the inductive hypothesis)}   \\
& = &  T(G;x,x) .
\end{eqnarray*}

The result follows.

Observe that, since $\med(G)$ is 4-regular, $\alti(G)$ has no proper 1-semiloops.
Furthermore, the only minors of it we need to form do not require 1-reduction, so these
minors are each 4-regular and so have no proper 1-semiloop too.
\eopf   \\

The Tutte polynomial evaluation $T(G;x,x)$ is just the Martin polynomial of $\med(G)$
(see, e.g., \cite{ellismonaghan04}).

One reason that Tutte invariants of alternating dimaps are more limited
than Tutte invariants for graphs is the non-commutativity of the minor
operations.  The definitions of such invariants for alternating dimaps
require the stated recursive relations to hold for reduction of \textit{any}
edge of the stated type,
which means that the invariant will need to be unperturbable by
some variations of the order of operations.

These observations raise the possibility that better invariants may
come from including an ordering of the edges in the object to which
the invariant applies.   \\

\noindent\textbf{Definitions}

An \textit{ordered alternating dimap}
is a pair $(G,<)$ where $G$ is an alternating dimap and $<$ is a
linear order on $E(G)$.

If $(G,<)$ is an ordered alternating dimap and
$\mu\in\{1,\omega,\omega^2\}$, then the
$\mu$\textit{-reduction} $(G,<)\mino{\mu}$ of $(G,<)$
is the ordered alternating dimap $(G\mino{\mu}e_0,<')$ where $e_0$
is the first edge in $E(G)$ under $<$ and the
order $<'$ on $E(G)\setminus\{e_0\}$ is obtained by simply removing
$e_0$ from the order $<$.

Tutte invariants and extended Tutte invariants
are defined for these objects by modifying the definitions of such invariants
for ordinary alternating dimaps as follows:
\begin{enumerate}
\item The definitions apply to ordered alternating dimaps, rather than
just to alternating dimaps.
\item All references to $G\mino{\mu}e$ are replaced by $(G,<)\mino{\mu}$,
for each $\mu$.
\item All universal quantification over edges is deleted (since there is
no choice of which edge to reduce, since it is always the first edge
in the ordering which must be reduced).
\item All reference to an edge $e$ is replaced by reference to the first
edge $e_0$ in the ordering.
\end{enumerate}
For example, the second condition in each of the definitions becomes:
if $e_0$ is a 1-loop, then $F((G,<)) = x F((G,<)\mino{1})$.

When $G$ is a general plane alternating dimap, the extended Tutte invariants
$T_c(G;x,y)$ and $T_a(G;x,y)$ we considered earlier actually
depend on the order in which the edges are considered.  So they pertain to ordered
alternating dimaps.  But, if $G$ has the form $\altc(H)$ (with analogous remarks applying
to $\alta(H)$), then the
order in which the edges $uv$ of $H$ are considered does not matter, and each time
a corresponding $(u,v)$ is reduced in $G$, it leaves behind an $\omega^2$-loop which can
be reduced at any time.  (Note also that, if we do not use $\omega$-reductions, we cannot
encounter those situations of non-commutativity for two edges that are consecutive in an a-face
or an in-star.)  For such cases, the invariants are well-defined without having to
specify an order on the edges at the beginning.   \\

\noindent\textbf{Problem}

Characterise (a) Tutte invariants, and (b) extended Tutte invariants,
of ordered alternating dimaps.   \\

\vspace{0.3in}

\noindent\textbf{Acknowledgements}

I am most grateful to Tony Grubman, Lorenzo Traldi and Ian Wanless for their comments.

\end{document}